\documentclass{amsart}
\usepackage{amssymb,latexsym}
\theoremstyle{definition}

\begin{document}

\title[Aut($G$) for $|G|=64$] {On the Groups and Automorphism Groups\\ of the 
Groups of order 64$p$\\ without a Normal Sylow $p$-Subgroup}

\author{Walter Becker}
\address{266 Brian Drive, Warwick, Rhode Island 02886}
\email{w\_becker@hotmail.com}
\author{Elaine W. Becker}
\address{266 Brian Drive, Warwick, Rhode Island 02886}
\email{elainewbecker@gmail.com}

\maketitle
\begin{abstract}
The groups of order 64$p$ without a normal sylow $p$-subgroup are listed, and their automorphism groups are also determined. As a by-product of our original effort to get these groups, we needed to determine the automorphism groups of those groups of order 64 with an odd-order automorphism. In view of the fact that we already had determined these groups and that these automorphism groups are not given explicitly in the literature, we have appended to this report these automorphism groups. In another project we were looking for new complete groups by following automorphism group towers up to completion when the computer memory allowed such followups. We did this for these groups of order 64. In another appendix we give the results of this work as applied to the groups of order 64. 
\end{abstract}
\section{Introduction}

The calculations reported here  were done in the mid 1990s in response to a 
then recent bibliographic survey by E. A. O'Brien and M. W. Short \cite{1}, 
in which it was stated that the number of groups of orders 192, 240, 
and 252 were the only orders below 256 for which the number of groups was 
still unknown. The groups of order 252, or more generally 36$p$, appear in 
the thesis of B. Malmrot \cite{2} as a prelude to his determination of the 
groups of order 72$p$, for $p > 3$. In an effort to fill in these gaps we 
have looked at the groups and automorphism groups of the groups of these 
orders. The work on the groups of orders 240 and 36$p$ ($p > 3$) will be 
given in later reports.\\

The number of groups of order 192 (or more generally 64$p$, $p > 2$) is very 
large, and from subsequent work of Bettina Eick and Hans Ulrich Besche 
\cite{10} (done around the same time as the calculations reported below were 
done) we have the results reported below in Tables 1a and 1b for the groups of order 
$64p$ for various choices of the odd prime $p$. For comparison with other 
orders we give similiar results for the groups of orders $16p$ and $32p$ 
from the earlier work of A. C. Lunn and J. K. Senior from the 1930s 
\cite{8}. Also included in Tables 1a and 1b are other results taken from the work of 
Besche and Eick reported in \cite{9}.\\

One way to verify the numbers in the following Table 1a for the $m=$ 3, 4, 5 
and 6 cases with a normal sylow $p$-subgroup would be by just counting the 
number of different ways the quotient groups of $C_2$, $C_4$, \ldots , $C_{64}$ 
arise in the subgroup lattice diagrams for the groups of orders 16, 32 and 64 
found in \cite{3}. The $C_2$ images (or equivalently the groups arising in the 
case of $p\equiv 1$ mod(2)) are called dimidiations in the older 
literature. If one just counts the number of normal subgroups of order 
$2^{m-1}$ in the Hall-Senior lattice subgroup diagrams for the groups of 
order $2^m$, one should have the number of groups of order $2^m*p$ arising 
for each group of order $2^m$ for the case when $p\equiv 1$ mod(2) but 
$p\not\equiv 1$ mod(4),\ldots . This is the way in which we initially had an idea of just how 
many groups of order $64p$ one would encounter in this enumeration process. 
This report will deal with those groups of order 64$p$ without a normal 
sylow $p$-subgroup \cite{4}.\\

We also give, in Tables 1c, 1c$^\prime $ and 1d, some partial results for the number of groups in the orders $2^np^2$ for $n \leq 8$. As one can see, the number of groups for $n\geq 5$ is very large, and a detailed study of the groups, order by order, is probably not feasible. If one is interested in certain subsets or classes of groups of these orders, this may or may not be feasible, e.g., by looking at groups of order $2^np^2$ whose action of the 2-group on the $p$-group is by a $D_4$ or a $Q_2$ action. For $n=5$, such a study was attempted in \cite{12} with only partial results. For the case of $n=6$, the number of groups even in this restricted subset increases very rapidly, yielding 374 cases for the $D_4$ action and 70 cases for the $Q_2$ action case. \\

\begin{tabular}{|c|c|c|c|c|c|c|c|c|c|}\hline
\multicolumn{10}{|c|}{Table 1a} \\ \hline
\multicolumn{10}{|c|}{Groups of order $2^m p$ for $m=3$, 4, 5, 6, 7 and 8} \\ 
\hline
\multicolumn{10}{|c|}{The case when the sylow $p$-subgroup is normal in 
$G$}\\ \hline
  & & \multicolumn{8}{|c|}{number of new or additional groups when $p\equiv 
1$ mod($2^n$)}  \\ \hline
$m$& direct  & $n=1$&$n=2$&$n=3$&$n=4$&$n=5$&$n=6$
&$n=7$ &$n=8$ \\ 
& products & & & & & & & & \\ \hline
3 & 5 & 7 & 2 & 1 & & & & &  \\
4 & 14& 28 & 9 & 2& 1& & & &\\
5 & 51&144 & 40& 9& 2& 1 & & &\\
6& 267&1120&243 &42&9&2&1& &   \\
7& 2,328&16,996& 2,180&262& 42&9&2&1&\\
8& 56,092&1,027,380&32,836&2,339&263&42&9&2&1\\ \hline 
\end{tabular}
\linebreak
\linebreak

\begin{tabular}{|c|c|c|c|c|c|} \hline
\multicolumn{6}{|c|}{Table 1b}\\ \hline
\multicolumn{6}{|c|}{Groups of order $2^np$} \\ \hline
\multicolumn{6}{|c|}{The cases without a normal sylow $p$-subgroup}\\
\multicolumn{6}{|c|}{(normal sylow 2-subgroup cases, no normal sylow subgroup cases)}\\ \hline
$n$ & $p=3$& $p=5$& $p=7$&  $p=17$& $p=31$ \\ \hline
3&  (2,1)& (2,0)     & (1,0)     & -   & -                      \\
4&    (6,4)& (1,0) & (1,0)&  -    &  -          \\
5&   (19,17)  &(2,1)   & (2,0)   &  -    & (1,0)   \\
6&  (70,86)& (5,5) & (9,0) &  -&  (1,0)         \\
7&  (309,536) &(13,21)   &(24,1) &(0,0)   & (2,0)               \\
8&  (1851,4912)& (49,104) & (77,4) & (1,0)  &   (5,0)           \\ \hline
\end{tabular}
\linebreak

The listing in Table 1c$^\prime $ is modeled on that for the groups of order $2^np$. As one goes to higher and higher orders (e.g., $2^np^2$ for $n \geq 5$) one needs a better way to list the results of the calculations. In these higher orders we also need to consider cases such as
$p\equiv 7$ mod(8), etc. or equivalently $p\equiv -1$ mod($2^n$). If in these higher orders one is only interested in the way in which the total number of groups varies with $p$, then the way these groups were listed in most of the early computations makes more sense and is shown in Table 1c. In the case of $n=5$, for example, we would need to consider the cases when $p\equiv $ 1, 3, 5, 7, 9, and 31 mod(32), as well as 15 and 17 mod(32). The most useful display, however, might be a combination of Tables 1c and 1c$^\prime $, showing the various quotient groups and how they act on the $p$-groups as a function of the prime $p$. These tables, however, can become very large and cumbersome as the order of the 2-group increases past $n=4$ or 5. \\

\begin{tabular}{|c|c|c|c|} \hline
\multicolumn{4}{|c|}{Table 1c (see \cite{13})}\\ \hline
\multicolumn{2}{|c|}{Groups of order $8p^2$}& \multicolumn{2}{|c|}{Groups of order $16p^2$} \\ \hline
prime & number of groups &prime & number of groups \\ \hline
$p=3$ & 50&$p=3$&197  \\
$p=5$ &52 &$p=5$ &221 \\
$p=7$ &44&$p=7$ &172  \\
$p\equiv 1$ mod(8)&60& $p\equiv 1$ mod(16)& 257\\
$p\equiv 3$ mod(8)&42& $p\equiv 3$ mod(8)& 167\\
$p\equiv 5$ mod(8)&52& $p\equiv 5$ mod(8)& 219\\
$p\equiv 7$ mod(8)&42& $p\equiv 7$ mod(8)& 169\\
& & $p\equiv 9$ mod(16) & 243\\ \hline
\end{tabular}
\linebreak

In our previous work on the groups of order 32$p$ we indicated how we estimated the number of groups of order 288 by using some counting arguments and the Hall-Senior Tables (or charts contained therein). The argument, there dealing with the 
$C_2\times C_2$ type images, says to just count the number of lines hitting the set of 
$C_2\times C_2$ images coming down from the normal subgroups of order 16, and this should give us the number of groups of order $32p^2$ with a $C_2\times C_2$ action on the $p$-group. When we went through and actually counted the number of such lines we did not get the exact same number as we got from determining the number of groups of order 288 which arise from the $C_2\times C_2$ actions of the order 32 groups on the $C_p\times C_p$ groups (for the cases of $p=3$ or 5). The numbers were sufficiently close so that one might guess that the differences arose from errors in the Hall-Senior charts. We tried a similar counting argument here for groups of order 64$p^2$ for $p=3$ and 5. Again some of the numbers are different, but again sufficiently close so that one might suspect the differences arise from omissions in the number of lines in the Hall-Senior charts. This might be an interesting question for the reader to look into. What we find is  indicated in Table 9. There are many such cases, and we give the number of such extensions as a function of the isoclinic class. The explicit numbers expected for each 2-group are given in Tables A5 and A6. We have not checked the $p=3$ cases other than to count the number of times the group $QD_8$ arises as a quotient group in the groups of order 64. The groups $\Gamma _2b$ and $\Gamma _2d$ of order 16, which appear in the Hall-Senior charts, only have a single line hitting them. Hence in these cases, according to our \textquotedblleft little empirical rule", we only get one new group for each occurrence of these order 16 actions (in order $64*5^2$) on the $p$-group, whereas in cases such as $C_2\times C_2$ we can have more than one group arising from that $C_2\times C_2$ image. These multiple lines seem to be associated with quotient groups that are direct products, e.g., $C_2\times C_2$, $C_4\times C_2$, $C_4\times C_4$, etc. Again this might be an interesting exercise for the reader to explain the reason for this apparent correspondence.\\

\begin{tabular}{|c|c|c|c|c|c|c|c|c|c|c|c|}\hline
\multicolumn{12}{|c|}{Table 1c$^\prime $} \\ \hline
\multicolumn{12}{|c|}{Groups of order $2^mp^2$ for $m=3$, 4, 5, and 6} \\ 
\hline
\multicolumn{12}{|c|}{The case when the sylow $p$-subgroup is normal in 
$G$}\\ \hline
  & & \multicolumn{10}{|c|}{number of groups for a given quotient group 
action}  \\ \hline
$m$& $p\equiv $1  &direct&   $C_2$&$C_2\times 
C_2$&$C_4$&$C_8$&(2,1)&$D_4$&$Q_2$&G(16) &$C_{32}$ \\
& mod($2^n$) &products& & & & & & & & & \\ \hline
3 &$n=1$& 10  & 21  &  6 & 2    & 1&0 &1 & 1& 0 &0 \\
  &$n=2$& 10  & 21  &  6 & 10& 1& 2&1 &1 & 0& 0\\
  &$n=3$& 10 &  21 & 6 & 10  &9 & 2& 1& 1& 0& 0\\ \hline
4 &$n=1$& 28 & 84  &35   &  9   & 2& 0&6 & 2& 1& 0\\
  &$n=2$& 28 & 84  &35 & 45& 2&14  & 6&2 &3 & 0\\
  &$n=3$& 28 & 84  &35 & 45    &16 & 14&6 &2 &11 &0 \\
  &$n=4$& 28 & 84  &35 & 45&16&14 &  6&2 &15 & 0\\ \hline
5 &$n=1$& 102 &432  &274& 40 &9 &0 &42 &11 &4 & \\
  &$n=2$& 102 &432  & 274& 200 & 9&112 &42  &11   &20$^*$  & 1\\
  &$n=3$& 102 &432  & 274& 40  & 81&112&42 &11 &  & \\ \hline
6 &$n=1$& 534 &3360 &3362   &241& 42&0 &374 & 70&36$^{**}$ & 0\\
  &$n=2$& 534 &3360 &3362&1213     &42 &1068 &374 &70 & 216$\dagger $&8 \\ \hline
\multicolumn{12}{|l|}{G(16) means a group of order 16. For $p\equiv 3$ mod(8) this group is $QD_8$ (in $16p^2$ case).} \\ \hline
\multicolumn{12}{|l|}{Note that this table does not explicitly contain the orders $p\equiv -1$ mod(8), but they} \\ 
\multicolumn{12}{|l|}{just differ from the $p\equiv 3$ mod(8) case by only two groups in the $16p^2$ cases, i.e., }\\ 
\multicolumn{12}{|l|}{delete $QD_8$ and add the three groups of order 16: $D_8$, $Q_4$, and $C_{16}$. } \\
 \hline
\multicolumn{12}{|l|}{$^*$ These 20 groups break up into 6 $C_4\times C_4$, 10 $C_4$Y$Q_2$, and 4 $<2,2|2>$ } \\ \hline
\multicolumn{12}{|l|}{$^{**}$ These 36 groups are all $QD_8$.} \\ \hline
\multicolumn{12}{|l|}{$\dagger $ The order 16 groups are: 63 $C_4\times C_4$, 125 $\Gamma _2b$ cases and 28 $\Gamma _2d$ cases.} \\ \hline
\end{tabular}
\linebreak

The number of groups for the case of $m=4$ and $n=4$ when the action is of order 16 is explicitly given in our 16$p^2$ article \cite{6}.\\

Another curiosity concerns the $C_4$ actions of the groups of order 64 on the group $C_p\times C_p$. In the case of $p=3$, we have 241 cases, according to the GAP/Small Group Library runs. In previous cases we find that the number of groups arising when $p=5$ is just five times the number arising in the $p=3$ case. This can be understood quite simply; one new set of groups comes from the action of the $C_4$'s on the group $C_{25}$ whose automorphism group is $C_5\times C_4$. The other four cases arise, as indicated in our 16$p^2$ group discussion, by the set of actions:
\begin{equation*}
a^5=b^5=(a,b)=a^c*a^2=b^c*b^x= \cdots =1,
\end{equation*}
where $x=1$,2,3, or 4. The case \textquotedblleft coming up" from the $p=3$ $C_4$ action corresponds to $x=3$. According to our GAP/Small Group Library runs there are 241 such $C_4$ cases arising in the groups of order 576, but there are 1213 cases when we run the order 1600 groups. There are thus an additional eight groups apparently arising here. This difference also might be an interesting exercise for the reader to explain why these extra 8 groups arise in order 1600.\\ 

As a by-product of the work on the groups of order 64$p$ with a normal sylow 
2-subgroup, we computed the automorphism groups for those groups of order 64 
possessing odd-order automorphisms. A brief outline of this work is given in 
Appendix I and Table A1. This material was not contained in an 
earlier report \cite{4}. In several cases these automorphism groups have a 
very simple representation, e.g., direct products of other well-known 
groups, three cases yield complete groups, and other cases have not yielded 
any simple interpretation. One of the side interests in \cite{4} was the 
finding of new complete groups. One way to do this is to find a group with a 
trivial center and then follow its automorphism group tower up to 
completion. In our attempts to come up with \textquotedblleft reasonable" or
\textquotedblleft simple" 
representations for the automorphism groups of the groups of order 64, we 
encountered several normal subgroups in these automorphism groups with a 
trivial center. In many cases these groups led us, by means of their 
automorphism group towers, to additional complete groups. In other cases, no 
such end, or termination point, was found. This work is presented in the 
appendix. This paper contains frequent references to groups of order 16, the 
automorphism groups of order $16p$, and groups of order 96. Convenient 
sources where these groups can be found
are papers \cite{6} and \cite{7}.\\

\begin{tabular}{|c|c|c|c|c|c|} \hline
\multicolumn{6}{|c|}{Table 1d}\\ \hline
\multicolumn{6}{|c|}{Groups of order $2^mp^2$ for $m=3$, 4, 5, 6, 7 and 8} \\ 
\hline
\multicolumn{6}{|c|}{The cases without a normal sylow $p$-subgroup}\\
\multicolumn{6}{|c|}{(normal sylow 2-subgroup cases, no normal sylow subgroup cases)}\\ \hline
$n$ & $p=3$& $p=5$& $p=7$&  $p=17$& $p=31$ \\ \hline
3& (6,4) &  (0,0)    & (2,0)     & 0     &  0                     \\
4& (13,17)   &(4,0)    & (2,1)     &  0    &  0          \\
5& (41,90)&(4,4)   &  (2,4)  &  0    & (2,0)  \\
6&  (152,510)&(10,24)  &(19,7)  & 0 & (2,$\geq 1$)         \\
7& (618,1072)$^*$   &   & (49,29)$*$&  0 &  (4,$\geq 2$)              \\
8&  &  &  &  (1,0) &      (10,$\geq 7$)       \\ \hline
\multicolumn{6}{|l|}{The group $1^5@C_{31}$ 's automorphism group is a 
complete}\\
\multicolumn{6}{|l|}{group of order 4950; hence the group $1^5@C_{31}$ only 
has }\\
\multicolumn{6}{|l|}{an outer automorphism of order 5. Hence there are no 
groups }\\
\multicolumn{6}{|l|}{of the form $\left(1^5@C_{31}\right)$@[2-group], but we 
can have groups}\\
\multicolumn{6}{|l|}{of the form $\left(1^5@C_{961}\right)$@[2-group] (number undetermined), and 
}\\
\multicolumn{6}{|l|}{of the form $1^5@C_{31} \times C_{31}$@[2-group].} \\ 
\hline
\multicolumn{6}{|l|}{$^*$For $2^7p$ when $p=3$ we have (309,536);} \\
\multicolumn{6}{|l|}{for $2^7p^2$ when $p=3$ we have at least (618,1072);} \\ 
\multicolumn{6}{|l|}{for $2^7p^2$ when $p=7$ we have at least (49,29).} \\ \hline
\end{tabular}
\linebreak\\

\section{The groups of order 64$p$ with a normal sylow 2-subgroup}

		The method of determining the groups of order 64$p$ with a normal sylow 
2-subgroup is straightforward. Just compute the automorphism groups for 
those groups of order 64 known to have an
odd-order automorphism. The particular groups that possess this property can 
be found by looking up the orders of these automorphism groups in the work 
of Hall and Senior \cite{3}. We then use CAYLEY to check the automorphism 
group, determine the sylow $p$-subgroup, and the actions of this subgroup on 
the generators of the group of order 64. In most cases the automorphism 
groups contain the prime only to a single power, so the above procedure works 
exactly as stated. In a few cases the sylow $p$-subgroup has order $p^2$ or 
$p^4$. In those cases we compute the subgroup lattice for these sylow
$p$-subgroups and determine the actions of the generators of each $C_p$ on 
the generators of the order 64 group. Using this method one should get all 
cases that will yield a group of order 64$p$ with a normal subgroup of 
order 64. This method may yield some duplicates which will need to be weeded out 
by hand, but these are only a few cases to look at and present no great 
burden. \\

The automorphism groups for the groups of order 64$p$ with a normal sylow 
2-subgroup yielded several new complete groups. In many other cases the 
automorphism groups of the groups of order 64$p$ have a simple breakdown into 
direct products of other well-known groups. The tabulated results of these 
calculations are given in Table 2. The complete groups of orders 168, 384, 
960, and 5760 are the same groups that arose in our other studies of groups 
of orders $8p$, $8p^2$, $16p$, and $16p^2$. The automorphism groups in many 
of the other cases involve well-known groups whose presentations are either 
well known (e.g., $S_4$) or which can be found in other published material 
(e.g., relations for 2-groups in \cite{3} or \cite{5}); others are 
explicitly written out here. Some of the groups appearing here as 
automorphism group factors are groups of order 192, without a normal sylow 
2-subgroup (e.g., Aut$(C_4 \times  C_2 \times  C_2)$), and explicit 
presentations for these groups are given below in the section devoted to 
these groups.\\

The automorphism groups for numbers 16 (= 192 \#1506), 19 (= 192 \#1507), 
and 63 (= 192 \#1009) were rerun in GAP to check on the identity of their 
automorphism groups. These automorphism groups are in fact, as suspected, 
isomorphic with the 1152 order factor being group number 155,478 of order 
1152 in the Small Group Library. Likewise the groups numbered 40 (= 192 \#1508), 
42 (= 192 \#1509), 55 (=192 \#1022) and 62 (= 192 \#1024) were rerun 
for the same purpose. The results indicate that the order 576 groups in 
these automorphism groups are also, as suspected, isomorphic and are 
isomorphic to order 576 number 8654 in the Small Group Library. \\

In one case ($1^6 @ C_3$) the order of the automorphism group is 
sufficiently large that very little other information about this group is 
available. The authors suspect that this group is a complete group.\\

In our initial listing of the groups of order 64$p$ with a normal sylow 2-subgroup we missed one of the groups for the case of $p=7$ (number 9 in Table 2b). The presentations for the groups of order 64$p$ of the form $1^6 @ C_7$ can be read from the following matrices:
\begin{equation*}
C_7 =\begin{pmatrix} 0&1&1\\
                     1&0&0\\
                     1&0&1 \end{pmatrix}\quad
C_7 =\begin{pmatrix} 0&1&1&0&0&0\\
                     1&0&0&0&0&0\\
                     1&0&1&0&0&0\\
                     0&0&0&0&1&1\\
                     0&0&0&1&0&0\\
                     0&0&0&1&0&1 \end{pmatrix}
\quad C_7=\begin{pmatrix}  1&0&1&1&1&1\\
                           1&1&0&0&0&1\\
                           0&0&0&1&1&0\\
                           0&0&1&1&1&0\\
                           0&0&0&1&0&0\\
                           1&1&0&0&1&0
     \end{pmatrix}
\end{equation*}
The middle $C_7$ matrix generates the group with the automorphism group of order 677,376 [448,\#1393]; the other $6\times 6$ matrix generates the one with an automorphism group of order 18,816 [448,\#1394].\\

The other $p=5$ and $p=7$ actions are given either in Table 2b, or in previous papers, e.g., the $C_5$ actions in \cite{6}, and the $C_{31}$ case in \cite{7}.\\

The results reported here can easily be extended to yield the groups of 
order $64p^2$ with a normal sylow 2-subgroup. These groups can be broken up 
into the following classes:\\

	     1. The group of order 64 having only a single odd-order automorphism 
of order $p$.\\

		a. A group of order $64p \times  C_p$.\\

		b. Take the group of order $64p$ and in its presentation
			change $C_3$ to $C_9$ (or more generally $C_p$ to $C_{p^2}$) without changing the actions on the order 64 group. This gives us groups of the type
			[64] @ $C_9$ (or [64] @ $C_{p^2}$).\\

             2. The sylow $p$-subgroup having order $p^2$ or higher. These cases 
are shown in Table 3.\\

In these cases one can have groups with the structure [64]$ @ (C_3 \times  
C_3)$ which may not be easily related to a group of order $64p$. Table 3 here shows that there are only a few cases of this type occurring here. The elementary abelian group of order 64 
has sylow $p$-subgroups of orders $3^4$ and $7^2$, along with the ones of 
orders 5 and 31. The group (2,$1^4$) has a sylow 3-subgroup of order $3^2$.\\

\section{ Groups of order 192 without a normal sylow subgroup}

A systematic way of determining the groups of order 192 without a normal 
sylow 2-subgroup would be to look for 2-group extensions of groups, $G$, of 
orders $3*2^n$ with $n < 6$, by 2-groups. This idea was tested out on the 
groups of orders 48 and 96, and this method yielded all of the groups without 
a normal sylow subgroup. If one recalls that there are 28 groups of order 
$16p$ produced via a $C_2$ action of the groups of order 16 on the $C_p$ 
group, it should be no surprise that there are 28 groups of the form:\\

$A_4$ @ [group of order 16].\\

In fact, the automorphism groups of these groups are just\\

$S_4 \times $ invariant factor appearing in Table 2a of \cite{6}.\\

A similar remark can be made for finding the groups of order 192 coming from the groups of order 24, extended by the groups of order 8. Here the group in question is SL$(2,3)$. 
Note that there are only two groups of order 8 with an automorphism of order 3: $C_2 \times  C_2 \times  C_2$ and $Q_2$. The first group extended by $C_3$ yields $A_4 \times  C_2$, and $Q_2$ yields SL$(2,3)$.\footnote{Note, however, in this case, Aut($G$)\ $\not\simeq\ $Aut(SL(2,3)) $\times $ Invariant factor, where this invariant factor is the same as that arising in the groups of order 8$p$ when we have a normal sylow $p$-subgroup. See Table 4b below.} All extensions of $A_4 \times  C_2$ by a group of order 8 will 
yield duplicates of groups produced by the $A_4$ @ [groups of order 16]. For 
the groups of order 16 there are four groups which have an automorphism of 
order 3, namely $(1^4)$, $(2,1^2)$, (2,2) and $C_4 \text{Y} Q_2$. Of these 
groups, $(2,1^2)$ will give a group of order 48 of the form  $A_4 \times  
C_4$, which will again duplicate results from extensions of the form $A_4$ @ 
[group of order 16]. The group $(1^4)$ has two extensions in order 48. They 
are $A_4 \times  C_2 \times  C_2$ and $(C_2 \times  C_2 \times  C_2 \times  
C_2) @ C_3$. The first case again duplicates results from the $A_4$ 
extensions. The second group will yield new groups in order 192.\\

Of the groups of order 48 that do not have a normal sylow subgroup, only 
$<2,3,4>$ was not generated by a simple semi-direct product from a group of 
order 24 (or lower) with a normal sylow
2-subgroup. The group $<2,3,4>$ can be obtained as a nonsplit extension of 
SL$(2,3)$ and $C_2$, with the presentation:\\
\begin{gather*}
		a^3=a*b*a*((b*a*b)^{-1})=c^4=c^2*(a*(b^{-1}))^2\\
			=a^c*a*b*(a^{-1})=b^c*b=1.
\end{gather*}
The action of the element [$c$] on SL$(2,3)$ is the same as giving rise to 
GL$(2,3)$.\footnote {That is, if we set $c^4\rightarrow c^2$ and omit the relation 
$c^2*(a*b^{-1})^2$ what we get is a presentation for GL$(2,3)$.}\\

Table 4 gives those groups of order 192 without a normal sylow subgroup 
that arise as split extensions (i.e., semi-direct products) of a group with 
a normal sylow 2-subgroup and a 2-group. Some of the groups here arise from 
two or more different extensions. In many cases these \textquotedblleft repeat cases" have 
been indicated for the reader.\\

For order 192 there are several cases of nonsplit extensions, which arise 
in exactly the same way as the group $<2,3,4>$ does in order 48. The first 
two such cases, not involving $<2,3,4>$, were pointed out to us by Dr. 
Antonio Vera Lopez \cite{11}. In Dr. Lopez's work classifying groups with 
13 conjugacy classes he found four groups of order 192 with just 13 classes: 
the two that we found, arising as split extensions, namely
\begin{equation*}
		       \text{Aut}(C_4 \times  C_2 \times  C_2) \qquad \text{ and} \qquad   
\text{Hol}(Q_2),
\end{equation*}
and two others. The other two cases are described by Dr. Lopez as nonsplit 
extensions of $[Q_2 $Y$ Q_2] @ S_3$. In Table 5 we list those nonsplit 
extensions giving rise to groups of order 192 that have presentations 
analogous to that yielding the group $<2,3,4>$ given above. The 
presentations given in Table 5 are modelled on those given in Table 4, in 
that the only modifications to those in Table 4 are finding the generator of 
the center of the kernel of the extension, and then mapping some
power of the element coming from the quotient group of the extension onto 
the generating element of the center of the kernel. The only possible 
choices for the kernels in these cases are displayed in Table 6.\\

The group SL$(2,3)$ presents an interesting case, and we show in Table 7 the 
relationship between the split and nonsplit extensions for the groups of 
order 192 arising from a normal SL$(2,3)$. For the $C_8$, $D_4$ and $Q_2$ 
groups, the quotient group from which the $C_8$, $D_4$ or $Q_2$ image comes 
is given instead of the extension itself. The notation $C_2$  $[C_4]$ 
means that the generator of order 4 in the group acts as an element of order 
2 on SL$(2,3)$. The numbers in the split extension column and the last 
column refer to the number of the group in Tables 4 and 5.\\

In Table 7 we have tried to show some of the groups of order 192 that occur as a semi-direct product of SL$(2,3)$ and the groups $D_4$ and $Q_2$, and how they are related to some other groups that are not semi-direct products. In each case, the quotient group
\textquotedblleft $D_4$" or \textquotedblleft $Q_2$" is modified or altered, and we show how this \textquotedblleft alteration" is mapped into the center of SL$(2,3)$. The best way to explain this is just by showing the presentation used here. In the first case, i.e. \#33, we have 
\begin{gather*}
a^3=aba(bab)^{-1}=c^4=c^dc=a^caba^{-1}=b^cb\\
=(a,d)=(b,d)=\begin {cases}d^2 \quad \qquad \qquad\qquad [33]\\ d^4=d^2(ab^{-1})^2\quad\,\,  [72]\end {cases} = 1.
\end{gather*} 
One should also note that some of our normal semi-direct product cases can be obtained from some apparent nonsplit extensions as well. To see this, look at 
group \#33 and group \#59, or \#34 and \#64: 
\begin{gather*}
a^3=aba(bab)^{-1}=d^2=c^dc=a^caba^{-1}=b^cb\\
=(a,d)=(b,d) =\begin{cases} c^4 \quad \qquad \qquad \qquad [33]\\ c^8=c^4(ab^{-1})^2\quad\,\,\, [59]\end{cases} = 1.
\end{gather*}
The nonsplit extensions using the form $c^8=c^4(ab^{-1})^2$ all yield groups which were obtained by semi-direct products involving other groups of order 24, 48 or 96.\\

The only nonsplit extensions in which the actions were greater than order 2 
are the $C_2 \times  C_2$ cases listed in Tables 5 and 7. No nonsplit 
extensions were found here for any other cases. All of the actions for the 
higher-order groups with nonsplit extensions were also of order 2. In our 
initial search for nonsplit extensions our list was incomplete, in part 
because we did not have a systematic method to search for all of these 
groups in a manner analogous to that for the groups with a normal sylow 
2-subgroup, and for those cases without a normal sylow subgroup that arise 
as a semi-direct product of lower-order groups. In this context see the 
comment below at the end of this section.\\

The reader should also note that some of the relations given in Table 4 for 
some of the groups of order 96 are not the same as those given in Tables 3a 
and 3b of \cite{7}. In some of these cases they represent relations derived 
from a degree 8 permutation representation for that group which has fewer 
generators than those in Tables 3a and 3b of \cite{7}. These relations and 
their permutation representations can be found in the work of J. Burns 
listed in Appendix 0 of \cite{4}. \\

The number of nonisomorphic groups listed in Tables 4 and 5 is 81. Comparing this 
with the number given in the Small Group Library of Besche and Eick, we find 
that they list 86 groups without a normal sylow subgroup. For the cases when the normal subgroup was either GL$(2,3)$ or the Coxeter group $<2,3,4>$ we got many duplicate cases. We show these cases in Tables A2 and A3 of the Appendix. In our attempts to identify which groups of order 192 we were missing we prepared Table A4, which shows the class/order structure for the groups we found, and cross-referenced them with the groups of order 192 in the Small Group Library. The missing groups of order 192 are shown in Table 8. From this we were led to the five groups of order 192 that we missed. The ones missing from our list are \#'s 949, 950, 954, 1489, and 1490 of order 192 in the Small Group Library \cite{10}.  \\

Some of these groups (e.g., \#950 in Table 8) were apparently missed because 
they possessed the same order structure and automorphism groups as others 
in Table 4. It is not clear why we missed the last one, \#954 of order 192.\\

The groups numbered 1490 and 1489 in Table 8 have the structural form
\begin{equation*}
\left[\left(Q_2 \times C_2\times C_2\right)@C_3\right] @ C_2.
\end{equation*}
For some reason in our initial working with the extensions of the order 96 group $\left(Q_2\times C_2\times C_2\right)@C_3$ it was believed that this order 96 group only gave rise to extensions with a normal sylow 2-subgroup. If one looks at the normal subgroups of these groups (\#1489 and \#1490) one finds that $\left[\left(Q_2 \times C_2\times C_2\right)@C_3\right] $ does appear as a normal subgroup in both groups. Subsequently we obtained the following presentations for these last two groups:
\begin{gather*}
a^4=b^4=a^2*b^{-2}=a^b*a=c^2=d^2=(c,d)=(a,c)=(b,c)=(a,d)=(b,d)=\\
e^3=a^e*b^{-1}=b^e*b^{-1}*a^{-1}=c^e*d=d^e*c*d=\\
\begin{cases}h^2= a^h*a=b^h*b*a^{-1}=c^h*c=d^h*c*d=e^h*e=1,\\
            h^4=a^h*h=b^h*b*a^{-1}=c^h*c=d^h*c*d=e^h*e=h^2*b^{-2}=1. \end{cases}
\end{gather*} 

The first group is just a semi-direct product of the order 96 group with a $C_2$. The second group, number 1489 in the Small Group Library, is a nontrivial central extension, obtained from group number 1490 by the replacement of $C_2$ by a $C_4$ and with the extra relation $[c]^2$ $*$ (center of 96 group)$^{-1}$. It might also be of interest to point out that in the first case we have
\begin{gather*}
 <a,b,e>\ \simeq \textup{SL}(2,3),\quad <a,b,e,h>\ \simeq \textup{GL}(2,3),\quad <c,d,e> \ \simeq A_4, \\
 <c,d,e,h>\ \simeq S_4   
\end{gather*}
and in the second case,
\begin{gather*}
 <a,b,e>\ \simeq \textup{SL}(2,3),\quad <a,b,e,h>\ \simeq\ <2,3,4>,
\quad <c,d,e>\ \simeq\ A_4, \\
 <c,d,e,h>\ \simeq\ A_4@C_4.  
\end{gather*}
An interesting observation also is that if in the second form (i.e., $h^4=\cdots)$, 
if we replace $h^2*b^{-2}$ with $h^2*d^{-2}$ we get an alternate presentation for the first of these two groups.\\

In the case of the last missing group, \#954 of order 192, we know that it has the structure SL$(2,3) @ D_4$. It has three normal subgroups of order 96, namely:
\begin{equation*}
\textup{GL}(2,3) \times C_2, \qquad <2,3,4> \times\ C_2, \qquad\text{and}\quad \textup{SL}(2,3)\times C_4.
\end{equation*}
In our attempts to get a \textquotedblleft reasonable" presentation for this group we tried several permutations or modifications of the actions on the group SL$(2,3)$ given in Table 5. The extension forms we tried were based upon one or the other assumed structures:
\begin{equation*}
\textup{SL}(2,3)@D_4, \qquad (\textup{SL}(2,3)\times C_4)@C_2.
\end{equation*}
The only result from this exercise was our getting several different/alternate presentations of other groups of order 192 given in these tables. The Small Group Library gives rise to a  presentation for this group on 7 generators with 19 relations. A slightly reduced set of relations can be obtained by eliminating the sixth and seventh generators, yielding the presentation:
\begin{gather*}
   a^2=
   b^4=
   c^2=
   d^3=
   e^4=
   (a,c)=
   (a*d^{-1})^2\\=
   (b,d)=
   (c,d)=
  (b,c)=
   (b,e)=
   (c,e)=
   e^2*b^{-1}*c*b^{-1}\\=
   a*c*b*a*b^{-1}=
   (a*e^{-1}*d)^2=
   (d^{-1}*e^{-1})^3=1.
\end{gather*}

\section*{Acknowledgements}
Almost all of the calculations reported here were done on the DEC computer in the
Department of Cognitive and Linguistic Sciences at Brown University. We give special
thanks to Dr. James Anderson for letting us use the computers. Margaret Doll helped
with the installation of CAYLEY on the computers. Some of the calculations are in
the Appendices. Thanks also go to Dr. John Cannon for giving us the CAYLEY program to
use and for answering questions about using it.\\

We also used the programming system GAP for checking our earlier CAYLEY results. In this
endeavor, Dr. Steve Linton was instrumental in helping us with the programming.\\

The bulk of this report was written by May 29, 1997. The comments related to the Small Group Library and other work by B. Eick et al. were added at a much later date (2006 to 2008). \\

\section{Appendix I\\
                The automorphism groups of the groups of order 64\\
                       with an odd-order automorphism}

In the course of the work done on determining the groups of order 192 with a 
normal sylow 2-subgroup we computed certain properties of the automorphism 
groups of the groups of order 64 which have an odd-order automorphism. The 
main purpose here was to determine the actions of the resulting operators of 
order 3 in the automorphism group on the group of order 64, thereby yielding 
the desired groups of order 192. In the course of this study we thought it 
would be interesting to determine presentations and other properties of 
these automorphism groups. A partial summary of this work is contained in 
Table A1. This table shows that in some cases these automorphism groups have a rather simple structure, while in other cases the orders are sufficiently large that direct computational procedures are not especially 
useful in determining some properties (e.g., their automorphism groups). In what follows let U64 denote the set of automorphism groups of the groups of order 64 which possess odd-order automorphisms.\\

The two main areas of interest here were to obtain some information on the 
structure of the groups in U64, and second to see if we could find any new 
complete groups related to groups in the set U64. As the reader can see from 
Table A1 we have been successful in finding several new complete groups. We 
have not been as fortunate in arriving at some simple structural properties 
for many of the groups in U64 as we would have liked.\\

The numbers appearing after the order of an automorphism group in Table A1 are to be
interpreted as follows. For example, 
for group \#3, 147,456(2) means that the automorphism group of the group number 3 of order 64 has order 147,456 and has center of order 2; in other cases, e.g., \#15, the ($1^2$) means that the center for this group is $C_2 \times  C_2$,\ldots .  The presentations for many of the automorphism groups of the groups of order 64 are given in Appendix II. Complete presentations for many of the groups in the automorphism tower sequences listed in Table A1 are given in Appendix III. Those automorphism towers leading to complete groups were originally planned for inclusion in a future paper devoted to the complete 
groups arising in our studies of automorphism groups. They are included here for now for the sake of completeness.\\

The groups in U64 for the most part do not have a normal sylow subgroup, 
which means that the structure of these groups takes the schematic form:
\begin{equation*}
			\text{Aut} = ( T1 \,\,@ \,\,C_3 )\,\, @\,\, T2,
\end{equation*}
where $T1$ and $T2$ are 2-groups. The automorphism groups of the groups $T1$ can 
be very large, but usually they have relatively small odd-order sylow 
subgroups; typically these groups have orders such as $3*2^n$, where $n$ is 
often a number in the range 10 to 20. It is therefore very easy to obtain a 
presentation for the group $T1 @ C_3$. But to complete the sequence one needs 
to know how $T2$ acts on the group $T1 @ C_3$, and these groups often have 
very large sylow 2 components which make finding the relevent 2-group action 
on the required normal subgroup very difficult. In most cases the group $T1 
@ C_3$ has a trivial center, and we have followed the
automorphism towers for these groups to obtain new complete groups.\\

There are three groups in U64 that have normal sylow 2-subgroups. These are 
the automorphism groups for the following groups of order 64: \#76, \#93 and \#153.\\ 

The groups \#3 and \#14 of order 64 have the structure:
\begin{equation*}
			\text{Aut} = ( T1\,\, @\,\, (C_3 \times  C_3) )\,\, @\,\, T2.
\end{equation*}
The automorphism group for \#187 of order 64 is more 
complicated and appears to be representable as:
\begin{equation*}
			\text{Aut} = (( 1^4 \,\,@\,\, C_3) \times  1^4 )\,\, @\,\, C_5 )\,\, 
@\,\, C_4.
\end{equation*}
The first five cases and numbers 14, 103, 104 and 105 in Table A1 have not 
been looked at for this breakdown. Numbers 1 and 2 appear to be too large 
to do this with CAYLEY.\\

\newpage

\section{Tables 2a to 9}
\begin{tabular}{|l|c|l|}\hline
\multicolumn{3}{|c|}{Table 2a}\\ 
\hline
\multicolumn{3}{|c|}{Groups of order 192 with a normal}\\
\multicolumn{3}{|c|}{sylow 2-subgroup and their automorphism groups}\\ 
\hline
    &       Group of order 192        &          Automorphism Group\\ 
\hline
       1.& $  A_4 \times 1^4                $& $         S_4\ \times\ $GL$(4,2)$ \\
       2. & $  1^4 @ C_3 \times 1^2         $& $           S_3 \times [5760] 
                $ \\
       3. & $ 1^6 @ C_3                $&           order is [23,224,320]    
   \\
       4. & $ A_4 \times (2,1^2)            $& $          S_4 \times 
\text{Aut}(C_4 \times C_2 \times C_2)     $ \\
\hline
\end{tabular}
\linebreak\\

\begin{tabular}{|l|c|l|}\hline
\multicolumn{3}{|c|}{Table 2a (continued)} \\ 
\hline
\multicolumn{3}{|c|}{Groups of order 192 with a normal}  \\
\multicolumn{3}{|c|}{sylow 2-subgroup and their automorphism groups}\\ 
\hline
    &       Group of order 192        &          Automorphism Group      \\ 
\hline
       33.& [32 \# 33] @ $C_3 \times C_2   $& $           C_2 \times$ [384] 
complete         \\
          & $C_2 \times$ Table 3a in \cite{7} \# 14   &                      
                  \\
       34.& $ C_2 \times$ ([32 \# 34] @ $C_3$) & $            C_2 \times$ 
[384] complete         \\
          & $ C_2 \times$ Table 3a in \cite{7} \# 15   &                     
                  \\
\hline
\end{tabular}
\linebreak\\

\begin{tabular}{|l|c|l|}\hline
\multicolumn{3}{|c|}{        Table 2a (continued)}                              
    \\ \hline
\multicolumn{3}{|c|}{Groups of order 192 with a normal} \\
\multicolumn{3}{|c|}{sylow 2-subgroup and their automorphism groups}\\ 
\hline
       61.& \# 183 a.  &  [36864] complete \\            
       62.& \# 183 b.  &  [576], $\textup{ncl}=16, z=1$ \\         

       63.& \# 183 c.  &  [1152], $\textup{ncl}=17$ \\             

       	    			&       & Aut(1152)=[2304] complete \\
       64.& \# 184   &  [384] complete \\             

       65.& \# 187    &  [61440] complete \\            

       66.& \# 241      & $   S_4 \times \text{Hol}(C_8) 
               $ \\
       67.& \# 242    & $  S_4 \times D_4 \times  C_2$ \\ 
         
       68.& \# 243   & $  S_4 \times \text{Hol}(C_8) $\\
              
       69.& \# 259     & $  C_2 \times \text{Aut}(Q_2 
\times C_2)          $ \\
       			   	&&	  Table 4 below \# 59.       \\
       70.& \# 260                    & $         C_2 \times \text{Aut}(Q_2 
\times C_2).         $ \\ \hline
\end{tabular}
\linebreak\\

Note that the numbers in the group 192 column refer to the sylow 2-subgroup; e.g., 
\#64 means that the order 192 group is \#64 @ $C_3$.

\begin{tabular}{|l|c|} \hline
\multicolumn{2}{|c|}{  Notes for Table 2a }  \\ \hline
\multicolumn{2}{|c|}{The complete group [2304] is the same as that obtained}\\
\multicolumn{2}{|c|}{from Aut[Aut$[(Q_2 \times  C_2 \times  C_2) @ C_3]]$. 
}\\ \hline
\multicolumn{2}{|c|}{Presentations for nondirect product cases.  }\\ \hline
\#27. [64\#30]& $ a^4=b^4=c^4=(b,c)*a^2=(a,b)=(a,c)=  $  \\
  	&	$d^3=(a,d)=b^d*(c^{-1})=c^d*b*c=1;          $ \\ \hline
\#29. [64\#36] &$ b^2=c^2=(c,b)*a^8=(a,b)=(a,c)=    $                   \\
         &$d^3=(a,d)=b^d*b*c*b=c^d*b*c*(a^{-4})=1; $   \\
\hline
\end{tabular}
\linebreak\\

\begin{tabular}{|l|c|} \hline
\multicolumn{2}{|c|}{  Notes for Table 2a (continued)}  \\ \hline
\multicolumn{2}{|c|}{The complete group [2304] is the same as that obtained   }\\
\multicolumn{2}{|c|}{from Aut[Aut$[(Q_2 \times  C_2 \times  C_2) @ C_3]]$. 
}\\ \hline
\multicolumn{2}{|c|}{Presentations for non-direct product cases.  }\\ \hline
\#47. [64\#144]&$ a^2=b^2=c^2=(a,b)^2=(a,c)^2=(b,c)^2=                 $\\
  		&$((a,b),c)=((a,c),b)=                                       $\\
  		&$d^3= a^d*c*b*c = b^d*b*c*b = c^d*b*c*a*c*b=1;              $\\ \hline
  	\#48. [64\#145]&$ a^2*(b^{-2})=c^4=(b,a)*a^2=                            
$\\
\hline
\end{tabular}
\linebreak\\

\begin{tabular}{|l|c|} \hline
\multicolumn{2}{|c|}{Notes for Table 2a (continued)}  \\ \hline
\multicolumn{2}{|c|}{The complete group [2304] is the same as that obtained   }\\
\multicolumn{2}{|c|}{from Aut[Aut[$(Q_2 \times  C_2 \times  C_2) @ C_3]]$. 
}\\ \hline
\multicolumn{2}{|c|}{Presentations for nondirect product cases.  }\\ \hline
  	\#61. [64\#183]&$ c^2=d^2=(c,a)*a^2=(a,d)*(b^{-2})*a^2=                  
$\\
  		&$(b,c)*b^2*(a^{-2})=(d,b)*b^2=(a,b)=(c,d)=e^3=                $\\
  		&$	a^e*a*(b^{-1})=b^e*a=c^e*((a^2*c*d)^{-1})=             $\\
  		&$	d^e*b*(c^{-1})*(b^{-1})=1;                             $\\
\hline
\end{tabular}
\linebreak\\

\begin{tabular}{|l|c|l|}\hline
\multicolumn{3}{|c|}{ Table 2b } \\ \hline
\multicolumn{3}{|c|}{ Groups of order $64p$ with a normal }\\
\multicolumn{3}{|c|}{sylow 2-subgroup and their automorphism groups}\\ 
\hline
   &   Group of order $64p$       &         Automorphism Group         \\ 
\hline
      &  $p=5$ (5 cases)   &                                              \\ 
\hline
       1.& $[1^4] @ C_5 \times 1^2$ &                  $S_3 \times$ [5760] 
complete  \\
       2.& $[1^4] @ C_5 \times C_4$ &                  $ C_2 \times $[5760] 
complete  \\
       3.& \# 104                   &           $C_2 \times $[320] complete  
  \\
\hline
\end{tabular}
\linebreak\\

\begin{tabular}{|l|l l|c|}\hline
\multicolumn{4}{|c|}{Table 3} \\ \hline
Group \# & \multicolumn{2}{|c|}{order of sylow $p$-subgroup of Aut$(G)$} & 
\# of 64$p$ groups \\ \hline
\#1 GL$(6,2) $& $3^4$ and $7^2$, & all others are $C_p$ & 3, 3 \\
\#2 $\left(2,1^4\right)$  & $3^2$, & all others are $C_p$ & 2 \\
$\left(2^2,1^2\right)$ & $3^2$& & 3 \\
$Q_2\times 1^3$ & $3^2 $ && 3 \\
$C_4 \text{Y} Q_2 \times 1^2$ & $3^2$ && 3 \\
\# 103 & $3^2$ &&2 \\
\#105 & $3^2$ && 2 \\
\#183 & $3^2$ && 3 \\ \hline
\end{tabular}
\linebreak \\

\begin{tabular}{|l|}  \hline
\multicolumn{1}{|c|}{Table 4        April 11, 1996       }\\
\multicolumn{1}{|c|}{ Groups of order 192 without a normal  }\\
\multicolumn{1}{|c|}{ sylow subgroup and their automorphism groups } \\
\multicolumn{1}{|c|}{   split extensions       } \\ \hline
      A. Groups of the form $A_4$ @ [order 16 group] (28 cases) \\
\multicolumn{1}{|c|}{ [$A_4$: $x^2=y^3=(y*x)^3=1$. }                      \\
\multicolumn{1}{|c|}{ $C_2$ actions: $(g,x)= y^g*y=1$]}                 \\
                                                                   \\
        	    There are 28 such cases and they take the same form \\
       as the groups of order $16p$. In these cases, however,    \\
       instead of acting on a group $C_p$, the action is by the same \\
       element of order 2, and the automorphism group is   \\
       $S_4 \times $ the invariant factor listed in Table 2a in \cite{6}. \\
                              \\
           	In this set of groups the following direct          \\
       products of groups of order 96 without a normal sylow    \\
       2-subgroup with $C_2$ can be found:                           \\
                                                                     \\
        	    \#'s 20 to 26 of Table 3b of [7].       \\ \hline
\end{tabular}
\linebreak\\

\begin{tabular}{|l|c|l|} \hline
\multicolumn{3}{|c|}{Table 4 (continued) }\\ \hline
\multicolumn{3}{|c|}{B. Groups of the form SL$(2,3)$ @ [group of order 8] 
(13 cases)}\\ \hline
\multicolumn{3}{|c|}{[SL$(2,3)$: $a^3=a*b*a*((b*a*b)^{-1})=1$; } \\ \hline
\multicolumn{2}{|l}{ actions:   }& \\ \cline{2-3}
$C_2$:   & \multicolumn{2}{c|}{ $a^x*a*b*(a^{-1})=b^x*b=1; $ }      \\  
\hline
$C_4$:   & \multicolumn{2}{c|}{ $a^x*b=b^x*b*a*(b^{-1})=1; $ }        \\ 
\hline
$C_2 \times  C_2$:   & \multicolumn{2}{c|}{ $a^x*a*(b^{-1})*(a^{-1})=$}  \\
   & \multicolumn{2}{c|}{  $ b^x*b*(a^{-1})*(b^{-1})= $ }      \\
\hline
\end{tabular}
\linebreak\\

\begin{tabular}{|l|c|l|} \hline 
\multicolumn{3}{|c|}{Table 4 (continued)}\\ \hline 
\multicolumn{3}{|c|}{B. Groups of the form SL$(2,3)$ @ [group of order 8]  
(13 cases)}\\ \hline 
\multicolumn{3}{|c|}{SL$(2,3)$: $a^3=a*b*a*((b*a*b)^{-1})=1$; } \\ \hline 
\multicolumn{3}{|c|}{The other cases yield only groups with a normal sylow  
2-subgroup. }\\ \hline 
\multicolumn{3}{|c|}{note: There are two $C_2 \times C_2$ images in  
Aut(SL(2,3)). }\\ 
\multicolumn{3}{|c|}{The following cases appeared to give duplicate groups.  
}\\ \cline{2-3} 
1.& \multicolumn{2}{c|}{The $Q_2$ case with a $C_2 \times C_2$ image }        
\\ 
    & \multicolumn{2}{l|}{yields one of the missing cases: \#950 (note added fall 2007).}      \\ 
  2. & \multicolumn{2}{c|}{One of the $C_2 \times C_2$ images from $D_4$  
appears to}      \\ 
   & \multicolumn{2}{l|}{give the same group as a $C_2$ case from $D_4$, \#34 above. }   
\\ 
  3. & \multicolumn{2}{c|}{The $C_4$ action from $C_8$ appears to give a }    
    \\ 
  & \multicolumn{2}{l|}{duplication of the $C_2$ image case from $C_8$, i.e., \#29 above. }     
\\ \hline 
\end{tabular} 
\linebreak\\ 
 
\begin{tabular}{|l|c|l|}\hline 
\multicolumn{3}{|c|}{ Table 4. (continued)    June 21, 1994 }\\ \hline 
\multicolumn{3}{|c|}{  C. Groups from $[C_4 \times  C_4] @ C_3$  by a group  
of order 4 (4 cases) }\\ 
\multicolumn{3}{|c|}{$ a^4=b^4=(a,b)=c^3=a^c*(b^{-1})*(a^{-1})= $} \\ 
\multicolumn{3}{|c|}{   $b^c*b^2*(a^{-1})=1; $ }\\ \hline 
& \multicolumn{2}{l|}{actions: } \\ 
$C_2$: & \multicolumn{2}{l|} {$a^x*b*a^2=b^x*a*b^2=c^x*c*a=1; $ }\\ 
  $C_4$:& \multicolumn{2}{l|}{$a^z*a*(b^{-1})=b^z*(b^{-1})*a^2=c^z*a*c=1;  
$}\\ 
  $C_2 \times  C_2$:& \multicolumn{2}{l|}{$a^x*b*a^2=b^x*a*b^2=c^x*c*a^2=$ }  
\\
\hline 
\end{tabular} 
\linebreak\\ 
 
\begin{tabular}{|l|l|l|}\hline 
\multicolumn{3}{|c|}{Table 4 (continued) June 21, 1994   }\\ \hline 
\multicolumn{3}{|c|}{D. Groups from  $[C_4 $Y$ Q_2] @ C_3$  (there are 6 cases  
here)}\\ \hline 
\multicolumn{3}{|c|}{$a^3=a*b*a*((b*a*b)^{-1})=c^2= $} \\ 
\multicolumn{3}{|c|}{$a^c*b*(a^{-1})*(b^{-1})=b^c*a*(b^{-1})*(a^{-1})=1;  
$}\\ \hline 
  &Actions: & \\ \cline{2-3} 
   &$C_2$: & \\ 
   &1. $ a^x*b*a*(b^{-1})=b^x*b=$ & \\ 
   &  \qquad $c^x*((a*b*c*a)^{-1})=1$ & \\
\hline 
\end{tabular} 
\linebreak\\ 
 
\begin{tabular}{|l|c|l|}   \hline 
\multicolumn{3}{|c|}{Table 4 (continued) June 21, 1994 }\\ \hline 
\multicolumn{3}{|c|}{ E. Groups from $[C_2 \times  C_2 \times  C_2 \times   
C_2] @ C_3$  (four cases here) }\\ \hline 
\multicolumn{3}{|c|}{ Relations for 48 group:        }\\ 
\multicolumn{3}{|c|}{$a^2=b^2=c^2=d^2=(a,b)=(a,c)=$} \\ 
\multicolumn{3}{|c|}{$(a,d)=(b,c)=(b,d)=(c,d)=e^3=a^e*b=b^e*a*b= $}\\ 
\multicolumn{3}{|c|}{$ c^e*d=d^e*c*d=1;      $}\\  \hline 
& $C_2$ action: & \\ 
  & $ a^f*c=b^f*c*d=c^f*a=d^f*a*b=e^f*e=1;$&\\
\hline 
\end{tabular} 
\linebreak\\ 
\begin{tabular}{|l|c|l|}\hline 
\multicolumn{3}{|c|}{ Table 4 (continued)  June 21, 1994 }\\ \hline 
\multicolumn{3}{|c|}{F.  Groups from the groups of order 96.     }\\ \hline 
\multicolumn{3}{|c|}{16 cases here, not all different.}\\\hline 
\multicolumn{3}{|c|}{$[C_8 $Y$ Q_2] @ C_3$ (four cases) } \\ \hline 
& \multicolumn{2} {c|} {Relations for the order 96 group.  [Table 3a in  
\cite{7} \# 12]} \\ 
&\multicolumn{2}{c|} {$b^3*(a^{-3})=b*a*b*((a*b*a)^{-1})=b*a^2*b*a^8=1; $}  
\\ \cline{2-3} 
     \#56. & $    c^2=a^c*a^2*(b^{-1})=b^c*a*(b^{-1})*a=1;$ & $    S_4  
\times D_4 \times C_2 $ \\ 
     \#57. & $    c^2=a^c*((a*b^2*a^2)^{-1})=          $ & $S_4 \times C_2  
\times C_2 \times C_2 $\\
\hline 
\end{tabular} 
\linebreak\\ 
\begin{tabular}{|l|c|l|}\hline 
\multicolumn{3}{|c|}{ Table 4. (continued)      June 21, 1994       }\\ \hline 
\multicolumn{3}{|c|}{F.  Groups from the groups of order 96.     }\\ \hline 
\multicolumn{3}{|c|}{16 cases here, not all different.}\\\hline 
            &  $[32 \# 41] @ C_3$   &                                         
   \\ 
         &         All cases here have normal sylow 2 subgroups!&        \\ 
        &         order 96 group is [Table 3a in \cite{7} \# 16]   &          
     \\ 
                        &&                                                  
\\ 
         & $[Q_2 \text{Y} Q_2] @ C_3$&                one case   \\
\hline 
\end{tabular} 
\linebreak\\ 
\begin{tabular}{|l|c|l|}\hline 
\multicolumn{3}{|c|}{ Table 4. (continued)      June 21, 1994       }\\ \hline 
\multicolumn{3}{|c|}{F.  Groups from the groups of order 96. (continued)      
}\\ \hline 
       &$ [Q_2 \text{Y} D_4] @ C_3$  [Table 3a in \cite{7} \# 19]  &        two  
cases   \\ \hline 
               &    Relations for order 96 group.        &                  
\\ 
        &$       c^2=c*b*(a^{-1})*b*c*(a^{-1})= $&                       \\ 
        &$          c*(a^{-1})*b*c*(b^{-1})*a=b^3*(a^{-3})=$&               
\\ 
        &$           b*a*b*((a*b*a)^{-1})=(b*a^2)^2=1; $&               \\  
\hline 
     \#64. &     $  d^2=a^d*a*b*(a^{-1})=b^d*b=$          & $S_4 \times D_4  
\times C_2$      \\ 
          &     $  c^d*((a^2*c*a)^{-1})=1;  $       &   \\ \cline{2-3} 
      &  $  d^2=a^d*(c^{-1})*a*(c^{-1})= $   & $S_4 \times D_4 \times C_2$    \\ 
          &  $  b^d*b^2*(a^{-1})=(c,d)=1;$       &    \\  
          &              looks like \# 49.   &      \\ \cline{2-3} 
     \#65. & $ d^2=a^d*a^2*(b^{-1})=b^d*(c^{-1})*b*c= $& $S_4 \times D_4  
\times C_2$      \\ 
          &     $      c^d*a*(c^{-1})*(a^{-1})=1; $       &                   
  \\ \hline 
\end{tabular} 
\linebreak\\ 
 
\begin{tabular}{|c|c|c|}\hline 
\multicolumn{3}{|c|}{ Table 4 (continued)     June 21, 1994} \\ \hline 
\multicolumn{3}{|c|}{G. Groups coming from the groups GL$(2,3)$ and  
$<2,3,4>$.} \\ \hline 
\multicolumn{3}{|c|}{Note that most of the cases arising here are duplicates of   
} \\ 
\multicolumn{3}{|c|}{previous listings. These alternate presentations, which  
}  \\ 
\multicolumn{3}{|c|}{might be of interest to the reader, are given in Tables A2 and A3.}\\   
\multicolumn{3}{|c|}{  }\\ 
\multicolumn{3}{|c|}{ Only the new cases of groups without a normal sylow  
subgroup } \\ 
\multicolumn{3}{|c|}{are given here.  }\\
\hline 
\end{tabular} 
\linebreak\\ 
\begin{tabular}{|l|c|l|}\hline 
\multicolumn{3}{|c|}{Table 5} \\ 
\multicolumn{3}{|c|}{ Groups of order 192 without a normal} \\ 
\multicolumn{3}{|c|}{  sylow subgroup and their automorphism groups:} \\ 
\multicolumn{3}{|c|}{  nonsplit extensions }                       \\  
\hline 
  \multicolumn{3}{|c|}{ 1. Groups coming from SL$(2,3)$. }                \\  
\hline 
    &  SL$(2,3): a^3=a*b*a*((b*a*b)^{-1})=1$;&                        \\  
\hline 
    &     actions:                             &                      \\ 
    &             $ C_2: a^x*a*b*(a^{-1})=b^x*b=1$;&                         
\\
\end{tabular} 
%\end{document} 
    \#76.        &    $<-2,4|2>$ with $c$ acting on SL$(2,3)$   & $C_2  
\times  D_4 \times  S_4$\\ 
                 & This can also be written as $\left[\left(D_4\text{Y}Q_2\right)@C_3\right]@C_2.$&  \\ & 
&                                                  \\   
\hline 
\end{tabular} 
\linebreak\\ 
\begin{tabular}{|l|c|l|}\hline 
\multicolumn{3}{|c|}{        Table 5}                                 \\ 
\multicolumn{3}{|c|}{ Groups of order 192 without a normal  }            \\ 
\multicolumn{3}{|c|}{  sylow subgroup and their automorphism groups:}          
\\ 
\multicolumn{3}{|c|}{  nonsplit extensions }                       \\  
\hline 
    \#77 *.  &  SL$(2,3) @ Q_2$  &                                $1^4  
\times  S_4$ \\ \cline{2-2} 
                   &    $c^4=d^4=c^d*c=c^2*d^2*(a*(b^{-1}))^2=1$; &         
\\ 
                   &    both $c$ and $d$ act as $C_2$ on SL$(2,3)$    &       
    \\ \hline 
                                   &&                                 \\
\hline 
\end{tabular} 
\linebreak\\ 
\begin{tabular}{|l|c|l|}\hline 
\multicolumn{3}{|c|}{        Table 5  (continued) }                              
     \\ 
\multicolumn{3}{|c|}{ Groups of order 192 without a normal  }            \\ 
\multicolumn{3}{|c|}{  sylow subgroup and their automorphism groups:}          
\\ 
\multicolumn{3}{|c|}{  nonsplit extensions }                       \\  
\hline 
              && \\ 
      &  $[Q_2 \text{Y} Q_2] @ C_3$   [Table 3a in \cite{7} \# 18]  & \\ 
       &        two cases (second 96 extension) &   \\ \hline 
               &   Relations for order 96 group.              &        \\
\hline 
\end{tabular} 
\linebreak\\ 
 
\begin{tabular}{|c|c|c|c|c|c|}\hline 
\multicolumn{6}{|c|}{  Table 6 } \\ \hline 
\multicolumn{6}{|c|}{Centers of groups of orders 24, 48, and 96 and  
nonsplit extensions }\\ \hline 
     order &     group           &    generator for    & $Z(G)$ & \# of  
&group\\ 
           &              &   center of group  &      &groups&  \#  \\  
\hline 
       24  &  SL$(2,3)$         &  $<(a*(b^{-1}))^2> $    &  $C_2$  &  8   
  &  $a$  \\ \hline 
       48  &$ (C_4 \text{Y} Q_2) @ C_3$      & $ <b*(a^{-1})*c> $   & $ C_4$  
  &  *   &     \\ \hline 
       96  &$ (Q_2 \times C_2 \times C_2) @ C_3$ &  $<a^2>$ or $ < b^2 >$   &  $C_2  
$ &  1   & 1489    \\ 
           &$   (C_8 \text{Y} Q_2) @ C_3   $ &     $ < a^3 > $       &   
$C_8$  &  **  &     \\ 
           &$   [32 \# 18] @ C_3   $ & $<a*(c^{-1})*b*(c^{-1})>$ &  $C_2$  &  
  1   & 79  \\ 
           &$   (Q_2 \text{Y} Q_2) @ C_3   $ &     $<a^2>$           &   
$C_2$  &  2   &80,81\\ 
           &$   (D_4 \text{Y} Q_2) @ C_3   $ &  $ <a*c*(a^{-1})*c>$    &   
$C_2$ & ***  &     \\ \hline 
\end{tabular} 
\linebreak\\ 
 
\begin{tabular}{|c|} \hline 
                           Notes for Table 6                        \\  
\hline 
                                                                      \\ 
            Note that some of the SL$(2,3)$ cases are just direct products,  \\ 
          e.g., $<2,3,4> \times C_4$ and $<2,3,4> \times C_2 \times C_2$.     
         \\ 
                                                                      \\ 
               a. Note that all of the $<2,3,4>$ extensions could also be     \\ 
                       nonsplit SL$(2,3)$ extensions. These have not been   \\ 
                       counted here as SL$(2,3)$ cases. See table below. \\ 
                                                                       \\ 
           *  All of the groups found here appear to be            \\ 
               duplicates of the SL$(2,3)$ extension by $C_8$.           \\ 
                                                                      \\ 
            ** The four cases here have the same order structure    \\ 
               and automorphism group as the split extension   \\ 
                $A_4 @ C_{16}$.                                       \\ 
                                                                       \\ 
               *** This case also arises as an extension from          \\ 
               SL$(2,3) @ D_4$  $(<-2,4|2>$ case).                  \\  
\hline 
\end{tabular} 
\linebreak\\ 
 
\begin{tabular}{|l|c|c|c|c|}\hline 
\multicolumn{5}{|c|}{Table 7} \\ \hline 
\multicolumn{5}{|c|}{ Nonsplit extensions arising from SL$(2,3)$ } \\  
\hline 
        2-group  &   action    &     split     &    nonsplit      &   \#    
\\ 
                 &             &   extension   &    extension      &        
\\ \hline 
          $C_8     $&$ C_2          $& \#29           &$     C_{16}           
  $&  71   \\ \hline 
        $C_4 \times  C_2  $&$ C_2[C_4]$      & \#31           &               
      &       \\ 
                 & $C_4[C_4]$      & \#36           &                   &     
    \\ 
                 & $C_2[C_2]$      & \#30  GL$(2,3)$  &  $<2,3,4> \times   
C_4$     &  69   \\ 
                 &             &        $ \times  C_4$   &                    
&       \\ 
                 & $C_2 \times  C_2$     & \#37           &   $  C_4 \times   
C_4$       &  78   \\ \hline 
        $C_2\times C_2\times  C_2 $&$ C_2$   & \#32 GL$(2,3)   $&$ <2,3,4>  
\times C_2 \times C_2$ &  67   \\ 
                 &             &     $\times C_2 \times C_2$ &                 
    &       \\ 
                 & $C_2 \times C_2$     & \#39           &$[<2,3,4> @ C_2]  
\times C_2$&  68   \\ \hline 
        $ D_4$      & $C_4[C_2]$      & \#33           & $C_4 @ C_4 [C_4$  
quot] &  72   \\ 
                 &             &        \#33   &    $(4,4|2,2)$      &  73    
\\ 
                 &             &        \#33   &       $Q_4[c]$       &  74   
  \\ 
                 &    $C_2[C_2]$     &        \#34   &      $ Q_4[d] $      &  75   
  \\ 
                 &  $C_2\times C_2$ &   \#40 &      &     \\ 
                 &   $D_4$          &   \#41 &       $<-2,4|2>$  & 76\\ \hline 
         $Q_2$   & $C_2$       & \#35    & $(C_4 @ C_4)$       & \#949 M   \\   
                 & $C_2\times C_2$ & \#950 M &      &             77  \\ \hline  
\end{tabular} 
\linebreak\\ 
\begin{tabular}{|c|c|c|c|c|c|c|c|c|} \hline 
\multicolumn{9}{|c|}{Table 8}\\ \hline 
Our & GAP & \multicolumn{6}{|c|}{order structure}& Aut(g) \\ \cline{3-8} 
\#  & \#  &  2 & 3 & 4 & 6 & 8 &12 &  \\ \hline 
35 &945& 3/3 & 8/1 & 76/7&24/3& 48/4 &32/4 & $(C_2\times C_2)\,\wr\,C_2  
\times S_4$ \\ 
- & 949 &3/3& 8/1 & 76/7& 48/4&24/3 &32/4 & $(C_2\times C_2)\,\wr\,C_2 \times  
S_4$ \\ 
77 & 948 & 3/3&8/1 & 76/7& 48/4&24/3 &32/4 & $1^4 \times S_4$ \\ 
- & 950 & 3/3&8/1 & 76/7& 48/4&24/3 &32/4 & $(C_2\times C_2)\,\wr\,C_2 \times  
S_4$ \\ 
60 & 181 & 31/4 & 32/1& 48/4 & 32/1& 48/4 & &  $C_2\times \#44[192]=384[17949]$ \\ 
79 & 180 & 7/3 & 32/1& 72/5& 32/1& 48/4&  & $C_2\times \#44[192]=384[17949] $ \\ 
- & 1490 & 31/4 & 32/1& 48/4 & 32/1& 48/4 & &  $C_2\times  
\text{Aut}(Q_2\times C_2)$ \\ 
- & 1489 & 7/3 & 32/1& 72/5& 32/1& 48/4& &  $C_2\times \text{Aut}(Q_2\times  
C_2)$ \\ \hline \hline 
- & 954 & 27/4& 8/1& 52/6& 24/3& 48/4&32/4 & $1^4\times S_4$ \\ \hline 
\multicolumn{9}{|c|}{Presentation for group number 949 in SGL:}\\ 
\multicolumn{9}{|c|}{ Just make the replacement in the presentation for our  
group \#35 (SGL =\#945)}\\ 
\multicolumn{9}{|c|}{ $c^2*d^2=1 \rightarrow a*d*a^{-1}*d*b=1$} \\ \hline 
\multicolumn{9}{|c|}{Presentation for group number 950 in SGL:}\\ 
\multicolumn{9}{|c|}{ Just make the replacement in the presentation for our  
group \#77 (SGL =\#948)}\\ 
\multicolumn{9}{|c|}{ $a^c*a*b*a^{-1}=1 \rightarrow (a,c)=(b,c)=1$} \\ 
\multicolumn{9}{|c|}{ This group can also be generated by the extension  
SL$(2,3)@Q_2$} \\ 
\multicolumn{9}{|c|}{with a $C_2\times C_2$ action; see Table 4B comment 1  
at end of table.}\\ \hline 
\multicolumn{9}{|c|}{Presentation for group number 79 above:}\\ 
\multicolumn{9}{|c|}{ Just make the replacement in the presentation for our  
group \#60}\\ 
\multicolumn{9}{|c|}{ $d^2=1 \rightarrow d^4= c*a^{-1}*c*b^{-1}=1$} \\  
\hline
\multicolumn{9}{c}{}\\ 
\end{tabular} 
\linebreak\\ 
\begin{tabular}{|c|c|c|c|}\hline 
\multicolumn{4}{|c|}{Table 9 } \\ \hline 
\multicolumn{4}{|c|}{Table for images arising in order $64p^2$ groups} \\ \hline 
\multicolumn{4}{|c|}{Case of $p=3$ (See Table A6 for a listing of the individual groups.)} \\ \hline 
image &SGL \# & Hall-Senior \# & Hall-Senior Breakdown \\ \hline 
$\Gamma _3a_2\simeq QD_8$ & 36 & 36 & $\Gamma _3[9]$, $\Gamma _8[3]$, $\Gamma _{14}[5]$, 
$\Gamma _{15}[8]$, $\Gamma _{17}[6]$, $\Gamma _{20}[2]$, $\Gamma _{21}[3]$\\ \hline 
\multicolumn{4}{|c|}{Case of $p=5$ (See Table A5 for a listing of the individual groups.)} \\ \hline 
image &SGL \# & Hall-Senior \# & Hall-Senior Breakdown \\ \hline 
$C_4\times C_2$ & 1068 & 757 & $\Gamma_1 [39]$, $\Gamma _2[222]$, $\Gamma _3[107]$, 
$\Gamma _4[163]$, \\
\hline   
\end{tabular} 
\linebreak\\ 
  
\section{Tables A1 to A6} 
\begin{tabular}{|c|c|c|c|c|}\hline 
\multicolumn{5}{|c|}{ Table A1 }\\ \hline 
\multicolumn{5}{|c|}{ Table of Automorphism Sequences for Groups of order 64  
}\\ \hline 
  number  &  Aut($g$)& Aut$^2(g)$ & Aut$^3(g)$ & Aut$^4(g)$\\ \hline 
     1    &   GL(2,6)  &             &              &                \\ 
     2    & 10,321,920 &             &              &                \\ 
     3    &  147,456(2)& \multicolumn{3}{c|}{not found ! runs out of memory}  
   \\ 
     4ag  & $C_2 \times$  21,504 &   21,504     &    43,008     & 86,016 $\ddagger$  \\ 
     5ag  & $C_2 \times $ 43,008 &  86,016 complete            &       &        
   \\
\hline 
\end{tabular} 
\linebreak\\ 
\begin{tabular}{|c|c|c|c|c|}\hline 
\multicolumn{5}{|c|}{ Table A1 (continued)}\\ \hline 
\multicolumn{5}{|c|}{ Table of Automorphism Sequences for Groups of order 64  
}\\ \hline 
  number  &  Aut($g$)& Aut$^2(g)$ & Aut$^3(g)$ & Aut$^4(g)$\\ \hline 
   150    & \multicolumn{2}{l|}{$C_4 \times  C_2 \times $ ([192] \# 3 in  
Table 2 above)}    &  &    \\ 
   153d   &  10,752    &   64,512  complete         &         &       \\ 
    155ik  & 96[table 3b \# 36] $\times $ 32[33] &     &       &         \\ 
   156k   & 96[table 3b \# 36] $\wr$ $C_2$ &  110,592 complete  & & \\ 
158ik  & 96[table 3b \# 36] $\times $ 32[33]  && &        \\ 
   159k   & 96[table 3b \# 36] $\times  1^4$   &&   &       \\ 
   162ik  & 96[table 3b \# 36] $\times $ 32[33]       &       &    &   \\ 
   173j   &    6144(1) &   49,152    &  786,432(1)  & $2^{19}*3$               \\ 
   181    &    1536    &    6,144    &   \multicolumn{2}{c|}{12,288  
complete}  \\ 
   183    &    9216    &   18,432    & \multicolumn{2}{c|}{36,864 complete}   
\\ 
   184c   &    1536    &\multicolumn{2}{c|}{258048 = $C_2 \times  C_2 \times  
$ 64,512}   &                \\ 
   187    &  15,360  complete        &              &       &         \\ 
   241    & $S_4 \times  \text{Hol}(C_8)$             &          &    &              
    \\ 
   242    & $C_2 \times  D_4 \times  S_4$             &        &      &       
           \\ 
   243    & $S_4 \times  \text{Hol}(C_8)$             &              &       
&          \\ 
   259    & $C_2 \times  \text{Aut}(Q_2 \times  C_2)$        &        &       
&                \\ 
   260    & $C_2 \times  \text{Aut}(Q_2 \times  C_2)$        &        &       
&    \\ \hline 
\end{tabular} 
\linebreak\\ 
\begin{tabular}{|c|l|} \hline 
\multicolumn{2}{|c|}{ Notes for Table A1   } \\ \hline 
   a.& The complete groups of order 86,016 are not isomorphic. The one for \\ 
      & \#4 is the same as for order 32 \#2.   \\ 
    b.& These two sequences of groups are isomorphic. \\ 
    c.& These order 1536 groups are all isomorphic. \\ 
    d.& This complete group of order 64,512 is different from that in the \\ 
      & other sequences (c type).                           \\ 
    e.& This automorphism group is isomorphic to  Hol$(C_8 \times  C_2  
\times  C_2)$.  \\
\hline 
\end{tabular}
\linebreak\\ 
\begin{tabular}{|c|c|c|}\hline 
\multicolumn{3}{|c|}{ Table A2 June 21, 1994} \\ \hline 
\multicolumn{3}{|c|}{ Groups coming from the group  GL$(2,3)$.} \\ \hline 
\multicolumn{3}{|c|}{Note that most of the following cases are duplicates of   }  
\\ 
\multicolumn{3}{|c|}{previous listings as is noted below. These alternate     
      } \\ 
\multicolumn{3}{|c|}{ presentations however might be of interest to the  
reader.    } \\ 
\multicolumn{3}{|c|}{Relations for GL$(2,3)$:} \\ 
\multicolumn{3}{|c|}{$b*a*b*((a*b*a)^{-1})=(b^2*(a^{-1}))^2=1$ }\\ \hline 
\multicolumn{3}{|c|}{} \\
\hline 
\end{tabular} 
\linebreak\\ 
\begin{tabular}{|c|c|c|}\hline 
\multicolumn{3}{|c|}{ Table A2 (continued)     June 21, 1994} \\ \hline 
\multicolumn{3}{|c|}{G. Groups coming from the group GL$(2,3)$.} \\ \hline 
no     &           GL$(2,3) \times 1^2$                        &\\ 
new   &         $ [\text{GL}(2,3) @  C_2] \times  C_2$  action [5]  &$     S_4  
\times $ \#33 order 32 \\ 
ones  &             same as  \# 47                                   &\\  
\cline{2-3} 
  in &         $C_2$ actions b. From  $C_4$ case.                             
  &\\ \cline{2-3} 
this  &          GL$(2,3) @  C_4$    action [1]     &$       S_4 \times 1^4  
$\\ 
  list &             same as  \# 37                                   &\\
\hline 
\end{tabular} 
\linebreak\\ 
\begin{tabular}{|c|c|c|}\hline 
\multicolumn{3}{|c|}{ Table A3      June 21, 1994} \\ \hline 
\multicolumn{3}{|c|}{G. Groups coming from the group  $<2,3,4>$.} \\ \hline 
\multicolumn{3}{|c|}{Note that most of the following cases are duplicates of   }  
\\ 
\multicolumn{3}{|c|}{previous listings as is noted below. These alternate     
      } \\ 
\multicolumn{3}{|c|}{ presentations however might be of interest to the  
reader.    } \\ 
\multicolumn{3}{|c|}{ } \\ 
\multicolumn{3}{|c|}{ Relations for $<2,3,4>$:  } \\ 
\multicolumn{3}{|c|}{$b*a*b((a*b*a)^{-1})=b*a^2*b*(a^{-2})=1$ } \\
\hline 
\end{tabular} 
\linebreak\\ 
\begin{tabular}{|c|c|c|}\hline 
\multicolumn{3}{|c|}{ Table A3 (continued)     June 21, 1994} \\ \hline 
\multicolumn{3}{|c|}{G. Groups coming from the group  $<2,3,4>$.} \\ \hline 
     & $             [<2,3,4> @ C_2] \times C_2$  action [2] &$  S_4 \times  
\text{Hol}(C_4 \times C_2)$\\ 
       & $              C_2 \times $  Table 3b in \cite{7} \# 34 & \\ 
    &              same as \#68 & \\ 
       & $             [<2,3,4> @ C_2 \times C_2$  action [3] &$  S_4  
\times$ \#33 order 32 \\ 
       & $              C_2 \times$ Table 3b in \cite{7} \# 35          &\\ 
    &              same as \# 66              &                   \\
\hline 
\end{tabular} 
\linebreak\\ 
\begin{tabular}{|c|c|c|c|c|c|c|c|c|c|}  \hline 
\multicolumn{10}{|c|}{Table A4           } \\ 
\multicolumn{10}{|c|}{Class Structure for Groups of Order 192} \\ 
\multicolumn{10}{|c|}{ Nonnormal Sylow Group Types          }\\ \hline 
       &  group \#  & \multicolumn{8}{|c|}{  order structure of the group }   
\\ 
   ncl &  Tables   & \multicolumn{8}{|c|}{} \\ \cline{3-10} 
       &  4 / SGL  &  2  &   3  &   4   &   6  &   8  &  12  & 16  &  24 \\  
\hline 
   11  &    55/184     &19/3 & 32/1 & 60/3  & 32/1 & 48/2 &      &     &      
\\ 
   11  &    45/185     &19/3 & 32/1 & 108/5 & 32/1 &      &      &     &      
\\
\hline 
\end{tabular} 
\linebreak\\ 
\begin{tabular}{|c|c|c|c|c|c|c|c|c|c|}  \hline 
\multicolumn{10}{|c|}{Table A4 (continued)} \\ 
\multicolumn{10}{|c|}{Class Structure for Groups of Order 192    } \\ 
\multicolumn{10}{|c|}{ Nonnormal Sylow Group Types          }\\ \hline 
       &  group \#  & \multicolumn{8}{|c|}{  order structure of the group }   
\\ 
   ncl &  Tables   & \multicolumn{8}{|c|}{} \\ \cline{3-10} 
       &  4 / SGL  &  2  &   3  &   4   &   6  &   8  &  12  & 16  &  24 \\  
\hline 
   23x &   35/945      & 3/3 &  8/1 & 76/7  & 24/3 & 48/4 & 32/4 &     &  
     \\ 
   23x &   77/948 & 3/3 &  8/1 & 76/7  & 24/3 & 48/4 & 32/4 &     &      
\\
\hline 
\end{tabular} 
\linebreak\\ 
 
\begin{tabular}{|c|c|c|c|c|c|c|c|c|c|}  \hline 
\multicolumn{10}{|c|}{Table A4 (continued)} \\ 
\multicolumn{10}{|c|}{Class Structure for Groups of Order 192} \\ 
\multicolumn{10}{|c|}{ Nonnormal Sylow Group Types }\\ \hline 
       &  group \#  & \multicolumn{8}{|c|}{  order structure of the group }   
\\ 
   ncl &  Tables   & \multicolumn{8}{|c|}{} \\ \cline{3-10} 
       &  4 / SGL  &  2  &   3  &   4   &   6  &   8  &  12  & 16  &  24 \\  
\hline 
     &           &     &      &       &      &      &      &     &     \\ 
   32  &   29/183   & 3/3 &  8/1 & 28/8  & 24/3 & 96/12& 32/4 &     &     \\
\hline 
\end{tabular} 
\linebreak\\ 
\begin{tabular}{|c|}  \hline 
Notes for Table A4.\\ 
    ncl means the number of conjugacy classes in the group.   \\ 
            In the first row above, the number 19/3 means there are       \\ 
                 19 elements of order 2 distributed in 3 classes,       \\ 
               and likewise for the other groups and element orders.     \\  
\hline 
    Notes:    x in column 1 means groups having the same class/order structure but \\ 
                 having different automorphism groups; see Tables 4 and 5.      \\ 
                                                                   \\ 
             23M is a missing group in my list, GAP number given.            \\ 
                                                                         \\ \hline 
                                                                         \\ 
    The groups of the form $A_4$ @ (group of order 16) are not explicitly    
\\ 
        written out in Table 4A, but the 2-group along with the generator      \\ 
       acting on $A_4$ as an element of order 2 is explicitly written out.  \\ 
       Thus one can reconstruct the group of order 192.                   \\  
\hline 
                                                                         \\ 
          * The relations used here for $(4,4|2,2)$ were:                   
\\  
                                                                         \\ 
             $a^2=b^4=(b,a)*(b^{-1})*a*(b^{-1})*a*b^2=$                       
   \\ 
                     $(b^{-1})*a*b*a*(b^{-1})*a*b*a=1   $                     
   \\ \hline 
                                                                         \\ 
         ** The relations used for $C_4 $Y$ Q_2$ were:                          
\\  
                                                                         \\ 
                   $c^4=d^2=e^2=(d,e)*c^2=(c,d)=(c,e)=1 $                   
\\ \hline 
\end{tabular} 
\linebreak\\ 
 
\begin{tabular}{|c|c|l|} \hline 
\multicolumn{3}{|c|}{Table A5. Details on Quotient Groups in Table 9, $p=5$ case} \\ \hline 
image & Isoclinic Class& \qquad \qquad \qquad groups \\ \hline 
$C_4\times C_2$ & $\Gamma _1$ & $21^4$ [2], $2^21^2$ [4], $31^3$ [5], $2^3$ [2], 321 [11]\\ 
 757& &  411 [5], $3^2$ [2], 42 [5], 51 [3]\\ \cline{2-3} 
 & $\Gamma _2$ & $c_1$[4], $c_2$[4], $d$[13], $e_1$[6], $e_2$[4], $f$[8], $g$[13],\\  
 & & $h$[4], $i$[11], $j_1$[11], $j_2$[11], $k$[11], $l$[6], $m_1$[10], \\ 
 & & $m_2$[10], $n$[4], $o$[14], $p$[17],$q$[18],$r_1$[11],$r_2$[5],\\ 
 & & $s$[5], $t$[4], $u$[4], $v$[5], $w_1$[3], $w_2$[3], $x$[3] \\
\hline  
\end{tabular} 
\linebreak\\ 
 
\begin{tabular}{|c|c|l|} \hline 
\multicolumn{3}{|c|}{Table A6. Details on Quotient Groups in Table 9, $p=3$ and 7 cases} \\ \hline 
image & Isoclinic Class& \qquad \qquad \qquad groups \\ \hline 
\multicolumn{3}{|c|}{$p=3$ case}\\ \hline 
$QD_8=\Gamma _3a_2$ & $\Gamma _3$ & $a_2,c_1,c_2,d_2,i_2,k,n_1,n_2,o_2$\\ 
                    & $\Gamma _8$ & $c_1,c_2,e$\\ 
                    & $\Gamma _{14}$& $a_2,a_4,b_1,b_2,c_2$ \\ 
                    & $\Gamma _{15}$ & $ a_2,a_3,b_2,b_4,b_5,b_6,c_2,c_3$ \\ 
                    & $\Gamma_{17}$& $a_1, a_2(\text{ twice }), a_3,b_1,b_2$\\ 
                    & $\Gamma _{20}$ & $a$ (twice)  \\ 
                    & $\Gamma _{21}$& $a_1,a_2,a_3$\\ \hline 
\multicolumn{3}{|c|}{$p=7$} \\ \hline 
$C_{16}$ & $\Gamma _1$ & 411, 42, 51 (twice), and 6\\ 
         & $\Gamma _2$ & $w_1,w_2, x (\text{ twice) }$ \\ \hline 
$\Gamma _3a_1=D_8$ & $\Gamma _3$ & $ a_1,c_1,d_1,i_1,k,n_1,o_1$\\ 
                   & $\Gamma _8$ & $a_1[2], a_2[2], a_3[2], b[2], c_1,c_2,d_1,d_2,e,f$\\ 
                   & $\Gamma _{14}$ & $a_1,a_3,b_1,c_1$ \\ 
                   & $\Gamma _{15}$ & $a_1,a_3,b_1,b_3,e_1 ?$ \\ 
                   & $\Gamma_{17}$ & $a_1,a_2,b_1$ \\ 
                   & $\Gamma _{19}$ &$ a_1[2], a_2[2]$ \\ 
                   & $\Gamma _{21}$ &$ a_1,a_2,a_3$ \\ 
                   & $\Gamma _{27}$ &$ a_1,a_2,a_3$ \\ \hline 
$\Gamma _3a_3=Q_4$& $\Gamma _3$ & $ a_3,c_2,d_1,i_3,k,n_2,o_1$\\ 
                   &$\Gamma _8$ & $d_1,d_2,f$ \\ 
                   &$\Gamma _{14}$ & $ a_3,a_5,b_2,c_3$\\ 
                   & $\Gamma _{15}$ & $a_4,b_7,b_8,c_4$\\ 
                   & $\Gamma _{17}$ & $ a_2,a_3,b_2$\\ \hline  
\multicolumn{3}{c}{} \\
\end{tabular}

\begin{tabular}{|l|c|c|c|c|}\hline 
\multicolumn{4}{|c|}{Table 7} \\ \hline 
\multicolumn{4}{|c|}{ Nonsplit extensions arising from SL$(2,3)@D_4$ } \\  
\hline 
    action  & semi-direct case    & \multicolumn{2}{|c|}{nonsplit extensions} \\ 
      &         & $c^8=c^4(ab^{-1})^2$    & $c^4=d^4=d^2(ab^{-1})^2$ \\ \hline     
$C_2[C_4]$ & 33 & 64 &  72 \\  
$C_2[C_2]$ & 34 & 50 & 74 \\ 
$(C_2\times C_2)(c,d)$ & [33] & [59] & 949 M \\ 
$(C_2\times C_2)(d,c)$ & 40  & [41] & 73\\ 
$D_4$ & 41 & [40] & 76 \\ \hline                   
\end{tabular}
\newpage

\section{Appendix II \\
Structure and relations for selected \\automorphism groups in Table A1}

\subsection{Group \#3. [order 147,456]}
		The group has the structure:
\begin{equation*}
			[ (4096) @ (C_3 \times  C_3) ] @ (C_2 \times  C_2).
\end{equation*}
There are 13 possible ($C_3 \times  C_3$) actions on the 4096 group. CAYLEY generated the following two possible choices for a presentation of this group of order 36,864:\\
\begin{gather*}
 T.1^2=T.2^2=T.3^2=T.4^2=T.5^2=T.6^2=T.7^2=T.8^2=(T.1*T.2)^2=\\
(T.1*T.4)^2=(T.1*T.6)^2=(T.2*T.3)^2=(T.2*T.4)^2=(T.2*T.6)^2=\\
(T.3*T.5)^2=(T.3*T.7)^2=(T.3*T.8)^2=(T.4*T.6)^2=(T.5*T.7)^2=\\
(T.5*T.8)^2=(T.6*T.7)^2=(T.7*T.8)^2=(T.1*T.2*T.5)^2=\\
T.1*T.3*T.5*T.1*T.5*T.3=(T.4*T.6*T.8)^2=\\
T.4*T.7*T.8*T.4*T.8*T.7=\\
T.1*T.2*T.7*T.8*T.2*T.7*T.1*T.8=(T.1*T.3*T.4*T.3)^2=\\
T.1*T.3*T.4*T.3*T.4*T.8*T.1*T.8=(T.1*T.3*T.6*T.3)^2=\\
T.1*T.4*T.5*T.4*T.5*T.7*T.1*T.7=\\
T.2*T.3*T.6*T.3*T.6*T.7*T.2*T.7=\\
T.2*T.4*T.5*T.8*T.2*T.4*T.8*T.5=\\
T.3*T.4*T.5*T.4*T.6*T.5*T.3*T.6=1
\end{gather*}
with the actions:
\begin{gather*}
T.9^3=T.10^3=(T.9,T.10)=\\
T.1^{T.9}*((T.4*T.7*T.4*T.8)^{-1})=\\
T.2^{T.9}*((T.1*T.2*T.8*T.1*T.2)^{-1})=\\
T.3^{T.9}*((T.1*T.4*T.6*T.8*T.1)^{-1})=\\
T.4^{T.9}*((T.1*T.3*T.5*T.7*T.8*T.1)^{-1})=\\
T.5^{T.9}*((T.4*T.5*T.6*T.5*T.7*T.4)^{-1})=\\
T.6^{T.9}*((T.4*T.7*T.4*T.5)^{-1})=\\
T.7^{T.9}*((T.1*T.7*T.2*T.5*T.6)^{-1})=\\
T.8^{T.9}*((T.2*T.3*T.4*T.6*T.8)^{-1})=1\\
\end{gather*}
and
\begin{gather*}
T.1^{T.10}*((T.1*T.2*T.4*T.7*T.2*T.7)^{-1})=\\
T.2^{T.10}*((T.4*T.5*T.2*T.6*T.5)^{-1})=\\
T.3^{T.10}*((T.1*T.2*T.3*T.8*T.1*T.2)^{-1})=\\
T.4^{T.10}*((T.3*T.1*T.4*T.3*T.4)^{-1})=\\
T.5^{T.10}*((T.1*T.2*T.7*T.1*T.5*T.2)^{-1})=\\
T.6^{T.10}*((T.1*T.5*T.2*T.5)^{-1})=\\
T.7^{T.10}*((T.4*T.5*T.4)^{-1})=\\
T.8^{T.10}*((T.1*T.2*T.3*T.1*T.8*T.2*T.8)^{-1})=1\\
\end{gather*}
or, alternatively, with this set of actions:
\begin{gather*}
T.9^3=T.10^3=(T.9,T.10)=\\
T.1^{T.9}*((T.4*T.3*T.4*T.6*T.5*T.6)^{-1})=\\
T.2^{T.9}*((T.4*T.3*T.4)^{-1})=\\
T.3^{T.9}*((T.3*T.8*T.2*T.8)^{-1})=\\
T.4^{T.9}*((T.1*T.4*T.8*T.1*T.4*T.7)^{-1})=\\
T.5^{T.9}*((T.1*T.7*T.2*T.5*T.7)^{-1})=\\
T.6^{T.9}*((T.2*T.7*T.2)^{-1})=\\
T.7^{T.9}*((T.2*T.6*T.7*T.2)^{-1})=\\
T.8^{T.9}*((T.3*T.6*T.3*T.8*T.4)^{-1})=1
\end{gather*}
and
\begin{gather*}
T.1^{T.10}*((T.1*T.5*T.2*T.6*T.5*T.6)^{-1})=\\
T.2^{T.10}*((T.3*T.7*T.1*T.3*T.7)^{-1})=\\
T.3^{T.10}*((T.3*T.4*T.7*T.4*T.5*T.7)^{-1})=\\
T.4^{T.10}*((T.3*T.6*T.3)^{-1})=\\
T.5^{T.10}*((T.3)^{-1})=\\
T.6^{T.10}*((T.1*T.5*T.1*T.4*T.5*T.6)^{-1})=\\
T.7^{T.10}*((T.1*T.8*T.1)^{-1})=\\
T.8^{T.10}*((T.1*T.4*T.7*T.4*T.8*T.1)^{-1})=1.
\end{gather*}		     
\begin{tabular}{|l|c|c|c|c|c|c|}\hline
\multicolumn{7}{|c|}{order structure of $[ (4096) @ (C_3 \times  C_3) ]$ group } \\ \hline
 order of elements    & 2 & 3  &   4   &  6   &   12   &   order of \\
 number of elements   &     927  & 2816 & 3168 & 20736 & 9216  &  center     \\
 number of classes    & 25  &  8   &  32   & 40   &  12  & = 4   \\ \hline
\end{tabular}
\linebreak\\

The automorphism group of this group of order 36,864 has order 221,184 = $2^{13}*3^3$.\\

CAYLEY gives the following set of relations on three generators for the automorphism group (of order 147,456):
\begin{gather*}
             a^6=b^6=c^2=(a*(b^{-1}))^3=a^2*b*(a^{-1})*(b^{-2})*(a^{-1})*b=\\
             (a^2*(b^{-2}))^2=(a*(c^{-1}))^4=(b*(c^{-1}))^4=\\
             a*b*(a^{-1})*(b^{-1})*c*(b^{-1})*(a^{-1})*b*a*c=(a^2*(c^{-1}))^4=\\
             a^2*(c^{-1})*a*(b^{-1})*(c^{-1})*(b^{-2})*(c^{-1})*(b^{-1})*a*(c^{-1})=\\
             a*b*(c^{-1})*b*(a^{-1})*(c^{-1})*(a^{-1})*(b^{-1})*(c^{-1})*(b^{-1})*a*(c^{-1})=\\
             (a*(b^{-1})*c*b*(a^{-1})*c)^2=(a*c*(a^{-1})*b*c*(b^{-1}))^2=\\
             (a*c*(a^{-1})*c)^3=(b*c*(b^{-1})*c)^3=\\
             a*b^2*c*b^2*c*(a^{-1})*c*b*(a^{-1})*c*(a^{-1})*b=\\
             a^2*c*(b^{-2})*(a^{-1})*c*a*c*a*(b^{-1})*(a^{-1})*b*c*(a^{-1})*b=1.
\end{gather*}
            This automorphism group has only 19 normal subgroups.\\
The automorphism group of this order 147,456 group has order 294,912. It has 119 normal subgroups and has a center of $C_2\times C_2$. GAP fails to get the automorphism group of this group of order 294,912.

\subsection{Group \#4. [order 22,008]} The automorphism group for the group $C_8\times C_2\times C_2\times C_2$ is $C_2\times [21,504]$. The presentation of the order 21,504 group is:
\begin{gather*}
a^2=b^2=c^2=d^2=e^2=f^2=h^2=(a,b)=(a,c)=(b,e)=(b,h)\\
=(c,d)=(c,f)=(d,f)=(e,h)=\\
a*b*c*b*c=a*b*e*a*e=a*b*h*a*h\\
=a*c*d*a*d=a*c*f*a*f\\
=(b*d)^3=(b*f)^3=(c*e)^3=(c*h)^3=(d*e)^4\\
=(d*e*d*h)^2=(d*e*f*e)^2=(e*f*h*f)^2=1.
\end{gather*}

\subsection{Group \#5. [order 86,016]} The automorphism group of $C_4 \times  C_4 \times  C_4$ is $C_2\times$ [43008]. A presentation for the {43008} factor here is:
\begin{gather*}
a^2=b^2=c^2=d^2=(a*b)^2=(a*c)^3=(a*d)^3=\\
(b*c)^3=(c*d)^4=(a*b*d*c*d)^2=  \\
a*b*d*a*b*d*a*b*d*b*a*d= \\
a*c*a*d*a*c*a*d*a*c*d*b*d*c*a*d=\\
a*b*c*a*d*a*c*d*c*b*d*c*b*d*b*d*c*b*d=1.
\end{gather*}

\subsection{Group \#7. [order 768]}
The order 384 group [\# 20095 in the small group library] here has the structure:
\begin{equation*}
                  [( Q_2 \text{Y}  Q_2) @  C_3] @ \left(C_2\times C_2\right).
\end{equation*}
This is a nonsplit extension via a $C_4\times C_2$ quotient.\\
Relations for this case are:
\begin{gather*}
                 b^4=a^3*(b^{-2})=a^2*c*(a^{-1})*c=a*c^2*a*(c^{-1})=\\
                 a*b*a*(b^{-1})*a*(b^{-1})=a*(c^{-1})*b*c*(a^{-1})*(b^{-1})=\\
                 b*c*(b^{-1})*c*(b^{-1})*c=\\
                 d^4=e^2=d^2*b^2=(e,d)=\\
                 a^d*((a*c*b)^{-1})=b^d*((a*c*a)^{-1})=c^d*((b*a*c)^{-1})=\\
                 a^e*((b*c*a)^{-1})=b^e*((a*c*a)^{-1})=c^e*((c*a*b)^{-1})=1.
\end{gather*}

\subsection{Group \#8. [order 1,536]}
This group has the structure:
\begin{equation*}
                  [\text{order 192}\times C_2\times C_2] @ C_2.
\end{equation*}
The order 192 group in this representation is group number 47 in Table 2a. As of now we have not been able to reconstruct a presentation for this automorphism group using this decomposition. Cayley returns the following presentation for this automorphism group:
\begin{gather*}
a^2=b^2=c^2=d^2=(a*b)^2=(a*c*b*c)^2=\\
a*c*d*c*d*c*a*d=(a*d)^4=(a*d*b*d)^2=\\
(b*c)^4=(b*d)^4=a*b*c*a*b*c*a*c*a*c=\\
b*c*b*d*b*d*c*d*b*d=1.
\end{gather*}
\subsection{Group \#14. [order 73,728]}
The group can be written as:
\begin{equation*}
[(2048) @ ( C_3\times  C_3)] @ \left(C_2\times C_2 \right).
\end{equation*}

The relations for this group of order 2048 are
\begin{gather*}
 T.1^2=T.2^2=T.3^2=T.4^2=T.5^2=T.6^2=T.7^2=T.8^2=\\
(T.1*T.2)^2=(T.1*T.3)^2=(T.1*T.4)^2=(T.1*T.5)^2=\\
(T.1*T.6)^2=(T.1*T.7)^2=(T.2*T.3)^2=(T.2*T.4)^2=\\
 (T.2*T.5)^2=(T.2*T.6)^2=(T.2*T.8)^2=(T.3*T.4)^2=\\
(T.3*T.5)^2=(T.3*T.6)^2=(T.4*T.5)^2=(T.4*T.6)^2=\\
 (T.5*T.6)^2=(T.7*T.8)^2=(T.3*T.6*T.8)^2=\\
 (T.4*T.5*T.7)^2=T.1*T.2*T.7*T.2*T.7*T.8*T.1*T.8=\\
 T.3*T.4*T.5*T.8*T.3*T.5*T.4*T.8=\\
 T.3*T.4*T.6*T.7*T.3*T.6*T.4*T.7=\\
 T.3*T.4*T.7*T.3*T.8*T.4*T.8*T.7=\\
 T.3*T.4*T.7*T.4*T.7*T.8*T.3*T.8=1.
\end{gather*}
The automorphism group of this 2-group of order 2048 has order $2^{31}*3^2*7$ = 135,291,469,824. The 3-group actions on the 2-group are:
\begin{gather*}
T.9^3=T.10^3=(T.9,T.10)=\\
T.1^{T.9}*((T.1*T.3*T.4*T.8*T.1*T.5*T.8)^{-1})=\\
T.2^{T.9}*((T.1*T.4*T.8*T.1*T.6*T.8)^{-1})=\\
T.3^{T.9}*((T.2*T.3*T.4*T.5*T.7*T.6*T.7)^{-1})=\\
T.4^{T.9}*((T.3*T.8*T.1*T.6*T.8)^{-1})=\\
T.5^{T.9}*((T.1*T.2*T.3*T.6)^{-1})=\\
T.6^{T.9}*((T.1*T.3*T.4*T.5*T.6)^{-1})=\\
T.7^{T.9}*((T.3*T.8*T.3)^{-1})=\\
T.8^{T.9}*((T.1*T.3*T.8*T.1*T.3*T.7)^{-1})=\\
T.1^{T.10}*((T.6*T.7*T.3*T.7)^{-1})=\\
T.2^{T.10}*((T.1*T.3*T.4*T.8*T.1*T.5*T.6*T.8)^{-1})=\\
T.3^{T.10}*((T.1*T.2*T.3*T.5*T.8*T.4*T.8)^{-1})=\\
T.4^{T.10}*((T.3*T.5*T.6*T.7*T.2*T.7)^{-1})=\\
T.5^{T.10}*((T.2*T.4*T.6)^{-1})=\\
T.6^{T.10}*((T.1*T.3*T.7*T.2*T.5*T.6*T.7)^{-1})=\\
T.7^{T.10}*((T.1*T.4*T.8*T.1*T.4)^{-1})=\\
T.8^{T.10}*((T.7*T.8)^{-1})=1.
\end{gather*}
A possible alternate representation of this automorphism group is:
\begin{equation*}
[1^6 @ C_3 \times  1^3] @ (C_2 \times  S_4), 
\end{equation*}
where the order 192 group is (\#3 in Table 2a). A set of relations for this form has not been obtained.\\

A set of relations for this automorphism group generated by CAYLEY is:
\begin{gather*}
a^2=b^2=c^2=d^2=e^2=f^4=(a,b)=(a,d)=(a,e)=(a,f)=\\
(b,d)=(b,e)=(b,f)=(c,e)=(c,f)=(d,f)=(a*c)^3=\\
(e*f)^3=(b*c)^4=(c*d)^4=(d*e)^4=(a*b*c)^2*b*a*c=\\
d*e*f^2*e*d*e*(f^{-2})*e=a*b*c*d*c*a*d*b*c*d*c=\\
c*(d*c*e)^2*d*e*c*d*e=\\
b*c*d*e*d*c*b*c*d*(f^{-1})*e*f*e*(f^{-1})*d*c=1.
\end{gather*}

\subsection{Group \#15. [order 49,152]}
This group's structure takes the form:
\begin{equation*}
[8192\,\, @\,\,  C_3] \,\,@ \,\,C_2,
\end{equation*}
where the order 8192 relations are:
\begin{gather*}
a^2=b^2=c^2=d^2=e^2=f^2=\\
(a*b)^2=(a*e)^2=(a*f)^2=(b*d)^2=\\
 (b*e)^2=(c*d)^2=(c*e)^2=a*b*e*f*e*b*f=\\
 (a*c)^4=(a*c*a*d)^2=(a*c*b*c)^2=\\
a*c*f*d*a*d*f*c=(b*c)^4=(b*f)^4=\\
 (b*f*d*f)^2=(c*f)^4=(c*f*d*f)^2=\\
 (c*f*e*f)^2=(d*e)^4=\\
(d*f)^4=b*c*b*c*f*d*e*d*e*f=1.
\end{gather*}
The order 8192 group can be generated with six generators, has 254 classes, center $C_2\times C_2$ and an automorphism group of order $2^{26} * 3 = 201,326,592$.\\

The action of the $C_3$ on this 2-group is
\begin{gather*}
h^3=\\
a^h*((c*a*c*d*e*f*e*f*d)^{-1})=\\
b^h*((f*e*f)^{-1})=\\
c^h*((b*c*b*f*d*f)^{-1})=\\
d^h*((a*f*c*a*f)^{-1})=\\
e^h*((c*d*e*d*f*b*f*c)^{-1})=\\    
(f,h)=1,
\end{gather*}
which will yield the normal subgroup of order 24,576. This automorphism group does not have a normal sylow 2-subgroup. A CAYLEY generated set of relations for this automorphism group of order 49,152 is:
\begin{gather*}
a^6=b^6=c^2=(a*(b^{-1}))^2=a^2*(b^{-2})*a*b*(a^-1)*(b^{-1})=\\
(a*c)^4=(b*c)^4=a*(b^{-2})*a*c*(a^{-1})*b*a*(b^{-1})*c=\\
a^2*b*a^2*(b^{-2})*(a^{-1})*(b^{-2})*(a^{-1})*b=\\
a^2*c*a*(b^{-1})*c*(b^{-2})*c*(a^{-1})*b*c=\\
(a^2*c)^4=a^2*c*b*a*c*b^2*c*a*b*c=\\
a*b*c*a*(b^{-1})*c*(a^{-1})*(b^{-1})*c*(a^{-1})*b*c=\\
(a*c*(a^{-1})*c)^3=(a*c*(b^{-1})*c)^3=\\
a^3*(b^{-3})*c*(b^{-2})*(a^{-1})*b*a^2*c=1.
\end{gather*}
\subsection{Group \#16. [order 49,152]}
This group can be represented as:
\begin{equation*}
(8192 @  C_3) @ C_2 = 8192 @ S_3
\end{equation*}
or
\begin{equation*}
(4096) @  C_3) @\left(C_2\times C_2\right).
\end{equation*}
\linebreak\\
Relations for the 8192 group and its extensions by  $C_3$:\\

Sylow 2-factor [order 8192 case]:
\begin{gather*}
a^2=b^2=c^2=d^2=e^2=f^2=h^2=\\
(a,b)=(a,c)=(b,d)=(c,d)=(d,e)=(f,h)=\\
(a*b*e)^2=(a*b*f)^2=(a*c*h)^2=(b*c*e)^2=\\
(c*d*f)^2=(d*e*h)^2=\\
a*b*c*b*c*d*a*d=(a*e)^4=a*e*f*h*e*a*h*f=\\
(a*f)^4=(a*f*a*h)^2=(a*f*e*f)^2=\\
(a*h)^4=b*d*f*h*b*d*h*f=(b*h*c*h)^2=1.
\end{gather*}

Aut(of 8192 group) has order $2^{27}*3^2 =1,207,959,552$.\\ The automorphism group of this group of order 8192 gives rise to the following four actions of order 3 on this 2-group. 
The group appearing as a normal subgroup of Aut(64\# 16) is the above 2-group with the $C_3$ action in case four below.\\

Case 1:\\
\begin{gather*}
k^3=\\
a^k*((c*f*d*e*b*f*e)^{-1})=\\
 b^k*((b*c*e*c*e)^{-1})=\\
    c^k*((a*b*f*c*e*a*f)^{-1})=\\
     d^k*((a*f*e*f*d)^{-1})=\\
    e^k*((f*e*h*c*h*e*f)^{-1})=\\
     f^k*((a*c*f*a*c)^{-1})=\\
     (k,h)=1. 
\end{gather*}
This group can also be generated by $<a,b,d,e,h,k>$. It has 161 classes and center = I. The automorphism group of this group has order $2^{17} * 3^2 = 1,179,648$.\\

Case 2:
\begin{gather*}
     			a^k*((b*f*d*f*b)^{-1})=\\
     			b^k*((h*e*b*e*h)^{-1})=\\
     			c^k*((f*b*c*b*f)^{-1})=\\
     			d^k*((a*c*h*e*b*h*d)^{-1})=(e,k)=\\
     			f^k*((a*d*h*d*a*f)^{-1})=\\
     			h^k*((c*e*f*e*c)^{-1})=1.
\end{gather*}
This group can also be generated by $<a,b,c,f,h,k>$. This group has 145 classes and center = I.
The automorphism group of this group has order $2^{20} * 3^2 = 9,437,184$.\\
Case 3:
 \begin{gather*}
a^k*((a*b*c*e*f*e*f)^{-1})=\\ 
b^k*((c*f*h*b*f*h*c)^{-1})=\\
c^k*((b*f*e*f*h*c*h)^{-1})=\\
d^k*((d*h*b*f*h*c*f)^{-1})=\\
e^k*((f*e*h*c*h*e*f)^{-1})=\\
f^k*((b*d*h*b*d)^{-1})=\\
h^k*((a*b*f*h*a*b)^{-1})=1.
\end{gather*}
This group can also be generated by $<a,b,d,f,h,k>$. This group has 145 classes and center = I. The automorphism group of this group has order $2^{20} * 3^2 = 9,437,184$. This is the same as Case 2.\\

Case 4: This is the group action that will give the order 8192*3 group appearing as a normal subgroup of Aut(64 \# 16).

\begin{gather*}
k^3=\\
a^k*((a*b*f*d*f*a*e)^{-1})=\\
b^k*((a*f*h*b*f*h*a)^{-1})=\\
c^k*((d*f*e*f*d)^{-1})=\\
d^k*((a*b*d*f*h*c*f*h*b)^{-1})=\\
e^k*((h*b*e*c*h)^{-1})=\\
f*k*((a*b*h*a*b)^{-1})=\\
h^k*((a*b*f*h*a*b)^{-1})=1.
\end{gather*}
This group can also be generated by $<a,b,c,d,e,f,k>$. This group has 161 classes and a center of order 4.
The automorphism group of this group has order $2^{13} * 3^2 = 73,728$.\\

This group does not have a normal sylow 2-subgroup.\\

A presentation for the automorphism group of order 49,152 is:
\begin{gather*}
a^6=b^6=d^2=(a*(b^{-1}))^3=(a*(c^{-1}))^3=(b*(c^{-1}))^2=\\
(b^2*c)^2=(b*c^2)^2=a^2*b*a*b^2*a*b=\\
a^2*b*((a*b^2*a)^{-1})*b=a^2*(b^{-1})*c*(a^{-2})*c*(b^{-1})=\\
(a*b)^2*((b*a)^{-1})^2=a*b*((a*c)^{-1})*c^2*(b^{-2})=\\
a*(b^{-1})*a*d*(b^{-1})*a*(b^{-1})*d=\\
a*(b^{-1})*d*c*(a^{-1})*b*d*(c^{-1})=\\
(a*d*(b^{-1})*d)^2=(a*d)^4=(b*d)^4=\\
a^2*b*(a^{-1})*d*(a^{-1})*b*a^2*d=\\
a*b*((b*a)^{-1})*d*((a*b)^{-1})*b*a*d=1.
\end{gather*}
\subsection{Group \#17. [order 6,144]}
This group has the structure:
\begin{equation*}
[(256 \,\,@ \,\, C_3)\,\, @ \,\,\left(C_2\times C_2\right)].
\end{equation*}
Four possible cases of a  $C_3$ action on the order 256 group were found by CAYLEY. Only one generates the correct order 768 group.\\

\underline{2-group relations:}
\begin{gather*}
a^2=b^2=c^2=d^2=e^2=f^2=(a,b)=(a,c)=(a,d)=(b,c)=(b,d)=\\
(b,e)=(c,d)=(c,f)=(e,f)=(a*b*f)^2=a*e*f*a*f*e=(c*d*e)^2=\\
d*e*f*d*f*e=1.
\end{gather*}
This 2-group is number 55805 of order 256 in the Small Group Library of Besche and Eick.
This 2-group has an automorphism group of order $589,824=2^{16}*3^2$.\\

\underline{$C_3$ actions:}

\underline{Case 1.} In this case the center equals I, and this action does not yield the correct group. According to GAP this is group number 1085112 of order 768 in the Small Group Library. This group's automorphism group has order 18432.
\begin{align*}
&a^h*e*a*e*d*c*a=\\
&b^h*e*c*e=\\
&c^h*e*a*e*d=\\
&d^h*e*d*e*c*b=\\
&e^h*f=f^h*f*e=1.
\end{align*} 

\underline{Case 2.}  This action generates a group that looks like case one, above, and is identified as [768 \# 1085110] in the Small Group Library. This group's automorphism group has order 36864. The automorphism group arising in this case is a complete group with 56 classes. This complete group may be the same as the one arising from Aut[64 \# 183) @ $C_3$] below.
\begin{align*}
&a^h*e*a*d*e=\\
&b^h*e*d*e*c*b=\\
&c^h*f*d*b*f*d*a=\\
&d^h*e*c*e*c*a=\\
&e^h*c*a*e*c*a=\\
&f^h*a*f*a=1.
\end{align*}

\underline{Case 3.} According to GAP this action generates the same group as case one above, namely order 768 \#1085112.
\begin{align*}
&a^h*f*d*b*f*c=\\
&b^h*c*b*a=\\
&c^h*e*c*a*e*c=\\
&d^h*b=\\
&e^h*d*f*d=\\
&f^h*f*c*a*e*c*a=1.
\end{align*}

\underline{Case 4.} This action looks like the one that generates the correct group of order 768, [768 \# 1085111]. This group's automorphism group has order 9216.\\

\begin{align*}
&h^3=\\
&a^h*e*c*e*c*b=\\
&b^h*e*c*e*c*b*a=\\
&c^h*e*d*e=\\
&d^h*e*d*e*c=\\
&e^h*a*f*a=f^h*f*c*e*c=1.
\end{align*}
The next step in the sequence: finding the $C_2\times C_2$ actions on this 768 group has not been done. The subgroup lattice for this case may be rather hard to determine. An alternate presentation for this automorphism group's factor of order 3072 is:
\begin{gather*}
 a^6=c^6=d^2=(a*(b^{-1}))^2=(a*(c^{-1}))^2=\\
 (a*(d^{-1}))^2=(b*(d^{-1}))^2=(c*(d^{-1}))^2=\\
 (a^2*b)^2=(a*b^2)^2=\\
a^2*(b^{-1})*(c^{-1})*b*(c^{-1})*(a^{-1})*b=\\
a^2*c^{-2}*a*c*(a^{-1})*(c^{-1})=\\
 a*(b^{-1})*c^{-2}*(a^{-1})*b*c^2=\\
 b*c*b*c^2*(b^{-1})*c^2=1.
\end{gather*}

\subsection{Group \#19. [order 12,288]}
This group has the structure:
\begin{equation*}
[ 1536\times C_2\times C_2] @ C_2.
\end{equation*}
The order 1536 group has as a normal subgroup:
\begin{equation*}
                  1^6 \,\,@ \,\, C_3 \text{ with 24 classes and trivial center},
\end{equation*}
and the quotient group 1536/192 is the elementary abelian group of order 8. Cayley gives the following alternate description of this automorphism group of order 12,288:
\begin{gather*}
a^2=b^2=c^2=d^4=e^4=f^2=(a,b)=(a,d)=(a,e)= \\
(b,d)=(b,e)=(c,d)=(c,e)=(c,f)=(d,f)=\\
a*d^2*f*a*f=d^2*e*(d^{-2})*e=d*e*d*(e^{-1})*(d^{-1})*(e^{-1})=\\
e^2*f*(e^{-2})*f=(a*c)^4=(a*c*b*c)^2=  \\
(b*c)^4=(b*c*b*f)^2=(b*f)^4=  \\
e*f*e*f*(e^{-1})*f*(e^{-1})*f= \\
b*e*f*b*e*f*b*(e^{-1})*f*(e^{-1})*b*f=\\
b*c*b*c*d*e*d*f*(e^{-1})*f*(d^{-1})*(e^{-1})*f=1.
\end{gather*}
\subsection{Group \#21. [order 1,536]}
This group has the structure:

\begin{equation*}
  C_2\times [(\text{order 96 group; \#36 in Table 3b of \cite{7}})\times C_2\times C_2 ] @ C_2. 
\end{equation*}
There are 5 possible actions here by the $C_2$ quotient group, three of which yield a group of the form $D_4\times 96$. The $C_2$ action on $C_2\times C_2$ is the same as in the case of a wreath product.The other two cases listed below appear to give the order 768 factor in the automorphism group for this group:
\begin{gather*}
a^4=b^4=(a*b)^3=(a*(b^{-1}))^3=(a^2,b^2)=c^2=d^2=\\
=(c*d)^4=(a,c)=(b,c)=\\
\begin{cases}
a^d*a=b^d*((a*b*a)^{-1})=1\\
\text{or}\\
a^d*a*((b*a^2*b)^{-1})=b^d*a*b*((a*b*a)^{-1})=1.
\end{cases}
\end{gather*}
The group in question here (of order 768) is [768, \# 1087581] in the Small Group Library. Its automorphism group has order 1536. The group generated by $<a,b>$ is [96, \#227] and has an automorphism group of order 576 [576, \#8654].

\subsection{Group \#22. [order 12,288]}

This automorphism group and that for group \#68 below appear to be isomorphic. This automorphism group has the structure:
\begin{equation*}
[(\text{order 256 group}\times C_2\times C_2)\,\, @ \,\, C_3]\,\, @ \,\,\left(C_2\times C_2\right).
\end{equation*}
The automorphism group of the order 3072 normal subgroup has order $2^{20} * 3^3$. \\

There are 4 possible  $C_3$ actions on the order 256 group. The automorphism group of the 256 group is $2^{19} * 3^2$. See group \#68 below for details.\\

For this automorphism group, Cayley gives the following presentation:
\begin{gather*}
 a^2=b^4=c^4=d^2=a*b*a*(b^{-1})=(a*b*d)^2=\\
 a*c^2*a*(c^{-2})=b^2*c*(b^{-2})*c=\\
 b*c^2*b*(c^{-2})=c^2*d*(c^{-2})*d=\\
 a*c*a*c*a*(c^{-1})*a*(c^{-1})=(a*d)^4=\\
 c*d*c*d*(c^{-1})*d*(c^{-1})*d=\\
 a*c*b*c*b*c*a*(c^{-1})*b*(c^{-1})*(b^{-1})*c=\\
 a*c*d*a*d*c*a*(c^{-1})*d*a*d*(c^{-1})=\\
 a*c*d*c*a*d*a*(c^{-1})*d*(c^{-1})*a*d=\\
 b*c*b*c*b*c*(b^{-1})*(c^{-1})*(b^{-1})*(c^{-1})*(b^{-1})*(c^{-1})=\\
 b*d*c*d*b*d*(c^{-1})*d*(b^{-1})*d*(c^{-1})*d=\\
 a*b*c*d*b*c*b*c*a*d*(c^{-1})*d*(c^{-1})=1.
\end{gather*}

\subsection{Group \#30. [order 12,288]}
This group has the structure:
\begin{equation*}
                  C_2\times [1024 @  C_3] @ C_2,
\end{equation*}
where the order 6144 factor has 44 classes and a trivial center. This order 3072 group appears to be the same as the 3072 factor in the automorphism group 64 \#93. See group 64 \#173 below for an alternate presentation of the following groups of orders 1024 and 3072. \\

\underline{relations for order 1024 group:}
\begin{gather*}
a^2=b^2=c^2=d^2=e^2=(a*b)^2=(a*c)^2=(b*c)^2=\\
a*d*e*d*a*e=(a*d)^4=(a*d*b*d)^2=(a*d*c*d)^2=\\
(a*e)^4=(a*e*b*e)^2=(a*e*c*e)^2=\\
(b*d)^4=\\
(b*d*c*d)^2=(b*e)^4=(b*e*c*e)^2=(c*d)^4=(c*e)^4=\\
a*b*c*d*b*e*c*e*d=1.
\end{gather*}
This 2-group has an automorphism group of order $2^{20}*3^2=9,437,184$, and center $C_2$.

\underline{possible $C_3$ actions:}\\
The following two actions seem to generate the same group of order 3072 which appears in this decomposition:
\begin{gather*}
f^3=\\
\begin{cases}
a^f*e*b*e*d*c*d*c*b=b^f*e*d*e*c*a*d*a=\\
c^f*e*b*a*e*c*a=d^f*c*d*c*e*a=e^f*d=1 \\
\text{or}\\
(a,f)=b^f*e*c*b*e=c^f*e*d*e*c*d*c*b=\\
d^f*d*c*a*e*c=e^f*b*e*d*e*b=1. \\
\end{cases}
\end{gather*}

\subsection{Group \#44. [order 6,144]}
The structure of this group is:
\begin{equation*}
[768 @ \left(C_2\times C_2\times C_2\right)],
\end{equation*}
where the 768 group has 36 classes and a center of $C_2\times C_2$.
This 768 group appears to be the same as the one appearing in group \#17 above and \#108 below. Cayley gives the following presentation for this automorphism group on five generators:
\begin{gather*}
a^2=b^2=d^4=c^2=e^2=(a*b)^2=a*d*a*(d^{-1})= \\
(a*e)^2=b*d*b*(d^{-1})=(b*e)^2=(d*(e^{-1}))^2=(a*c)^3=\\
c*d^2*c*(d^{-2})=(b*c)^4=(b*c*d*c)^2=  \\
c*d*c*d*c*(d^{-1})*c*(d^{-1})= \\
(c*d*c*e)^2=(c*e)^4=a*b*c*a*b*c*b*a*c=\\
b*c*b*e*c*b*e*c*e*b*c*e=1.
\end{gather*}
\subsection{Group \#68. [order 12,288]}

This automorphism group appears to be isomorphic to that for group \#22 above.\\

This group does not have a normal sylow 2-subgroup. The structure looks like:
\begin{equation*}
[(256\times C_2\times C_2)\,\, @\,\,  C_3]\,\, @\,\,\left( C_2\times C_2\right),
\end{equation*}
where the order 256 2-group is generated by
\begin{gather*}
a^2=b^2=c^2=d^2=(a*b)^2=(c*d)^2=(a*c)^4=\\
(a*c*a*d)^2=(a*c*b*c)^2=(a*d)^4=\\
(a*d*b*d)^2=(b*c)^4=(b*c*b*d)^2=(b*d)^4=1.
\end{gather*}
This group is [256,\#8935] in the Small Group Library.
The order of the automorphism group of this 2-group of order 256 is $2^{19} * 3^2 = 4,718,592$. CAYLEY gives the following four possible $C_3$ actions on this 2-group: \\

\underline{Case 1.} This case has a trivial center. This action gives group 768 \#1083695, with an automorphism group of order 36,864 and center I:
\begin{align*}
&a^e*((a*d*a*b*d)^{-1})=b^e*((a*d*b*d)^{-1})=\\
&c^e*((a*b*c*b*d*a*d)^{-1})=d^e*((a*b*c*a*c*d*b)^{-1})=1.
\end{align*}

\underline{Case 2.} This case has center $C_2 \times  C_2$ and is the desired 3072 group. This corresponds to 768 \#1083694 and has an automorphism group of order 36,864 that is not isomorphic to that in Case 1:
\begin{align*}
&a^e*((d*a*b*d)^{-1})=b^e*((a*b*c*b*c)^{-1})=\\
&c^e*((a*c*a*b*d*b)^{-1})=d^e*((a*b*d*a*c*b*d)^{-1})=1.
\end{align*}

\underline{Case 3.} This action generates the same group as in Case 2:
\begin{align*}
&(a,e)=b^e*((a*c*a*b*c)^{-1})=\\
&c^e*(d^{-1})=d^e*(c*d)^{-1})=1.
\end{align*}
		
\underline{Case 4.} This case is isomorphic to Case 2:
\begin{align*}
&a^e*((d*a*b*d)^{-1})=b^e*((a*b*c*d*b*c*d)^{-1})=\\
&c^e*((a*b*c*b*d*a*c)^{-1})=d^e*((a*c*a*b*d*b)^{-1})=1,
\end{align*}
with the presentation for the $\left[(256) \times C_2\times C_2)@C_3\right]$ subgroup given by:
\begin{gather*}
a^2=b^2=c^2=d^2=(a*b)^2=(c*d)^2=(a*c)^4=\\
(a*c*a*d)^2=(a*c*b*c)^2=(a*d)^4=\\
(a*d*b*d)^2=(b*c)^4=(b*c*b*d)^2=(b*d)^4=\\
 e^3=f^2=h^2=(f,h)=(a,f)=(b,f)=(c,f)=(d,f)=\\
(a,h)=(b,h)=(c,h)=(d,h)=\\
a^e*((a*d*a*b*d)^{-1})=b^e*((a*d*b*d)^{-1})=\\
c^e*((a*b*c*b*d*a*d)^{-1})=d^e*((a*b*c*a*b*c*d)^{-1}=\\
f^e*h=h^e*f*h=1.
\end{gather*}

The order of the automorphism group of the order 1024 group is $2^{36}*3^3$. \\

The automorphism group of the resulting order 3072 group has order $2^{20}*3^3=28,311,552$, with 223 conjugacy classes and a trivial center. The next step would be to find the $C_2\times C_2$ actions on the resulting order 3072 group. This, in view of the order of the automorphism group of this group, is difficult to do computationally.\\

Cayley gives the following presentation for this automorphism group of order 12,288 on five generators:
\begin{gather*}
 a^2=b^2=c^2=d^2=e^4=(b*c)^2=(b*d)^2=\\
 (b*e)^2=(c*d)^2=a*e^2*a*(e^{-2})=(c*e)^3=\\
 d*e^2*d*(e^{-2})=a*c*a*c*e*a*(e^{-1})=\\
 (a*b)^4=a*b*a*b*c*(e^{-2})*c=(a*d)^4=\\
 d*e*d*e*d*(e^{-1})*d*(e^{-1})=c*d*e*c*d*e*d*c*e=\\
 a*c*d*a*c*d*a*d*c*a*d*c=\\
 a*d*a*e*d*e*a*d*a*(e^{-1})*d*(e^{-1})=1.
\end{gather*}

\subsection{Group \#76. [order 6144]}
This group has a normal sylow 2-subgroup. The order of the automorphism group of the sylow 2-subgroup is $2^{34}*3^2*7$. A presentation in terms of the sylow 2-subgroup and a  $C_3$ quotient is given by:
\begin{align*}
 &a^2=b^2=c^2=d^2=e^2=(a*b)^2=(a*c)^2=\\
 &(b*c)^2=(d*e)^2=(a*d)^4=(a*d*a*e)^2=\\
&(a*d*b*d)^2=(a*d*c*d)^2=(a*e)^4=\\
 &(a*e*b*e)^2=(a*e*c*e)^2=(b*d)^4=\\
 &(b*d*b*e)^2=(b*d*c*d)^2=(b*e)^4=\\
&(b*e*c*e)^2=(c*d)^4=(c*d*c*e)^2=(c*e)^4=\\
 &f^3=      \\
 &a^f*((b*d*e*b*d*a*c*e)^{-1})=\\
&b^f*((a*c*d*e*a*b*d*e)^{-1})=   \\
&c^f*((a*b*c*d*c*d*e*a*e)^{-1})=\\
&d^f*((a*b*d*b*e*a)^{-1})=   \\
&e^f*((a*b*e*a*d*e*b)^{-1})=1
\end{align*}
or
\begin{align*}
 & a^f*((b*d*e*b*d*a*c*e)^{-1})=\\
 & b^f*((c*e*b*e)^{-1})=   \\
 & c^f*((d*b*d)^{-1})=   \\
 & d^f*((a*b*c*d*b*c*e*a*d)^{-1})=\\
 & e^f*((a*b*c*e*a*b*d*c)^{-1})=1.
\end{align*}
The automorphism group has order 4,718,592 $(2^{19}*3^2)$ and a trivial center.\\

Cayley gives the following presentation for this automorphism group
 (Aut[64 \#76]) using only four generators:
\begin{gather*}
a^2=b^2=c^2=d^6=(b*c)^2=a*b*a*d*b*(d^{-1})=  \\
a*d*a*d*a*(d^{-2})=(a*b)^4=(a*b*a*c)^2=(a*c)^4=\\
(b*d*c*(d^{-1}))^2=(c*d*c*(d^{-1}))^2=\\
c*d^3*c*(d^{-3})=(a*c*(d^{-1}))^3=  \\
a*b*(d^{-1})*b*d*a*(d^{-1})*b*d*b= \\
a*c*d*a*c*d*a*c*(d^{-2})=1.
\end{gather*}

Another alternate 4-generator version of Aut(\#76) is
\begin{gather*}
a^2=b^2=c^2=d^2=(a,b)=(a,c)=(a,d)=(c,d)=\\
a*c*e*a*c*e^{-1}=b*c*b*e*c*e^{-1}=\\
e^6=a*b*c*b*e^{-1}*c*e=b*e*b*e*b*e^{-2}=\\
(b*c*b*d)^2=(b*d)^4=(d*e*d*e^{-1})^2=\\
d*e^3*d*e^{-3}=(b*d*e^{-1})^3=\\
(b*d*e)^2*b*d*e^{-2}=1.
\end{gather*}

\subsection{\underline{Group \#81. [order 12,288]}}
See Group \#30 above.\\

\subsection{\underline{Group \#82. [order 12,288]}}
 This automorphism group has the structure [512@$C_3$]@$D_4$. The automorphism group of 512@$C_3$ has
order 589,824.\\

\subsection{\underline{Group \#93. [order 6,144]}}
The structure of this group is
\begin{equation*}
C_2\times [1024\,\,\, @\,\,  C_3].
\end{equation*}
The automorphism group of the 1024 group has order $2^{20}*3^2$. See Group \#30.

\subsection{\underline{Group \#103. [order 36,864]}}
A presentation for this automorphism group is:
\begin{gather*}
 a^2=b^4=c^4=d^2=e^4=a*b*a*(b^{-1})=\\
 a*c*a*(c^{-1})=a*e*a*(e^{-1})=\\
 b*d*(b^{-1})*d=a*b^2*d*a*d=\\
 b^2*c*(b^{-2})*c=b^2*e*(b^{-2})*e=\\
 b*c^2*(b^{-1})*(c^{-2})=\\
 b*e*b*(e^{-1})*(b^{-1})*(e^{-1})=\\
c^2*d*(c^{-2})*d=c^2*e*(c^{-2})*e=\\
 d*e^2*d*(e^{-2})=\\
b*c*b*c*(b^{-1})*(c^{-1})*(b^{-1})*(c^{-1})=\\
 b*c*(d^{-1})*c*(b^{-1})*(c^{-1})*(d^{-1})*(c^{-1})=\\
 b*c*e*c*(b^{-1})*(c^{-1})*(e^{-1})*(c^{-1})=\\
 c*d*c*d*(c^{-1})*d*(c^{-1})*d=\\
 c*d*(c^{-1})*(e^{-1})*(c^{-1})*d*c*e=\\
 c*e*c*e*(c^{-1})*(e^{-1})*(c^{-1})*e=\\
 d*e*d*e*d*(e^{-1})*d*(e^{-1})=\\
 b*d*e*b*d*(e^{-1})*d*(b^{-1})*(e^{-1})=1.
\end{gather*}

\subsection{\underline{Group \#104. [order 61,440]}}
A presentation on four generators for this automorphism group is:
\begin{gather*}
 a^2=b^4=c^4=a*c*a*(c^{-1})=\\
 a*d*a*(d^{-1})=b^2*(d^{-2})=\\
 b^2*c*(b^{-2})*c=b*c^2*(b^{-1})*(c^{-2})=\\
 a*b*a*b*a*(b^{-1})*a*(b^{-1})=\\
 a*b*a*b*d*(c^{-2})*(d^{-1})=\\
 b*d*b*d*(b^{-1})*(d^{-1})*(b^{-1})*(d^{-1})=\\
 c*d*c*(d^{-1})*(c^{-1})*(d^{-1})*c*(d^{-1})*(c^{-1})*(d^{-1})=\\
 a*b*a*b*c*b*a*c*b*a*(b^{-1})*c*(b^{-1})*(c^{-1})=\\
 b*c*b*c*d*c*b*c*d*(c^{-1})*(d^{-1})*\\
(b^{-1})*(d^{-1})*(c^{-1})*(d^{-1})=1.
\end{gather*}

\subsection{\underline{Group \#105. [order 23,040]}}
A presentation for this complete group is:
\begin{gather*}
 a^4=b^2=c^4=d^2=a*b*(a^{-1})*b=a^2*c*(a^{-2})*c=\\
a*c^2*a*(c^{-2})=a^2*b*d*(a^{-2})*b*d=\\
 a*d*a*d*(a^{-1})*d*(a^{-1})*d=a*d*c*d*a*d*(c^{-1})*d=\\
 (b*c)^4=(b*c*b*(c^{-1}))^2=(c*d*(c^{-1})*d)^2=\\
 a*c^2*(d^{-1})*(b^{-1})*(d^{-1})*(a^{-1})*(c^{-2})*(b^{-1})=\\
 a*c*a*c*d*(c^{-2})*b*d*(a^{-1})*b*c=\\
 a*c*b*c*d*a*d*(a^{-1})*b*c*b*c=1.
\end{gather*}
\subsection{\underline{Group \#108. [order 12,288]}}
This group has the structure
\begin{equation*}
768 @ (D_4\times C_2).
\end{equation*}
The order 768 group appears to be the same one as appearing in groups \#17 and \#44 above. A set of relations for this form for the automorphism group has not yet been found. \\

An alternate presentation for this automorphism group is:
\begin{gather*}
 a^2=b^4=c^2=a*b*a*b*a*(b^{-1})*a*(b^{-1})=\\
 a*c*b^2*c*a*c*(b^{-2})*c= \\
b^2*c*(b^{-1})*c*(b^{-2})*c*(b^{-1})*c=\\
a*b^2*a*b*c*a*c*b*c*a*c=\\
a*b^2*c*b*c*a*c*(b^{-1})*c*(b^{-2})=\\
 b*c*b*c*b*c*(b^{-1})*c*(b^{-1})*c*(b^{-1})*c=\\
 a*b*a*b^2*c*a*c*a*(b^{-1})*a*c*a*c=\\
  a*b*c*b*a*b*c*b*a*(b^{-1})*c*(b^{-1})*a*(b^{-1})*c*(b^{-1})=1.
\end{gather*}
\subsection{\underline{Group \#109. [order 1,536]}}
This group has the structure
\begin{equation*}
[\text{order 96, 11 classes}]\,\, @\,\, (1^4 \text{group}),
\end{equation*}
where the automorphism group of the order 96 group has order 576. It appears that this extension is a nonsplit extension,
since the straightforward application of the $(C_2\times C_2\times C_2\times C_2)$ action on the order 96
group does not give the correct automorphism group.\\

An alternate presentation for this automorphism group is:
\begin{gather*}
 a^2=b^2=c^2=d^2=e^4=(a*b)^2=(a*d)^2=a*e*a*(e^{-1})=\\
 (b*c)^2=b*e*b*(e^{-1})=(d*e)^3=(a*c)^4=\\
a*c*a*c*d*(e^{-2})*d=(b*d)^4=b*d*b*e*d*(e^{-1})*d*e=\\
 c*e*c*e*c*(e^{-1})*c*(e^{-1})=c*e*d*(e^{-1})*c*(e^{-1})*d*e=1.
\end{gather*}
\subsection{\underline{Group \#153. [order 10,752]}}
This automorphism group has the structure
\begin{equation*}
            (1^9 \,\,@ \,\,C_7)\,\, @\,\,  C_3.
\end{equation*}
The action of the group $C_7 @  C_3$ on the elementary abelian group of order 512 is obtained by considering the  $C_3$ and the $C_7$ to act on each of the three triplets $C_2\times C_2\times C_2$ as
follows:
\begin{gather*}
a^2=b^2=c^2=(a,b)=(a,c)=(b,c)=\\
h^7=k^3=h^k*h^5=\\
a^h*b*c=b^h*a=c^h*b=\\
(a,k)=b^k*c=c^k*a*b*c=1.
\end{gather*}
A Cayley-generated presentation for this automorphism group is:
\begin{gather*}
 a^3=b^3=c^3=(a*(b^{-1}))^2=\\
 a*b*a*c*b*(a^{-1})*c*(a^{-1})*b*c=\\
 a*c*a*c*(a^{-1})*(c^{-1})*(a^{-1})*c*(a^{-1})*(c^{-1})=\\
 a*b*(c^{-1})*b*c*a*c*(b^{-1})*(c^{-1})*(a^{-1})*c=\\
 a*c*(a^{-1})*(b^{-1})*(c^{-1})*(a^{-1})*(c^{-1})*b*a*c*b*(c^{-1})=1.
\end{gather*}
\subsection{\underline{Group \#173. [order 6,144]}}
This automorphism group does not have a normal sylow 2-subgroup. It does however have a normal subgroup of order 3072 that does have a normal sylow 2-subgroup. The automorphism group thus can be represented as
\begin{equation*}
                  [1024 \,\,@ \,\, C_3]\,\, @\,\, C_2.
\end{equation*}
There are three different possible groups of order 3072 that can be used to express the automorphism group of the form here. The automorphism groups of these groups have the orders 294,912 (2 cases) or 589,824. An alternate presentation on four generators for the 1024 group is:
\begin{gather*}
 a^2=b^2=c^2=d^2=(a*b)^2=(a*c)^4=(a*c*b*c)^2=\\
 (a*d)^4=(a*d*b*d)^2=(b*c)^4=(b*d)^4=(c*d)^4=\\
 a*c*d*c*d*a*d*c*d*c=b*c*d*c*d*b*d*c*d*c=\\
 a*b*c*a*b*d*a*b*c*b*a*d=\\
a*c*d*a*c*d*a*d*c*a*d*c=1.
\end{gather*}
The corresponding  $C_3$ action to yield the order 3072 group is:
\begin{gather*}
e^3=a^e*b=b^e*b*a=c^e*d*a*c*a*d=\\
d^e*b*c*a*d*a*c*b=1
\end{gather*}
and yields a group that has 46 conjugacy classes and an automorphism group of order 589,824.\\

A presentation for the full automorphism group on four generators is:
\begin{gather*}
 a^4=b^4=c^4=d^2=a^2*b*(a^{-2})*b=a^2*c*(b^{-2})*(c^{-1})=\\
  a^2*(c^{-1})*(b^{-2})*c=a^2*d*(a^{-2})*d=\\
 a*b^2*(a^{-1})*(b^{-2})=b^2*d*(b^{-2})*d=\\
 (b*(c^{-1}))^3=c^2*d*(c^{-2})*d=\\
 a*b*a*b*(a^{-1})*(b^{-1})*(a^{-1})*b=\\
 a*b*c*b*(c^{-1})*(a^{-1})*b*c=\\
  a*c^2*a*(c^{-2})*(b^{-2})=a*d*a*d*(a^{-1})*d*(a^{-1})*d=\\
 b*c*(b^{-1})*d*(b^{-1})*(c^{-1})*b*d=\\
 b*d*b*d*(b^{-1})*d*(b^{-1})*d=\\
 c*d*c*d*(c^{-1})*d*(c^{-1})*d=\\
 a^2*b*c*d*(c^{-1})*d*b*d=\\
 a*b*d*a*(c^{-1})*a*c*d*(a^{-1})*b=\\
 a*b*d*b*(a^{-1})*d*(a^{-1})*(b^{-1})*d*(b^{-1})*a*d=1.
\end{gather*}

\subsection{\underline{Group \#181. [order 1,536]}}

This automorphism group does not have a normal sylow 2-subgroup. The largest normal sylow 2-subgroup has order 256 (\#55999). Of the three normal subgroups of order 768, only one has a normal sylow 2-subgroup. A possible representation of this automorphism group is therefore
\begin{equation*}
[(256)\,\, @\,\,  C_3]\,\, @\,\, C_2.
\end{equation*}
The order 768 group has the following presentation:
\begin{gather*}
a^2=b^2=c^2=d^2=e^2=f^2=(a*b)^2=(a*c)^2=\\
(a*d)^2=(a*e)^2=(b*c)^2=(b*e)^2=(b*f)^2=\\
 (c*d)^2=(c*f)^2=(d*e)^2=(d*f)^2=(e*f)^2=\\
 (a*f)^4=(b*d)^4=(c*e)^4=\\
a*b*c*d*b*d*e*f*a*f*c*e=\\
h^3=\\
a^h*f*a*f*c*a=b^h*f*a*f=\\
c^h*f*a*f*b*a=\\
d^h*f=\\
e^h*b*d*b=\\
f^h*c*e*c=1.
\end{gather*}
This group of order 768 (\#1085205) appears to be a characteristic subgroup of the complete group of order 12,288 below.\\

The smallest order normal subgroup with a factor of three in its order is 192. This order 192 group is \#61 in Table 2a, which is \#1023 in the Small Group Library. The quotient group aut(\#181)/(normal subgroup of order 192) is $D_4$. (This order 192 group has 9 classes and is \#183 @  $C_3$.)\\

A presentation for this automorphism group of order 1536 is:
\begin{align*}
&b^4=c^4=d^2=b*d*(b^{-1})*d=(c*d)^2=\\
&a^3*d*(a^{-1})*d=a^2*b*a^2*(b^{-1})=\\ 
&(a*(c^{-1}))^3=b*c^2*b*c^2=\\
& a*b*c*b*( c^{-1}) * (a^{-1} ) * ( b^{-1})=\\ 
&a*c*a*c*b^2*a*(c^{-1})=1.
\end{align*}
This group has 33 classes and the order structure:\\

\begin{tabular}{|c|c|c|}\hline
order of element& number of elements& number of classes\\ \hline   
2& 159&10\\	   
3& 128&1\\   
4& 480&13\\	   
6& 384&3\\   
8& 384&5\\ \hline	 
\end{tabular}
\linebreak\\
The group presentation below gives rise to the automorphism group of this order 1536 group and has order 6144. This is the presentation used as the starting point for the automorphism tower here (see Appendix III):
\begin{gather*}
 b^4=c^4=d^4=b*d*b*(d^{-1})=(c*(d^{-1}))^2=\\
 a^4*(d^{-2})=a^2*b*a^2*(b^{-1})=a*b*d*a*d*(b^{-1})=\\
 a*b*(d^{-1})*a*(d^{-1})*(b^{-1})=\\
a*(b^{-1})*d*a*(b^{-1})*(d^{-1})=\\
 (a*(c^{-1}))^3=b^2*d*(c^{-1})*(d^{-1})*(c^{-1})=\\
 b*c^2*b*(c^{-2})=a*b*c*b*(c^{-1})*(a^{-1})*(b^{-1})=\\
 a*c*a*c*(b^{-2})*a*(c^{-1})=1.
\end{gather*}

\subsection{\underline{Group \#183. [order 9,216]}}
The automorphism group does not have a normal sylow 2-subgroup. The largest-order normal 2-group is the
elementary abelian group of order 256. A group of order 2304 has a normal sylow 2-subgroup. This group has
the structural form:
\begin{equation*}
                  [(256) @ ( C_3\times  C_3)]  @ \left(C_2\times C_2\right).
\end{equation*}
Possible presentations for this group of order 2304 are:
\begin{gather*}
a^3=b^3=c^3=(a*(b^{-1}))^2=(a*(c^{-1}))^2=(b*(c^{-1}))^2=\\
(a*c*b)^2=\\
(a,a1)=(b,a1)=(c,a1)=\\
a1^2=b1^2=c1^2=d1^2=\\
(a1,b1)=(a1,c1)=(a1,d1)=(b1,c1)=(b1,d1)=(c1,d1)=\\
d^3=\\
a1^d*b1=b1^d*a1*b1=c1^d*d1=d1^d*c1*d1=\\
(a,b1)=(b,b1)=(c,b1)=      \\
(a,c1)=(b,c1)=(c,c1)=\\
(a,d1)=(b,d1)=(c,d1)=
\end{gather*}
coupled with one of the following actions of the order-three element $d$ on the generators $a$, $b$, and $c$. Cases 1, 2 and 3 all seem to generate this order 2304 group. The automorphism group of this group has order $2^{12}*3^4*5$. Case 4 has a different number of classes (64) and order structure [255 elements of order 2 in 35 classes, 1088 elements of order 3 in 8 classes, 960 elements of order 6 in 20 classes, and 255 elements of order 2 in 35 classes]. Case 4's automorphism group has order $2^{15}*3^4*5^2$.\\
Case 1.
\begin{equation*}
a^d*c*b=b^d*((b*c)^{-1})*a=(c,d)=1.
\end{equation*}
Case 2.
\begin{equation*}
a^d*b*((c*a)^{-1})=b^d*b*a=c^d*a*b=1.
\end{equation*}
Case 3.
\begin{equation*}
(a,d)=b^d*c*a=c^d*c*((a*b)^{-1})=1.
\end{equation*}
Case 4.
\begin{equation*}
a^d*b*((c*a)^{-1})=b^d*(c^{-1})*b*(a^{-1})=c^d*c*((a*b)^{-1})=1.
\end{equation*}

A three-generator set of relations for this automorphism group is:
\begin{gather*}
 a^2=b^{12}=c^2=(a*b*a*(b^{-1}))^2=a*(b^{-2})*c*a*c*b^2=\\
 (a*c)^4=a*b^2*a*c*a*c*(b^{-2})=(a*b*a*(b^{-1})*c)^2=\\
 a*b*c*a*b*c*a*c*(b^{-1})*c*(b^{-1})=a*b*a*b*c*b^4*c*b^2=\\
  (b^2*(c^{-1}))^4=(b*(c^{-1}))^6=\\
  a*b*c*(b^{-1})*a*c*b*c*(b^{-1})*c*b*c*(b^{-1})*c=1.
\end{gather*}
\subsection{\underline{Group \#187. [order 15,360]}}
The structure of this automorphism group, which is a complete group, has the form:
\begin{equation*}
( [(1^4 @  C_3)\times 1^4 ] @ C_5 ) @ C_4.
\end{equation*}
A presentation for the $\left[\left(1^4@C_3\right) \times 1^4 \right ]@C_5 $ normal subgroup of order 3840 is:
\begin{gather*}
a^3=b^3=c^3=(a*(b^{-1}))^2=(a*(c^{-1}))^2=(b*(c^{-1}))^2=\\
(a*c*b)^2=\\
a1^2=b1^2=c1^2=d1^2=(a1,b1)=(a1,c1)=(a1,d1)=\\
(b1,c1)=(b1,d1)=(c1,d1)=\\
(a,a1)=(b,a1)=(c,a1)=(a,b1)=(b,b1)=(c,b1)=\\
(a,c1)=(b,c1)=(c,c1)=\\
d2^5= \\
(a,d2)=b^{d2}*((b*c)^{-1})*a=c^{d2}*b*((c*a)^{-1})=\\
a1^{d2}*a1*c1*d1=b1^{d2}*a1=c1^{d2}*a1*b1*c1*d1=d1^{d2}*a1*d1=1.
\end{gather*}
The group $<a,b,c>$ is $1^4 @  C_3$ and is the group of order 48 whose automorphism group has order 5760. The actions of $C_5$ on each of the two pairs of $1^4$ groups should be  as in the group of order 80 with automorphism group of order 960.\\

In the representation above we have the order-5 automorphism acting on the group of order 48 directly. We have not yet found a consistent representation with the form
\begin{equation*}
                  [1^4\times (1^4\,\, @\,\,  C_3)\,\, @\,\, C_5] @ C_4.
\end{equation*}\\

The $C_4$ acts as an element of order 2 on the  $C_3$ and as an element of order 4 on the $C_5$. Relations based upon this normal subgroup structure have not yet been found.\\

A presentation of this automorphism group using only three generators is:
\begin{gather*}
(a*(b^{-1}))^2=a^3*(b^-3)*(a^{-1})*b= \\
 a*b*a*b*(a^{-1})*(b^{-1})*(a^{-1})*(b^{-1})= \\
 a*b*a*c*(b^{-1})*c*b*(c^{-1})=  \\
a*b*c*b*(c^{-2})*b*c=  \\
b^2*c*(b^{-1})*(c^{-2})*b*(c^{-1})=\\
a^2*(c^{-1})*(a^{-1})*(b^{-1})*a*c^2*(a^{-1})*(c^{-1})=\\
a*c*(a^{-1})*(c^{-2})*(a^{-1})*b*c*a*(c^{-1})=1.
\end{gather*}

\subsection{\underline{Relations for the group of order 1536 appearing in the}\\
   \underline{automorphism groups of numbers 144, 145, 147, 148 and 184}}
\begin{gather*}
a^2=b^2=c^2=d^2=e^2=f^2=(a,b)=(a,c)=(a,d)=(a,e)=(a,f)=\\
(b,c)=(b,d)=(b,e)=(b,f)=(c,d)=(c,e)=(c,f)=(d,e)=(d,f)=(e,f)=\\
g^3=a^g*b=b^g*a*b=c^g*d=d^g*c*d=e^g*f=f^g*e*f=\\
h^2=(a,h)=(b,h)=(c,h)=(d,h)=(e,h)=(f,h)=(g,h)=\\
k^2=a^k*b=c^k*d=e^k*f=(h*k)^4=g^k*g=1.
\end{gather*}
The structure of this group is:
\begin{equation*}
            [ (1^6\,\, @\,\,  C_3)\times C_2\times C_2 ]\,\, @\,\, C_2
\end{equation*}
with the action of $C_2$ on the order 192 group as follows:
\begin{align*}
&C_2 \text{ acts on pairs of } C_2\times C_2 \text{ as in } C_2 \wr C_2.\\
              &C_2 \text{ inverts the generator of } C_3.\\
&\text{The action of }C_2 \text{ on the other }C_2\times C_2 \text{ is as a wreath product}.\\
\end{align*}
The group $(1^6 @  C_3)\,\, @\,\, C_2$ has for its automorphism group the
complete group of order 64,512 appearing in \#144,... above.\\

\subsection{\underline{Some common 2-groups arising in the above automorphism groups}}
The following 2-groups arising in the above discussions seem to be isomorphic:
\begin{quotation}
            the 1024 groups in cases 64\#'s 30, 81, 93, 173
                  (yes, Dr. Newman confirmed these).
                  The 3072 groups $(1024)\,\, @\,\,  C_3$ in cases
                        \# 30, 81 and 91 appear isomorphic but
                        \#173 is different.
                        The automorphism groups for \#30 and \#81
                              have the form $(3072) @ C_2$ and
                              may be isomorphic.
\end{quotation}
Also, 
\begin{quotation}
            the 256 groups and the 768 groups that occur in the automorphism groups of \#'s 17, 44, and 108 of order 64 are all isomorphic.
These automorphism groups have the structures:\\
     $[[(256)\,\, @\,\,  C_3]\,\, @\,\, (C_2\times C_2)]\times C_2,$\\
     $[(256) \,\,@\,\,  C_3] \,\,@ \,\,(C_2\times C_2\times C_2),$ and\\
                        $[(256) \,\,@\,\,  C_3] \,\,@ \,\,(D_4\times C_2)$\\
                  with the same $(256\,\, @\,\,  C_3)$ normal subgroup.
In other cases, such as \#'s 104, 105, 153, and 187, among others, presentations on three or four generators are available.
\end{quotation}
\newpage

\section{Appendix III \\ Automorphism Towers for Groups in Table A1}
The presentations for these automorphism towers were originally computed by CAYLEY runs in the mid 1990s. When we started revising this report only a paper copy of these runs was on hand.\footnote{Subsequently the computer file which contained computer transfers of these presentations was located. Even so we still checked these presentations with runs using GAP.} In order to check these presentations in their current published form, the following was done. This set of presentations was scanned in from a computer transfer listing, but the scanning was not perfect; hence a fair amount of hand editing was involved in getting this listing for these presentations. The listings given previously for the automorphism groups of the groups of order 64 were taken directly from computer readable files, and hence are more reliable than the ones given here. In those cases (the automorphism groups for the groups of order 64) one also did some transfers from an original text editor output from CAYLEY to an MS Word format and then to the LaTeX form, so even there some minor editing was done. These presentations were run recently in GAP in order to verify that these presentations gave the correct orders and in many cases the correct class structures for these previously obtained results. One word of caution should be noted here. Namely, in the original CAYLEY runs, the presentations from one member of the automorphism tower were used as the input for the next group assuring the continuity of the automorphism tower structures. With GAP this was not done for the following reason. In CAYLEY the output was reproducible while in GAP this is not the case. The way things were done with CAYLEY was the following. One ran the initial group in CAYLEY and asked for, among other things, the automorphism group of the input group. In many cases these CAYLEY generated automorphism groups had large orders and many generators. We then asked CAYLEY to select a subset of this initial generating set that would also give this automorphism group, e.g., say generators numbers 1, 3, 6, 8 and 10. We then ran CAYLEY a second time and asked CAYLEY to get a presentation for the automorphism group using this specific generating set. Then further calculations on the automorphism group were carried out using this presentation. In the automorphism towers this process was then repeated for each subsequent group in the series until either the tower terminated or became too large for CAYLEY to continue up the sequence. In some cases these runs took several hours, depending upon the degree of the automorphism group and the number of generators needed to specify the group, so doing things without these \textquotedblleft built-in breaks" would be rather impractical. One really does need to know just how many generators one needs to specify the various groups, and this seems to be the best way to do that. Furthermore the use of permutations is not as efficient as the use of a presentation input. From our experience with CAYLEY, many calculations using permutations become very time-consuming or impossible, whereas with a presentation input many of these calculations are doable and are considerably faster.\\

In order to check these presentations, the following was done. We scanned the printed file into the computer. We then took the edited and hopefully corrected scanned presentations and made the following substitutions: $a\rightarrow f.1$, $b\rightarrow f.2$, etc., and = in the presentation was replaced by a comma, and we rewrote $a^{-2}\rightarrow$ f.1$\wedge-2$, etc. to get a GAP type input format. We then ran these edited presentations in GAP. Some errors were found in these initial runs, and hopefully the corrections have been incorporated back into the original LaTeX manuscript. A comment or two is in order here as well. In the CAYLEY runs, the orders of the centers and the order/class structure were obtained reasonably quickly. With GAP the determinations of the centers of the groups (mostly of orders 1, 2, 4, or 8) were very slow. This is probably a built-in problem with GAP since GAP converts everything to permutations before doing any calculations. In using CAYLEY it was found that using presentations rather than permutations as input resulted in much faster run times, and in fact some calculations probably could not have been done using a permutation input for many groups.\footnote{It was suggested to us that we should check these runs using GAP, possibly to see if there were any errors in our original presentations. As a matter of general information we encountered many more errors (i.e., typos) due to the rewriting of these presentations in GAP in the process of checking them than were found in the original listings.}  \\

Some of the presentations could have been simplified, e.g., words such as 
\begin{equation*}
a*b*d*a*b*d*a*b*d*b*a*d = (a*b*d)^3 * b*a*d
\end{equation*}
and many others below. These words were left intact to avoid any additional typing/editing errors here.\\

We now proceed to give the presentations for the groups in these automorphism group tower sequences.\\

\subsection{\underline{Group \#5.} Automorphism tower for $C_4 \times  C_4 \times  C_4$}  
a. The automorphism group of $C_4 \times  C_4 \times  C_4$ is
\begin{align*}
&a^2=b^2=c^2=d^2=(a*b)^2=(a*c)^3=(a*d)^3=\\
&(b*c)^3=(c*d)^4=(a*b*d*c*d)^2=\\
&a*b*d*a*b*d*a*b*d*b*a*d=\\
&a*c*a*d*a*c*a*d*a*c*d*b*d*c*a*d=\\
&a*b*c*a*d*a*c*d*c*b*d*c*b*d*b*d*c*b*d=1.
\end{align*}
b. The second group in this sequence is complete:
\begin{align*}
&a^2=b^2=c^2=d^2=e^2=(a*b)^2=(a*e)^2=\\
&(a*c)^3=a*c*a*e*d*e=(b*c)^3=(b*e)^4=\\ 
&(c*d)^4=a*b*d*a*c*e*b*d*e*d=\\
&a*b*d*a*b*d*a*b*d*b*a*d= \\
&a*c*a*d*a*c*a*d*a*c*d*b*d*c*a*d=\\ 
&a*b*c*a*b*d*c*e*c*b*d*c*b*d*b*c*d*c*e=1.
\end{align*}
\subsection{\underline{Group \#7}} Automorphism sequences for this group start out with $C_2\ \times$  (order 384 group).\\
a. Aut(64 \# 7) is $C_2$ cross the following order 384 (\# 20095) group:
\begin{align*}
&d^2=(a*b^{-1})^2=(a*c^{-1})^2=(b*c^{-1})^2= \\
&a^2*(b^{-2})*(a^{-1})*b = \\
&a^2*(c^{-2})*(a^{-1})*c=\\
& a^2*(d^{-1})*c*a*(d^{-1})=(a*c*b)^2= \\
&(a*(d^{-1})*(b^{-1}))^2=\\
&a*(c^{-1})*(b^{-1})*(d^{-1})*(a^{-1})*(b^{-1})*a*(d^{-1})=1.
\end{align*}

b. The second entry in the sequence is Aut(384) with order $2^8 * 3^2$ = 2304:
\begin{align*}
&a^2=d^2=e^2=a*c*a*(c^{-1} )=( a*e)^2=\\
&b^2*( c^{-2})=b*c*(b^{-1} )*(c^{-1})=b*e*(b^{-1} )*e=\\
&c*e*(c^{-1})*e=(d*e)^2=( a*b)^3=\\
&a*b*d*b*a*(b^{-1})*d*(b^{-1})=(a*d)^4=\\
&b^2*c*d*c^{-1}*b^{-2}*d=(b*d^{-1})^4=\\
&(c*(d^{-1}))^4=a*b*a*d*(b^{-1})*d*a*d*(b^{-1})*d=1.
\end{align*}

c. The third entry in the sequence has order
 $2^{12} * 3^3$ = 110,592:
\begin{align*}
&c^2=d^2=a^2*b^{-2}=a*b*a^{-1}*b^{-1}=\\
&a*d*a^{-1}*d=a*e*a^{-1}*e^{-1}=(c*d)^2=d*e*d*e^{-1}=\\
&(b*d^{-1})^3=a^2*b*c*b^{-1}*a^{-2}*c=\\
&b*e*b*e*(b^{-1})*(e^{-1} )*(b^{-1})*(e^{-1})=\\
&b*d*b*e*b*e*d*(b^{-1})*(e^{-1})=\\
&a*c*b*(e^{-1} )*c*( a^{-1} )*c*e*(b^{-1})*c=\\
&a^2*e*c*e*c*e*(b^{-1})*c*b*e*c=\\
&a*c*a*c*(a^{-1} )*c*(a^{-1})*c*b*c*(b^{-1})*c=\\
&b*c*e*(b^{-1})*c*e*c*(b^{-1})*(e^{-1})*(b^{-1} )*c*e=\\
&(b*(c^{-1}))^6=\\ 
&a^2*c*a*c*(a^{-1})*c*(b^{-1})*c*(a^{-1})*b*c*a*c=1.
\end{align*}

d. The fourth entry has order $2^{13}*3^3$ = 221,184:
\begin{align*} 
&d^2=e^2=f^2=a^2*(b^{-2})=a^2*(c^{-2})=\\
&a*b*(a^{-1})*(b^{-1})=a*c*(a^{-1})*(c^{-1})=\\ 
&a*e*(a^{-1})*e=a*f*(a^{-1})*f=\\
& b*c*(b^{-1})*(c^{-1})=c*f*(c^{-1})*f=(d*e)^2=\\
&(d*f)^2=(c*(e^{-1}))^3=a^2*c*d*(c^{-1})*(a^{-2})*d=\\ 
&b*(d^{-1})*( f^{-1})*b*(f^{-1})*(c^{-1})*(d^{-1})*(c^{-1} )=\\ 
&b*e*b*e*(b^{-1})*e*(b^{-1})*e= \\
&b*e*(b^{-1})*f*(b^{-1})*e*b*f=\\
&(b*f*(b^{-1} )*d)^2=(e*f)^4 =\\
&b*c*e*b*c*e*c*b*e=\\
&b*(e^{-1})*(f^{-1} )*(e^{-1} )*( f^{-1} )*(b^{-1})*( f^{-1})* (e^{-1}) *\\
& \qquad 	(f^{-1})*(e^{-1})=\\
&a*d*a*d*(a^{-1})*d*(a^{-1})*d*c*d*(c^{-1})*d=\\
&b*(d^{-1})*b*(d^{-1})*b*(d^{-1})*b*(d^{-1})*c*\\
& \qquad (d^{-1})*c*(d^{-1})= \\
&a*(b^{-1})*(d^{-1})*b*(d^{-1})*(b^{-1})*(d^{-1})*(a^{-1})*\\
& \qquad  d^{-1}*a^{-1}*d^{-1}*b^{-1}*a*d^{-1}=1.
\end{align*}
\newpage
e. The last entry in this sequence is complete and has order = $2^{14}*3^3$ = 442,368:
\begin{align*}
&d^2=f^2=h^2=a^2*(b^{-2})=a^2*(c^{-2})=\\
& a^2*(e^{-2})=a*b*(a^{-1})*(b^{-1})=\\
&a*c*a^{-1}*( c^{-1} )=a*d*(a^{-1})*d=\\
&a*e*a*(f^{-1})=a*e*(a^{-1})*(e^{-1})=\\
&a*f*(a^{-1} )*f=a*k*(a^{-1} )*(k^{-1} )=\\ 
&b*c*( b^{-1}) * ( c^{-1})=b*d*(b^{-1})*d= \\
&b*e*(b^{-1})*(e^{-1} )=b*f*(b^{-1})*f= \\
&c*d*( c^{-1})*d=c*e*( c^{-1}) * (e^{-1} )=\\ 
&c*f*(c^{-1})*f=(d*h)^2=e*f*(e^{-1})*f=\\
& (f*h)^2=b*(c^{-1})*(k^{-1})*d*k=\\
&b*k*d*(k^{-1})*( c^{-1})=\\
&a*(k^{-1})*h*k*(a^{-1} )*h=\\
&b*e*(d^{-1})*(e^{-1})*(d^{-1})*(c^{-1})= \\
&b*( e^{-1} ) *h* ( c^{-1} ) *e*h=\\
&b*( k^{-2} )*( c^{-1} )*k^2=e*k^2*(e^{-1} )*( k^{-2} )=\\
&e*( k^{-4}) * ( f^{-1} )=\\
&b*k*b*(k^{-1})*(b^{-1})*(k^{-1})*(c^{-1})*k=\\
& b*e*k*b*k*e*(k^{-1})*(b^{-1})*(k^{-1})=\\
&a*h*a*h*(a^{-1})*h*(a^{-1} )*h*c*h*(b^{-1})*h=\\
& b*(h^{-1})*b*(h^{-1})*b*(h^{-1})*b*(h^{-1})*b*\\
& \qquad (h^{-1})*c*(h^{-1})=\\
&e*k*e*k*e*k*f*(k^{-1})*(e^{-1})*(k^{-1})*(e^{-1} )*k=\\
&a*h*(a^{-1})*h*e^{-1}*h*a^{-1} *\\
& \qquad b*f*h*a*h*b^{-1}*h=1.
\end{align*}
\newpage

\subsection{\underline{Group \#8}.	$C_8 \times  C_8$ case}
	For the first step, see table of aut(64) groups. The
	next group in the sequence is a group of order 24,576 =
	$2^{13} * 3$ with the presentation:

\begin{align*}
&a^2=b^2=c^2=d^2=e^2=(a*c)^2=\\
&(a*d)^2=(a*e)^2=(a*f)^2=(c*d)^2=\\
&(c*e)^2=(c*f)^2=(d*e)^2=f^4=\\
&b*c*d*b*d*c=b*c*e*b*e*c=\\
&b*c*( f^{-1} )*b*f*c=d*f^2*d*(f^2)= \\
&e*f^2*e*f^{-2}=d*f*d*f*d*f^{-1}*d*f^{-1}=\\
&d*f*e*(f^{-1})*d*(f^{-1})*e*f= \\
&e*f*e*f*e*(f^{-1})*e*(f^{-1})= \\
&a*b*f^2*b*a*b*( f^{-2} )*b= \\
&a*b*a*b*a*b*a*(f^{-1} )*c*b*(f^{-1})*c=\\
&(a*b*c*b*c*b)^2=(b*c)^6=1. 
\end{align*}
This order 24,576 group has for its automorphism group one of order $2^{28}*3^2*7=16,911,433,728$ with a center $C_2$.
­\subsection{\underline{Group \#14}}
a.	For the automorphism group, see table of aut($g$) for order 64 groups. 
      This is the order 294,912 group. The relations for the
      294,912 group are
\begin{align*}
&a^2=b^2=c^2=d^2=e^2=f^2=h^2=(a*b)^2= \\
&(a*d)^2=(a*e)^2=(a*f)^2=(a*h)^2=(a*k)^2=\\ 
&(b*d)^2=(b*e)^2=(b*f)^2=(c*e)^2=(c*f)^2= \\
&(d*f)^2=(d*h)^2=d*k*d*(k^{-1})=(e*h)^2=(f*h)^2=\\ 
&f*k*f*(k^{-1})=k^4=b*d*k*b*k=d*h*k*h*(k^{-1})=\\
&(a*c)^3=a*c*a*h*c*h=a*c*a*k*c*(k^{-1})= \\
&(b*c)^4=(b*h)^4=(d*e)^4 =(d*e*f*e)^2= \\
&(e*f)^4=a*b*c*a*b*c*b*a*c=d*e*f*e*k*e*f*e*k=\\ 
&d*e*d*f*k^2*e*(k^{-2})*f=(b*c*h*b*h*c)^2=\\
& b*d*e*h*b*k*e*h*d*e*(k^{-1})*e= \\
&c*d*c*e*d*c*e*d*e*c*d*e=\\
&e*k*e*k*e*k*e*(k^{-1} ) *e* ( k^{-1} )*e* (k^{-1} )=1.
\end{align*}
The center of this group is generated by $<d*f*k, c*d*e*(k^{-1} )>$. The next three groups in the
automorphism tower sequence have orders $2^{21}*3^{3}=56,623,104, 2^{24}*3^{3}=452,984,832,$ and
$2^{25}*3^{3}=905,969,664.$\\
\newpage
b. A presentation for the quotient group 294,912/Z is:\\ 
\begin{align*}
&a^2=b^2=c^2=d^2=e^2=(a*b)^2=(a*d)^2=\\
&(a*e)^2=(a*f)^2=(b*d)^2=(c*d)^2=(d*e)^2=\\
&f^4=(a*c)^3=a*c*a*e*c*e=a*c*a*f*c*( f^{-1} )=\\
&b*e*f*e*b*f=b*f^2*b*(f^{-2} )=(b*c)^4=(b*e )^4=\\
&b*f*b*f*b*(f^{-1})*b*(f^{-1} )=a*b*c*a*b*c*b*a*c=\\
&(b*c*e*b*e*c)^2=\\
&b*( f^{-1} )*d*f*d*f*b*( f^{-1} )*d*( f^{-1} )*d*f=\\
&d*f*d*f*d*f*d*( f^{-1} )*d*( f^{-1} )*d*(f^{-1} )=\\
&(d*f*d*(f^{-1}) )^3=\\
&d*f^2*d*f^2*d*(f^{-2} )*d*(f^{-2} )=\\
&a*b*c*d*e*(f^{-1})*d*f*d*f^{-1} *c*b*e=\\
&a*c*f*d*f^2*c*d*f^2*d*f^{-1}*c*f*d*f=1.
\end{align*}

The automorphism group of this group is a complete group of order 147,456. A presentation for this
complete group is:\\

\begin{align*}
&a^2=b^2=c^2=d^2=e^2=(a*b)^2=(a*d)^2=\\ 
&(a*e)^2=(a*f)^2=(b*d)^2=(c*d)^2=\\
&(d*e)^2=d*h*d*(h^{-1})=f^4=(a*c)^3=\\
&a*c*a*e*c*e=a*c*a*f*c*(f^{-1})=\\
&b*e*f*e*b*f=b*f^2*b*(f^{-2})=h^6=\\
&a*b*e*b*h*e*a*(h^{-1})=(a*h*c*(h^{-1}))^2=\\
&(a*h*e*h)^2=(a*(h^{-1})*b*c)^2=(b*c)^4=\\
&b*c*(h^{-2})*b*h^2*c=(b*e)^4=\\
&b*f*b*f*b*(f^{-1})*b*(f^{-1})=(b*h*b*(h^{-1}))^2=\\
&(b*h*c*h)^2=(c*(h^{-1})*e*h)^2=\\
&a*b*c*a*b*c*b*a*c=a*b*c*h*b*c*h*b*e=\\
& a*b*h*f*(h^{-1})*c*h*(f^{-1})*(h^{-1})=\\
& a*d*f*h*e*f*(h^{-1})*e*f=\\
&b*h*f*h*f^2*h*(f^{-1} )*h=1.
\end{align*}
\newpage
This group can be generated by $<e,f,h>$ with the following presentation:
\begin{align*}
&a^2=b^4=a*b^2*a*b^2=c^6=a*b*a*b*a*(b^{-1})*a*(b^{-1})=\\
&a*b*a*b*c^2*a*(b^{-1})*a*(b^{-1})*(c^{-2})=\\
&a*b*a*(b^{-1}) *(c^{-1} )*a*(c^{-1})*b^2*c*a*c=\\
&a*b*a*c*a*c*b^2*(c^{-1})*a*(c^{-1})*(b^{-1})=\\
&a*(b^{-1} )*c*a*(c^{-1} ) *a*c^2*(b^{-1} )*( c^{-1} )*a*c=\\
&a*(b^{-1})*c*b^2*c*a*(b^{-1})*(c^{-1})*b^2*(c^{-1})=\\
&(a*c^2)^4=b^2*c*b^2*(c^{-1})*b^2*c*b^2*(c^{-1})=\\
&(b*c^2*(b^{-1})*c^2)^2=\\
&a*b*a*c^2*b^2*c^2*b^2*c^2*(b^{-1})=\\
&a*b*a*c^2*(b^{-1})*c^2*(b^{-1})*c^2*(b^{-1})*c^2=\\
&a*b^2*c*a*b*(c^{-1})*a*b^2*c*a*b*(c^{-1})=\\
&a*b^2*c*b*c*(b^{-1})*c^2*a*b*(c^{-1})*(b^{-1})*c=\\
&a*b*c*b*(c^{-1})*b*(c^{-1})*b*a*c*(b^{-1})*\\
&\qquad \qquad c*(b^{-1} )*( c^{-1} )=1.
\end{align*}
This complete group has the following order structure:\\

\begin{tabular}{|c|c|c|} \hline 
order of element& number of elements& number of classes\\ \hline
2&2847& 35\\	   
3&1664 &3\\	   
4&45792 &64\\	   
6&41856 &23\\	   
8&16896 &6\\	   
12&32256 &9\\	   
24&6144 &1\\ \hline	 
\end{tabular}
\linebreak\\

\subsection{\underline{Group \#15}} The order 196,608 group in this sequence has the presentation:
\begin{align*}
&c^2=d^2=a^2*b^2=b*d*(b^{-1})*d=(c*d)^2=\\
&a^2*(b^{-4})=a*b*a*b*(a^{-1})*(b^{-1})*(a^{-1})*(b^{-1})=\\
&a*b*(a^{-1} )*d*(a^{-1})*(b^{-1} ) *a*d=\\
&(a*(b^{-1})*c*(b^{-1}))^2=(a*c)^4=\\
&a*d*a*d*(a^{-1} )*d*a^{-1}*d=\\
&a^2*c*a^2*c*(b^{-2})*c*(b^{-2})*c=\\
&a*b*c*b*a*c*( a^{-1} )*(b^{-1} )*c* (b^{-1}) * ( a^{-1} ) *c=\\
&(a*c*(a^{-1})*c)^3=(b*c*(b^{-1})*c)^3=\\
&a*( b^{-1} )*c*(a^{-1} ) *c* (b^{-1} ) *c* ( b^{-1} ) *c* (a^{-1} ) *b*c*a\\
& *c*b*c*b*c=1.
\end{align*}
The center of this group is generated by $<a*d*(a^{-1})*d,(a*c*b*c*b*c)^2>.$\\

This group has 354 conjugacy classes and the order structure: \\

\begin{tabular}{|c|c|c|} \hline 
order of element& number of elements& number of classes\\ \hline
2& 6143& 106\\ 
3& 512& 1\\	   
4& 81920& 181\\	   
6& 36352 &36\\	 
8& 43008 &22\\
12& 28672 &7 \\ \hline
\end{tabular}
\linebreak\\
The automorphism group of this order 196,608 group has order 100,663,296, and a center of order 8 ($C_2\times C_2\times C_2$). The quotient group here has a center of order 32.
\subsection{\underline{Group \#16}}The group of order 294,912 has the following presentation
on five generators:
\begin{align*}
&a^2=b^2=c^2=d^4=e^4=(a*b)^2=(b*c)^2=\\
&(b*d)^2=a*c*d*a*(d^{-1})*c=a*c*e*a*(e^{-1})*c=\\
&(c*d)^3=c*e^2*c*(e^2)=d^2*e*(d^2)*(e^{-1})=\\
&(b*e)^4=(b*e*b*(e^{-1}))^2=(b*e*c*e)^2=\\
&b*e*(d^{-1})*e*b*(e^{-1})*(d^{-1})*(e^{-1})=\\
&c*e*c*e*c*(e^{-1})*c*(e^{-1})=\\
&d*e*d*e*(d^{-1})*e*(d^{-1} )*e=\\
&d*e*d*(e^{-1} )*(d^{-1} )*e*d*( e^{-1} )=\\
&(a*c*a*c*(e^{-1}))^2=\\
&a*e*b*(e^{-1})*a*e^2*b*(e^2)=\\
&a*e*d*e^2*a*e*(d^{-1} ) * ( e^2 ) =\\
&c*d^2*c*e*c*(d^2)*c*(e^{-1})=(a*c)^6=\\
&c*d*e*d*c*e*c*d*e*(d^{-1})*c*e=\\
&a*c*a*d*c*e*c*e*(d^{-1})*(e^{-1})*a*e*c*(e^{-1})=1.
\end{align*}

This group has the following order structure:\\

\begin{tabular}{|c|c|c|} \hline 
order of element& number of elements& number of classes\\ \hline 
2&6463 &79	\\   
3&5120 &3 \\ 
4&106176 &156\\   
6&78848 &20	\\   
8&18432& 4	 \\ 
12&79872& 11\\ \hline
\end{tabular}	 
\linebreak\\

\subsection{\underline{Group \#17}} The following is the second order 3072 group in this automorphism group sequence. The relations are:
\begin{align*}
&a^2=b^2=c^2=d^2=f^2=e^3=(a*b)^2=\\
&a*e*a*(e^{-1})=(b*c)^2=b*e*b*(e^{-1})=\\
&(c*d)^2=(c*f)^2=(d*f)^2=(a*b*d)^2= \\
&(a*b*f)^2=a*c*d*a*d*c=a*c*f*a*f*c=\\
&(c*e)^3=(d*e)^3=d*e*d*f*e*f=\\
&a*c*a*e*c*a*(e^{-1})*c*e*a*c*(e^{-1})=1.
\end{align*}
The next group has order 294,912 = $2^{15} * 3^2$ and has a center of order 2. Its presentation on six generators and class structure are:
\begin{align*}
&a^3=b^2=c^2=d^2=e^4=f^2=(b*d)^2=\\
&d*e*d*(e^{-1} )=a*d*e*(a^{-1} )*( e^{-1 })*d= \\
&a*d*f*(a^{-1} )*f*d=a*e*(a^{-1} )*b*(e^{-1} )*b=\\
&(b*c*(e^{-1}))^2=(e*f)^3=\\ 
&a*b*c*(e^{-1})*c*(e^{-1})*(a^{-1})*b=\\ 
&a*b*e*c*e*c*(a^{-1})*b=(a*c)^4=\\ 
&b*c*b*d*e*c*(e^{-1})*d=(b*f)^4=\\
&(c*d*f*d)^2=(c*f)^4=(d*f)^4=(b*e*f)^3=\\ 
&a*b*( a^{-1} )*c*a*e*b*( e^{-1}) * ( a^{-1}) *c=\\
&a*f*c*f*b*f*c*f*( a^{-1} ) *b=\\
& c*d*c*f*c*f*d*f*c*f =\\
&a*b*a*f*c*b*(a^{-1} )*c*b*f*b=\\ 
&a*d*a*(e^{-1})*f*(e^{-1} )*d*f* (e^{-1 })*d*a*f=1.
\end{align*}
This order 294,912 group has an automorphism group of order 1,179,648 and a center $C_2\times C_2$.\\

\begin{tabular}{|c|c|c|} \hline 
order of element& number of elements& number of classes\\ \hline 
2&4351&39\\
3&6272&3\\
4&89856&84\\
6&83840&22\\
8&36864&11\\
12&73728&14\\ \hline
\end{tabular}
\newpage

\subsection{\underline{Group \#19}} This tower gives an order of $2^{17}*3^2=2,359,296$ for the next group in the tower.

\subsection{\underline{Group \#22}} The order 196,608 = $2^{16} * 3$ in this sequence has the following presentation:
\begin{align*}
&a^4=b^4=c^2=d^2=e^2=f^2=(a*d)^2=\\
&b*c*(b^{-1})*c=(c*e)^2=(c*f)^2=(d*e)^2=\\&
(d*f)^2=(e*f)^2=a^2*c*(a^{-2})*c=\\
&a*d*c*d*( a^{-1} )*c=b*e*f*(b^{-1})*f*e=\\ 
&a^2*b^2*( a^{-2} )*( b^{-2} ) =a*c*e*a*e*c*( b^{-2} )=\\
&a*d*(b^{-2})*d*a^{-1}*b^2=a*e*b^{-2}*e*a^{-1}*b^2=\\
&(a*f*(a^{-1})*e)^2=(a*f*(a^{-1})*f)^2=\\
&a*b^2*(a^{-1})*b*a*(b^{-2})*(a^{-1} )*b=\\
&a*e*f*(a^{-1})*(b^{-1} )*a*f*e*(a^{-1})*b=\\
&a^2*b*a^2*b*(a^{-2})*(b^{-1})*(a^{-2})*(b^{-1})=\\
&a^2*b*a*d*(b^{-1})*(a^{-2})*(b^{-1})*d*(a^{-1})*b=\\
&a^2*b*(a^{-1})*f*(b^{-1})*(a^{-2})*(b^{-1})*(a^{-1})*f*(b^{-1})=\\
&a^2*b^{-1}*a^{-1}*b^{-1}*f*a^{-1}*f*d*b^{-1}*d*b^{-1}=\\
&a^2*e*b*d*a*(b^{-1})*e*(b^{-1} )*(a^{-1})*d*(b^{-1} )= \\
&a*b*a*b*d*(b^{-1} )*d* (b^{-1} )*(a^{-1} ) *(b^{-1} ) * ( a^{-1}) *b =\\
&a*b*a*(b^{-1 })*(a^{-1})*d*c*(b^{-1})*a*b*c*d=\\
&a*b*a*c*(a^{-1} )* (b^{-1} )*(a^{-1} )*( b^{-1} ) *a*c* ( a^{-1} )* b=\\ 
&a*b*( a^{-1} )*b*e*b*a* (b^{-1} ) * (a^{-1} ) * ( b^{-1} )*e* ( b^{-1} )=\\ 
&a*d*b*a*d*(b^{-1} )*d*(a^{-1})*(b^{-1} ) *d*a^{-1}*b=\\ 
&a*e*b*a*e*b*e*( a^{-1} ) * (b^{-1} ) *e* ( a^{-1} ) *b=1.
\end{align*}
This group has 741 conjugacy classes and the order structure:\\

\begin{tabular}{|c|c|c|} \hline 
order of element& number of elements& number of classes\\ \hline 
2&7487&187\\
3&512&1\\
4&86720&425\\
6&48640&63\\
8&36864&48\\
12&16384&16\\ \hline
\end{tabular}
\linebreak\\

\subsection{\underline{Group \#27}} This is the second term in the automorphism tower. 
The presentation for this group of order 294,912 = $2^{15} * 3^2$ is
\begin{align*}
&a^2=b^2=c^2=e^2=f^2=g^2=h^3=(a*c)^2=\\
&a*d*a*d^{-1}=(a*e)^2=(a*f)^2=(a*g)^2=\\
&a*h*a*(h^{-1})=b*d*b*(d^{-1})=(b*e)^2=(c*e)^2=\\
&(c*f)^2=(c*g)^2=c*h*c*(h^{-1})=d^4=\\
&d*e*(d^{-1})*e=d*f*(d^{-1})*f=d*g*(d^{-1} )*g=\\
&d*h*(d^{-1})*(h^{-1})=(e*f)^2=(e*g)^2=(g*h)^2=\\
&a*b*a*b*(d^{-2} )=a*b*a*f*b*f=b*g*(h^{-1})*b*h*g=\\
&(c*d)^3=(e*h)^3=(b*c)^4=(b*c*b*g)^2=\\
&(b*g)^4=(e*(h^{-1} )*f*h)^2=( f*g)^4 =\\
&(f*h*f*(h^{-1}))^2=f*g*f*h*f*h*g*f*(h^{-1})=\\
&a*b*a*g*f*g*b*g*f*g=\\
&a*b*a*(h^{-1} )*f*h*b* (h^{-1} ) *f*h=\\
&b*c*d^2=*c*b*c*(d^{-2})*c=\\
&b*c*d*c*h*b*c*h*b*h*(d^{-1} )*c=1.
\end{align*}
There are 1008 conjugacy classes in this group and the order structure of this group is \\

\begin{tabular}{|c|c|c|} \hline 
order of element& number of elements& number of classes\\ \hline 
2&13375&370\\
3&4256& 3\\
4&117696&445 \\
6&70496&139\\
12&89088&50\\ \hline
\end{tabular}
\newpage

\subsection{\underline{Group \#43}}The automorphism group tower is:\\
 a. aut($g$): order 12,288. 
This version of Hol($C_8\times C_2\times C_2$) is the presentation given in the holomorph article.
arXiv:math/0609571
\begin{align*}
&a^2=b^2=c^2=d^2=e^2=f^2=\\
&(a*b)^2=(a*c)^2=(a*d)^2=\\
&(a*e)^2=(a*f)^2=(b*c)^2=\\
&(b*d)^2=(b*e)^2=(b*f)^2=\\
&(c*d)^2=(c*e)^2=(c*f)^2=\\
&(d*e)^2=(d*f)^2=(e*f)^2=\\
&g^4=h^4=g^2*h*g^-2*h=g*h*g*((h*g*h)^{-1})=\\
&a^{-1}*g^{-1}*a*g=b^g*a*b=c^g*a*b*c=\\
&a^{-1}*h^{-1}*a*h=b^h*a*b*c=c^h*a*c=\\
&d^{-1}*g^{-1}*d*g=e^g*d*e=f^g*d*e*f=\\
&d^{-1}*h^{-1}*d*h=e^h*d*e*f=f^h*d*f=\\
&j^8=j^4*a=\\
&a^{-1}*j^{-1}*a*j=b^{-1}*j^{-1}*b*j=c^{-1}*j^{-1}*c*j=\\
&d^j*a*d=e^j*b*e=f^j*c*f=\\
&g^{-1}*j^{-1}*g*j=h^{-1}*j^{-1}*h*j=\\
&k^2=j^k*j=\\
&a^{-1}*k^{-1}*a*k=b^{-1}*k^{-1}*b*k=c^{-1}*k^{-1}*c*k=\\
&d^{-1}*k^{-1}*d*k=e^{-1}*k^{-1}*e*k=f^{-1}*k^{-1}*f*k=\\
&g^{-1}*k^{-1}*g*k=h^{-1}*k^{-1}*h*k=1.
\end{align*}
A 5-generator presentation for this group obtained with CAYLEY is:
\begin{align*}
&a^2=b^2=c^2=a*d*a*(d^{-1})=a*e*a*(e^{-1})=\\
&b*d*b*d^{-1}=c*d*c*d^{-1}=(d*e)^2=(d*e^{-1})^2=\\
&e^4=(a*b)^3=a*b*c*b*c*a*c= \\
&b*c*b*e*b*c*b*(e^{-1})=(b*e)^4=\\ 
&(b*e*b*(e^{-1}))^2=(c*e)^4=(c*e*c*(e^{-1}))^2=\\
& d^8=b*c*e*c* ( e^{-1} )*b*e*c* ( e^{-1} ) *c=\\
&c*d^4*e^2*c*(e^{-2})=\\ 
&a*b*a*e^2*b*a*( e^{-2} ) *b*(e^{-2} )=1.
\end{align*}
b. Order 98,304 case:
\begin{align*}
&a^2=b^2=c^2=f^2=k^2=x^2=a*d*a*(d^{-1})=\\
&a*e*a*(e^{-1})=(a*f)^2=a*h*a*(h^{-1})=(a*k)^2=\\
&(a*x)^2=b*d*b*(d^{-1})=(b*f)^2=b*h*b*(h^{-1})=\\
&(b*k)^2=c*d*c*(d^{-1})=(c*f)^2=c*h*c*(h^{-1})=\\
&d^4=(d*(e^{-1}))^2=d*f*(d^{-1})*f=d*h*(d^{-1})*h=\\
&d*k*d^{-1}*k=d*x*d^{-1}*x=e^4=(e*f^{-1})^2=\\
&e*h*e^{-1}*h^{-1}=e*k*e^{-1}*k=e*x*e^{-1}*x=\\
&(f*h)^2=(f*k)^2=(f*x)^2=h^4=h*k*(h^{-1}1)*k= \\
&h*x*(h^{-1} )*x= (k*x)^2=(a*b)^3=a*b*a*x*b*x=\\ 
&a*c*a*x*c*x=c*e^2*c*(e^{-2})=a*b*c*b*c*a*c=\\
&b*c*b*e*b*c*b*(e^{-1})=b*c*k*c*b*k*(e^{-2})=(b*e)^4=\\ 
&c*d^2*e*c*(h^{-2})*(e^{-1})=1.
\end{align*}
\subsection{\underline{Group \#44 }} 
The first two groups in this automorphism tower have the presentations:\\

a. aut($g$):\\
\begin{align*} 
&a^2=b^2=c^2=e^2=(a*b)^2=a*d*a*(d^{-1})=(a*e)^2=\\
& b*d*b*(d^{-1})=(b*e)^2=d^4=(d*(e^{-1}))^2=(a*c)^3= \\
&c*d^2*c*(d^{-2})=(b*c)^4=(b*c*d*c)^2=\\ 
&c*d*c*d*c*(d^{-1})*c*(d^{-1})=(c*d*c*e)^2=\\ 
&(c*e)^4=a*b*c*a*b*c*b*a*c=\\ 
&b*c*b*e*c*b*e*c*e*b*c*e=1.
\end{align*}
\newpage
b. Second member. This group has order 393,216 = $2^{17}*3$ and center $C_2 \times  C_2$.
\begin{align*}
&a^2=b^2=e^2=f^2=h^2=(a*b)^2=a*c*a*(c^{-1})=\\
&(a*e)^2=(a*f)^2=(a*h)^2=b*c*b*(c^{-1})=(b*e)^2=\\
&(b*f)^2=c^4=c*e*(c^{-1})*e=d^4=d*f*(d^{-1})*f=\\
&(f*h)^2=b*d^2*b*(d^{-2})=c^2*d*(c^{-2})*(d^{-1})=\\
&c^2*f*(c^{-2})*f=c*e*(d^{-1} )*e*(c^{-1})*d=\\
&c*e*h*(c^{-1})*e*h=b*d*h*d^2*h*(d^{-1} )=\\
&a*b*d*a*h*(d^{-1} )*b*h=( a*b*h*d)^2=\\
&a*(d^{-1})*c*d*a* (d^{-1} ) * ( c^{-1} )*d=\\
&c*d*c*d*(c^{-1})*d*(c^{-1})*d=\\
&c*d*c*(d^{-1})*(c^{-1})*d*(c^{-1})*(d^{-1})=\\
&c*d*c^{-1}*f^{-1}*c^{-1}*d^{-1}*c*f^{-1}=\\
&c*f*c*f*(c^{-1})*f*( c^{-1} )*f=\\
&c*f*(c^{-1})*h*(c^{-1})*f*c*h=(c*f*e*f)^2=\\ 
&c*h*c*h*( c^{-1} )*h*( c^{-1})*h=(e*f)^4=\\ 
&c*d*(c^{-1} ) *(h^{-1} ) *( c^{-1} )*(d^{-1} ) *c*d*(h^{-1} ) *(d^{-1} )=\\
&a*d*a*d*a*d*a*(d^{-1})*a*(d^{-1})*a*(d^{-1})=1.
\end{align*}

\subsection{\underline{Group \#81}} The group of order 6144 in this tower has the following presentation:
\begin{align*}
&a^2=b^2=c^2=d^2=(a*d)^2=(a*b)^3=(b*c)^3=\\
&(b*d)^3=(a*c)^4=(a*c*d*c)^2=(c*d)^4=\\
&(a*b*a*c)^3=a*b*a*c*a*c*b*a*d*b*c*a*c*a*b*d=\\
&a*b*a*c*b*d*c*a*d*b*c*a*b*d*c*d=\\
&a*b*a*c*d*c*b*a*d*b*c*d*c*a*b*d=1.
\end{align*}
The second group in this sequence has order 24,576 = $2^{13} * 3$ and the presentation:
\begin{align*}
&a^2=b^2=c^2=d^2=(a*d)^2=a*e*a*(e^{-1})=\\
&b*e*b*(e^{-1} )=c*e*c* (e^{-1} )=e^4=( a*b)^3=\\
&(b*c)^3=(b*d)^3=(a*c)^4=(a*c*d*c)^2=\\
&(c*d)^4=(d*e)^4=(d*e*d*(e^{-1}))^2=\\
&a*c*a*c*d*e^2*d*(e^{-2} )=\\
&a*c*a*d*e*d*c*d* ( e^{-1} )*d=( a*b*a*c)^3=\\
&b*c*b*d*c*b*(e^{-1})*d*e*b*c*d=\\
&a*b*a*c*b*d*c*a*d*b*c*a*b*d*c*d=1.
\end{align*}
The next factor has order $49,152 = 2^{14}*3$ and presentation:
\begin{align*}
&a^2=b^2=c^2=d^2=f^2=(a*d)^2=a*e*a*(e^{-1})=\\
&b*e*b*(e^{-1})=e^4=e*f*(e^{-1})*f=(a*b)^3=\\
&(a*b*f)^2=a*c*a*e*c*(e^{-1})=(a*c*f)^2=\\
&a*e^2*f*a*f=(b*c)^3=(b*d)^3=(a*c)^4=\\
&(a*c*d*c)^2=( c*d)^4=c*d*e*d*c*d*(e^{-1} )*d=\\
&c*d*e^2*d*(e^{-1})*c*(e^{-1})=(d*e)^4=(d*f)^4=\\
&a*c*d*f*d*c*d*a*f*d=b*c*e*c*b*d*e*f*d*f=\\
&b*c*b*d*c*b*e*d*(e^{-1})*b*c*d=(b*d*f)^4=1.
\end{align*}
The last group in this sequence is a complete group of order $98,304 = 2^{15}*3$ and presentation:
\begin{align*}
&a^2=b^2=c^2=d^2=(a*d)^2=a*e*a*(e^{-1})=\\
& b*e*b*(e^{-1})=c*e*c*(e^{-1})=e^4=e*f*e*(f^{-1})=\\
&(a*b)^3=a*b*f*a*b*(f^{-1})=a*c*f*c*a*( f^{-1} )=\\
&a*e^2*f*a*f=a*f*a*(f^{-3})=(b*c)^3=(b*d)^3=\\
&(a*c)^4=(a*c*d*c)^2=(c*d)^4=(d*e)^4=\\
&(d*e*d*(e^{-1}))^2=(d*f)^4=a*c*a*c*d*e^2*d*(e^{-2})=\\ 
&a*c*a*d*e*d*c*d*(e^{-1})*d=\\
&b*d*b*e*d*b*( f^{-1})*d* (e^{-1} ) *( f^{-1} )=\\ 
&c*d*f*d*c*e*d*f*d*(e^{-1})=\\
& b*c*b*d*c*b*d*b*f*c*d*( f^{-1} )=1.
\end{align*}

\subsection{\underline{Group \#82}} A presentation for group number 82's automorphism group is:\\

 $ a^2= b^2= c^2= d^3= b*c*b*c= e^4= a*e^{-1}*a*e^{-1}=$\\
 $ a*c*a*c= c*e^2*c*e^{-2}= d^{-1}*c*d^{-1}*c*d^{-1}*c=$\\
 $ d^{-1}*e*d^2*e*d^{-1}*e= e*b*e*b*d*e^{-1}*d^{-1}=$\\
 $ d*e^{-1}*d^{-1}*b*e*b*e= b*a*d^{-1}*a*b*a*d^{-1}*a=$\\
 $ c*d^{-1}*a*d*c*d^{-1}*a*d= d^{-1}*c*d*e*d^{-1}*c*d*e^{-1}=$\\
 $ c*d*a*d^{-1}*c*d*a*d^{-1}= b*d*b*d*b*d*b*d=$\\
 $ a*d*a*d^{-1}*a*d*a*d^{-1}=c*e^{-1}*c*e^{-1}*c*e^{-1}*c*e^{-1}=$\\
 $ a*d*e*d^{-1}*a*d*e^{-1}*d^{-1}= b*a*b*a*b*a*b*a=$\\
 $ c*e^{-1}*c*d^{-1}*c*e^{-1}*c*d^{-1}*e^{-1}*d^{-1}=c*d^{-1}*e^2*d*c*d^{-1}*e^2*d=$ \\
 $ c*e*c*d^{-1}*c*e*c*d^{-1}*e*d^{-1}=b*d*a*d^{-1}*a*b*d*a*d^{-1}*a=$\\
 $ e*d^{-1}*a*d*a*e*d^{-1}*a*d*a= $\\
 $ d^{-1}*c*b*d*c*b*a*d^{-1}*a*c=$\\
 $ c*d*a*d*a*c*d^{-1}*c*a*d*a*d=$\\
 $ d*a*d^{-1}*a*e^{-2}*d*a*d^{-1}*a*e^{-2}=$\\
 $ b*d*c*e*c*d^{-1}*b*d*c*e^{-1}*c*d^{-1}=$\\
 $ b*a*b*d^{-1}*c*a*d*b*a*b*d^{-1}*c*a*d=$\\
 $ a*d*a*d*e^2*d*a*d*a*d*e^{-2}*d =1.$

This automorphism group has the structure $[512 @ C_3]@D_4$. All but five of the normal subgroups of this automorphism group are 2-groups. Those that are not have orders 1536 (one case), 3072 (one case), and three cases with orders 6144. The quotient groups aut($g$)/(order 512 groups) are either small group (24,12), i.e., $S_4$ (seven cases), or $(24,8)=(4,6|2,2)$ or $C_3@D_4$ (one case). A presentation for the order 1536 group is:
\begin{gather*}
   a^2= g^2= f^3= c^4= b^{-1}*e*b*e= c^{-2}*d^2= b^4= c*d*c*d^{-1}.=e^2*b^2=\\ 
  a*c*a*c^{-1}= b*c*b^{-1}*c^{-1}= a*b*a*b^{-1}=a*e*a*e^{-1}= \\
a*d*a*d^{-1}= b^{-1}*d^{-1}*b^{-1}*d^{-1}=\\
  d^{-1}*e^{-1}*d^{-1}*e^{-1}= a*f*a*f^{-1}= c*f*b^{-1}*f^{-1}=\\
  c^{-1}*e^{-1}*c^{-1}*e^{-1}= a*g*a*g= b*g*b^{-1}*g=c*a*g*c^{-1}*g=\\
  d*a*g*d^{-1}*g= e*a*g*e^{-1}*g= c*b*f*c*f^{-1}=\\
  f^{-1}*a*d^{-1}*f*e^{-1}*b^{-1}= b*f*d*f^{-1}*d*e=\\
f^{-1}*g*f^{-1}*g*f^{-1}*g =1.  
 \end{gather*}
The sylow 2-subgroup of this order 1536 group, according to GAP, is generated by the 
following presentation:
\begin{gather*}
    a^4 = b^4 = c^2 = d^2 = e^2 = f^2 = g^2 = (a, b) =\\ 
    (b^{-1} * d)^2 = (c * d)^2 = a^{-1} * e * a * e = \\
    b^{-1} * e * b * e =  (c * e)^2 = (d * e)^2 =  \\
    a^{-1} * f * a * f = (c * f)^2 =  (d * f)^2 = \\
    (e * f)^2 = (c * g)^2 = (d * g)^2 = (e * g)^2 = \\
    (f * g)^2 = a * c * a^{-1} * b^2 * c = \\
    a^{-1} * d * a^{-1} * b^2 * d = g * e * a^{-1} * g * a =\\ 
    c * b^{-1} * a^2 * c * b =  f * e * b^{-1} * f * b = \\
    g * e * b^{-1} * g * b = 1.
\end{gather*}
The automorphism group of this order 1536 group is $589,824 = 2^{16}*3^2$. The next one in the automorphism tower sequence has order 1,179,648 and is a complete group. We were unable to generate presentations for these last two groups using GAP. The automorphism group orders were computed with MAGMA.\\
\newpage

\subsection{\underline{Group \#93}} Below are the presentations for the two members of this sequence.\\

a. order 3072:
\begin{align*}
&a^3=b^3=c^3=d^3=(a*(c^{-1}))^2=(b*(c^{-1}) )^2=\\
&(b*d^{-1})^2=(c*d^{-1})^2=(a*b)^3=(b*d*c)^2=\\
&(a*(b^{-1})*a*(d^{-1}) )^2=\\
&a*(b^{-1})*c*a*(b^{-1})*(c^{-1})*(a^{-1})*b=\\
&(a*d^{-1})^4=a*b*c*d^{-1}*a*d*a*d^{-1}*b^{-1}=1.
\end{align*}

b. order $294,912 = 2^{15}*3^2$:\\
\begin{align*}
&a^4=b^3=c^6=d^3=e^4=a*b*( a^{-1}) * ( b^{-1}) =\\
&a*e*(a^{-1})*(e^{-1})=a^2*c*(a^{-2})*(c^{-1})=\\ 
&a*c*a*c*( a^{-1} )*c=a*d*( a^{-1} )*( c^{-1} )*d*c*d=\\ 
&a*c*a*( c^{-1} )* (a^{-1} )*c*a*( c^{-1} )=\\ 
&a*c*e*c*a*(c^{-1} )*( e^{-1} )*( c^{-1} )=\\ 
&a*c*(e^{-1} )*c*a*( c^{-1} )*e*( c^{-1})=\\
&(a*d*a*d^{-1})^2=(a*d*a^{-1}*d^{-1})^2=\\
&a*e^2*c^{-1}*e*c*e*c^{-1}=b*c^3*b^{-1}*c^{-3}=\\
&b*(c^{-1})*(d^{-1})*(e^{-1})*(b^{-1})*(e^{-1})*c*(d^{-1})=\\
&b*d*b*( d^{-1}) * ( b^{-1})*( c^{-1}) * ( d^{-1}) *c=\\
&b*d*(b^{-1})*e*(d^{-1})*b*d*(e^{-1})=\\
&b*d*(b^{-1})*(e^{-1})*(d^{-1})*b*d*e=\\
&(b*(d^{-1}) )^4=\\
&b*(d^{-1})*(c^{-1})*e*(b^{-1})*e*(d^{-1})*c=\\
&b*e^2*b^{-1}*d^{-1}*e^{-2}*d=(b*e^{-1})^4=\\
&c*e*d*(e ^{-1})* (c^{-1}) * ( e^{-1}) * (d^{-1}) *e=\\
&c*e^{-1}*d*e*c^{-1}*e*d^{-1}*e^{-1}=\\
&(d*(e^{-1}) )^4=1.
\end{align*}
\newpage
\subsection{\underline{Group \#103}} The following are the higher-order presentations for the groups in this automorphism tower.\\

Group of order $73,728 = 2^{13} * 3^2$:
\begin{align*}
&a^4=b^4=c^4=d^2=b^2*(e^{-2})=a*d*(a^{-1})*d=\\
&a*e*(a^{-1}) *(e^{-1} )=c*d*( c^{-1} ) *d=a^2*b*a^2*b=\\
&a^2*c*a^2*c=a^2* (e^{-1} )*d*( e^{-1} )*d=\\ 
&a*c*a*(c^{-1})*(a^{-1})*(c^{-1})= \\
&a*e*c*(e^{-1} )*(a^{-1} )*( c^{-1})=b^2*d*b^2*d=\\ 
&a*b*a*b*( a^{-1})* (b^{-1}) *( a^{-1}) * ( b^{-1})=\\
&a*b*c*b*a^{-1}*b^{-1}*c^{-1}*b^{-1}=\\
&a*b*d*b*( a^{-1} )*(b^{-1})*d* (b^{-1} )=\\
& a*b*e*b*(a^{-1} )*(b^{-1})*(e^{-1} )*(b^{-1} )=\\
& b*c*b*c* ( b^{-1} ) * ( c^{-1} ) * (b^{-1} )*c=\\
&b*d*b*d*(b^{-1})*d*(b^{-1})*d=\\ 
&b*d*(b^{-1} )*(e^{-1} ) *(b^{-1} )*d*b*e=\\
&b*e*b*e*(b^{-1})*(e^{-1})*(b^{-1})*(e^{-1})=\\ 
&a*(b^{-1} )*( e^{-1} ) * ( b^{-1} ) * ( c^{-1} ) * ( a^{-1} ) *b*e*b*c=\\
& b*c^2*(b^{-1})*d*(b^{-1})*c^2*b*d=\\ 
&a*b*d*b*c*a*b*d* (b^{-1} ) *c*a* ( b^{-1} ) *d* (b^{-1} ) * ( c^{-1} )=1.
\end{align*}
The number of classes in this group is 119 and its order structure is:\\

\begin{tabular}{|c|c|c|}\hline
order of element& number of elements& number of classes\\ \hline 
2&2127&28\\
3&1088&2\\
4&22704&56\\
6&14784&16\\
8&16128&5\\
12&16896&11\\ \hline
\end{tabular}
\linebreak\\
\newpage
The next group in this tower is a group of order $147,456 = 2^{14} * 3^2$
with the presentation:
\begin{align*}
&a^4=b^4=c^2=d^4=b^2*e^{-4}=a*c*a^{-1}*c=\\
&a^2*b*a^{-2}*b=a^2*d*a^{-2}*d=a*b^2*a^{-1}*b^{-2}=\\
&(a*(b^{-1})*(e^{-1}))^2=\\
&a*d*a*(d^{-1})*(a^{-1})*(d^{-1} )=\\
&(a*e*b^{-1})^2=a*e^2*a^{-1}*e^{-2}=\\
&a*(e^{-1} ) *b*e*(a^{-1} ) *(b^{-1} )=b^2*c*b^2*c=\\
&b^2*d*b^2*d=b^2*e*b^2*(e^{-1})=c*d^2*c*d^2=\\
&a*b*a*b*(a^{-1})*(b^{-1})*(a^{-1})*(b^{-1})=\\
&a*b*c*b*(a^{-1})*(b^{-1})*c*(b^{-1} )=\\
&a*b*d*b*a^{-1}*b^{-1}*d^{-1}*b^{-1}=\\
&a*(d^{-1})*b*e*(d^{-1})*(b^{-1})*a*(e^{-1})=\\ 
&b*c*b*c*(b^{-1})*c*(b^{-1})*c=\\
& b*c*(b^{-1})*(d^{-1})*(b^{-1})*c*b*d=\\
& b*d*b*d*(b^{-1})*(d^{-1})*(b^{-1})*d=\\
&b*e*c*d^2*( e^{-1} ) * (b^{-1} ) *c=\\ 
&c*d*c*d*c*(d^{-1} ) *c* (d^{-1} )=\\
&a*c*d*a*c* (d^{-1} ) *c* ( a^{-1} ) * (d^{-1} )=1.
\end{align*}
The number of classes in this group is 106 and it has the following order structure:\\

\begin{tabular}{|c|c|c|}\hline
order of element& number of elements& number of classes\\ \hline 
2&2383&20\\
3&1088&2\\
4&37296&41\\
6&26048&18\\
8&42240&10\\
12&32256&12\\
24&6144&2\\ \hline
\end{tabular}
\newpage

The last group in this tower has order $294,912 = 2^{15} * 3^2$ and presentation:
\begin{align*}
&b^2= e^2= a^4=a*b*a^{-1}* b=a*d*a^{-1}* d^{-1}= \\
&a*e*a^{-1}* e=a*h*a^{-1}* h^{-1}=b*c*e*c^{-1}=\\
&b*c^{-1}*e*c=(b*e)^2=b*g*b*g^{-1}=e*h*e*h^{-1}=\\
&g^4=a^2*c*a^{-2}*c=a^2*d^{-1}*b^{-1}*d^{-1}*b^{-1}=\\
&a^2*g*a^{-2}*g=a*d^{-1}*g^{-1}*a^{-1}*d*g=\\
&a*g*a*g^{-1}*a^{-1}*g^{-1}=a*g*c^{-1}*h^2*c^{-1}=\\
&a*h*g*h^{-1}*a^{-1}*g=e*g^2*e*g^{-2}=(a*c)^2*(a^{-1}*c^{-1})^2= \\
&a* c*h^{-1}*c^{-2}*h^{-1}*c^{-1}*g^{-1}=\\
&b* h*c^{-1}*g*e*g*c*h^{-1}=(d*h)^2 *(d^{-1}*h^{-1})^2=\\
&(e*g)^2 *(e*g^{-1})^2=a*c*b*h^{-1}*b*h^{-1}*c*e*g=\\
&a*c*h^{-1}*d^{-1}*h^{-2}*d^{-1}*h^{-1}*c*a^{-1}=\\
&a*c^{-1}*h*c*d*a^{-1}*c^{-1}*h^{-1}*d^{-1}*c^{-1}=\\
&c*d*h^{-1}*c^{-1}*h^{-1}*c^{-1}*d*h*c^{-1}*h=\\
&a*c*d*h*d*g*h*c^{-1}*h^{-1}*d^{-1}*h=1.
\end{align*}
The number of classes in this group is 152, and it has the order structure:\\

\begin{tabular}{|c|c|c|}\hline
order of element& number of elements& number of classes\\ \hline 
2&3087&21\\
3&1088&2\\
4&60144&50\\
6&40896&33\\
8&82176&17\\
12&76800&23\\
16&18432&1\\
24&12288&4	\\ \hline
\end{tabular}
\newpage

­\subsection{\underline{Group \#108}} For this automorphism sequence we have found the following three groups with the following presentations:\\

a. order 12,288, and center $C_2$:
\begin{align*}
&a^2=b^4=c^2=a*b*a*b*a*(b^{-1})*a*(b^{-1})=\\
&a*c*b^2*c*a*c*b^2*c=\\
&b^2*c*(b^{-1} )*c*b^2*c*(b^{-1} )*c=\\
&a*b^2*a*b*c*a*c*b*c*a*c=\\
&a*b^2*c*b*c*a*c*( b^{-1})*c*b^2=\\
&b*c*b*c*b*c*(b^{-1})*c*(b^{-1} ) *c* (b^{-1} ) *c=\\
&a*b*a*b^2*c*a*c*a*(b^{-1})*a*c*a*c=\\
&a*b*c*b*a*b*c*b*a*b^{-1}*c*b^{-1}*a*b^{-1}*c*b^{-1}=1.
\end{align*}
This group has 93 classes and the order structure:\\

\begin{tabular}{|c|c|c|}\hline
order of element& number of elements& number of classes\\ \hline 
2&847&32\\	   
3&128&1\\   
4&5040&42\\	   
6&2432&10\\   
8&2304&4\\   
12&1536&3\\ \hline
\end{tabular}	 
\linebreak\\

b. The next group has order $49,152 = 2^{14}*3$. The center is elementary abelian 
of order 8. This group appears to be a direct product of $C_2 \times  C_2$ and a 12,288 order group. The order 49,152 group has the following presentation:
\begin{align*}
&a^2=b^4=c^2=d^2=e^2=f^2=(a*d)^2=(a*e)^2=\\
&(a*f)^2=b*d*(b^{-1})*d=b*e*(b^{-1})*e=\\
&b*f*(b^{-1})*f=(c*f)^2=(d*e)^2=(d*f)^2=\\
&(e*f)^2=a*b^2*a*b^2=a*b*a*b*a*b^{-1}*a*b^{-1}=\\
&a*c*b*c*a*c*(b^{-1})*c=(b*c*e*c)^2=(c*e)^4= \\
&a*c*d*b^2*c*a*c*d*c=a*c*d*c*a*e*c*d*c*e=\\ 
&a*c*e*b^2*c*a*c*e*c=\\
&b^2*c*(b^{-1})*c*b^2*c*(b^{-1})*c=\\
&b^2*c*d*c*b^2*c*d*c=(a*b*a*c*a*c)^2= \\
&a*c*a*d*c*a*d*c*d*a*c*d=\\ 
&b*c*b*c*b*c*(b^{-1} )*c*(b^{-1}) *c*(b^{-1} )*c=\\
&a*b*a*b*c*a*b*a*b*c*d*c*d*c=\\
& a*b*a*c*b*a*c*a*d*c*b*e*c*e*(b^{-1})*c*d*c=1.
\end{align*}
This group has 464 classes with the order structure:\\

\begin{tabular}{|c|c|c|}\hline
order of element& number of elements& number of classes\\ \hline 
2&3583&163\\
3&128 &1\\
4&23040&228\\
6&10112&43\\
8&144	&16\\
12&6144&12\\ \hline
\end{tabular}
\linebreak\\
The automorphism group of the above order 49,152 group has order $2^{27}*3^2=1,207,959,552$. GAP did not return/find the center of this order $2^{27}*3^2$ group.

b. A representation for the order 12,288 group that appears in the 
above 49,152 order group is:
\begin{align*}
&a^2=b^2=c^2=d^2=e^3=f^2=(a,b)=(a,d)=(b,d)=\\
&(c,e)=(c,f)=(a*b*f)^2=(a*e*f)^2=\\
& a*(e^{-1})*f*a*f*e=b*e*f*b*f*(e^{-1} )=\\ 
&(d*e)^3=d*e*f*d*f*(e^{-1})=(a*c)^4=\\
& (a*c*b*c)^2=(a*c*d*c)^2=(b*c)^4=(c*d)^4=\\ 
&a*b*c*d*c*b*c*d*(e^{-1})*c*a*e=\\
& a*c*a*d*c*e*d*c*f*d*f*c*d*(e^{-1})=1.
\end{align*}

c. The automorphism group of the order 12,288 group in the order 49,152 group 
has order $2^{16}*3$ = 196,608 and has for a center the elementary abelian group of 
order 8. A presentation for this group is

\begin{align*}
&a^2=b^2=c^2=d^3=e^2=f^2=g^2=h^2=(a*c)^2=\\
&(a*e)^2=b*d*b*(d^{-1})=(b*e)^2=(b*f)^2=(b*g)^2=\\
&(b*h)^2=(c*h)^2=(e*g)^2=a*b*g*a*g*b=\\
&a*d*a*e*(d^{-1})*e=a*d*a*g*(d^{-1})*g=\\
&a*d*a*h*(d^{-1})*h=a*d*f*a*f*(d^{-1} )=\\
&(a*(d^{-1})*f)^2=a*f*a*g*f*g=a*f*a*h*f*h=\\
&b*c*b*g*c*g=(c*d)^3=c*d*f*c*f*(d^{-1} )=\\
&c*d^{-1}*e*d*c*e=(a*b)^4=(a*b*a*h)^2=\\
&a*b*e*f*a*e*f*b=(a*h)^4=(a*h*e*h)^2=\\
&(b*c)^4=(e*h)^4=(e*h*g*h)^2=(g*h)^4=\\
&a*b*h*g*h*a*h*g*h*b=\\
&b*c*b*d*c*b*f*c*f*b*c*d^{-1}=1.
\end{align*}
This group has 864 classes and an order structure of:\\

\begin{tabular}{|c|c|c|}\hline
order of element& number of elements& number of classes\\ \hline 
2&	8575 &242\\	   
3&	128 &1	\\   
4&	97920 &465	  \\ 
6&	28544 &79	   \\
8&	24576 &40	   \\
12&	36864 &36	 \\ \hline
\end{tabular}
\linebreak\\

­\subsection{\underline{Group \#144}} The complete group factor in this tower.
\begin{align*}
&a^2=b^2=c^2=d^2=(a*b)^3=(a*d)^3=(b*c)^4=\\
&(a*b*d*b*d)^2=(a*b*a*d)^3=(a*b*c)^4=(c*d)^6=\\
&a*b*c*a*b*c*d*b*a*c*b*a*c*d=\\ 
&a*b*c*a*b*d*b*a*c*b*a*d*c*d=\\
&(a*b*c*a*d*a*c)^2=(a*d*b*c)^4=1. 
\end{align*}
­\subsection{\underline{Group \#153}} The first group in this automorphism tower has order 10,752 and presentation:
\begin{align*}
&a^3=b^3=c^3=(a*(b^{-1}) )^2=\\
&a*b*a*c*b*( a^{-1})*c*(a^{-1} )*b*c=\\
&a*c*a*c*(a^{-1})*( c^{-1} )*(a^{-1})*c*( a^{-1} )*( c^{-1} )=\\
&a*b*( c^{-1})*b*c*a*c*(b^{-1})*( c^{-1} )*( a^{-1}) *c=\\
&a*c*a^{-1}*b^{-1}*c^{-1}*a^{-1}*c^{-1}*b*a*c*b*c^{-1}=1.
\end{align*}
This group has 48 classes and the order structure:\\

\begin{tabular}{|c|c|c|}\hline
order of element& number of elements& number of classes\\ \hline 
2&	511&	29\\
3&	896&	2\\
6&	6272& 14\\
7&	3072&	2\\ \hline
\end{tabular}
\linebreak\\
The next group in this tower has order 64,512 and presentation:
\begin{align*}
&a^4=b^3=c^6=(a*b*a*(b^{-1}))^2=\\
&(a*b*(a^{-1})*(b^{-1}))^2=\\
&a*c*a*c*( a^{-1} ) * ( c^{-1} )*( a^{-1} )*( c^{-1} )=\\
&a*c* (a^{-1} )*c*a*(c^{-1})*a*( c^{-1} )=\\
&a*c^3*(a^{-1})*c^3=a^2*b*a^2*b*a^2*b=\\
&a*b*c*( b^{-1})*c*( b^{-1})*( a^{-1})*b^{-1}*c^{-1}*b=\\ 
&a*( b^{-1})*a*c*( a^{-1})*( c^{-1})*b*c*( a^{-1})*c^{-1}=\\ 
&a^2*b*c*b*c*( a^{-1} )*(b^{-1} )*( a^{-1})* (b^{-1})*c=1.
\end{align*}
This group has 48 classes and the order structure:\\

\begin{tabular}{|c|c|c|}\hline
order of element& number of elements& number of classes\\ \hline  
2&	703& 11\\	   
3&	2816& 5\\	   
4&	1344& 4\\	   
6&	30464& 17\\	   
7&	3072& 2	\\   
12&	10752& 4	  \\ 
14&	9216 &2	   \\
21&	6144 &2\\ \hline
\end{tabular}	 
\linebreak\\

\subsection{\underline{Group \#156}} A permutation representation and a presentation for the group $\left[\text{96 
number 2301}\right] \wr C_2$ is
\begin{align*}
&A = (1,2,3,4)(5,7)\\
&B = (1,2)(5,8,6,7)\\
&C = (1,9) (2,10) (3,11) (4,12) (5,13) (6,14) (7,15) (8,16)\\
&A^4=B^4=C^2=(A*B)^3=(A*(B^{-1})^3=\\
&A^2*B^2*( A^{-2} )* ( B^{-2} )=\\
&A*C*A*C*(A^{-1})*C*(A^{-1})*C=\\
&A*C*B*C*(A^{-1} )*C*(B^{-1} ) *C=\\
&B*C*B*C* (B^{-1} )*C*(B^{-1} )*C=1.
\end{align*}
The automorphism group of this group has order
$110,592=2^{12}* 3^3$ and a presentation on four generators:
\begin{align*}
&a^4=b^2=c^4=d^6=a*c^2*a^{-1}*c^2=\\
&a*d*b*(d^{-1})*(a^{-1})*b=a*d^2*(a^{-1})*(d^{-2})=\\
&(a*(d^{-1})*b^{-1})^2=b*c*(d^{-1})*b*d*(c^{-1})=\\
&a*b*a*b*(a^{-1})*b*(a^{-1})*b=a*b*c^2*b*( c^{-1})*( a^{-1} ) *c=\\
&a*b*d*c^2*(a^{-1})*b*d=a*c*a*c*a*( c^{-1})*a*( c^{-1} )=\\
&a*d^{-1}*c*d*a^{-1}*d^{-1}*c^{-1}*d=\\
&a^2*d^{-1}*c^{-1}*d*c*a^{-1}*b*d^{-1}*c^{-1}*a^{-1}*b=1.\\
\end{align*}
\newpage
­\subsection{\underline{Group \#173}} This automorphism tower has had the following members' presentations determined:\\

a. Group of order 6144:\\

\begin{align*}
&a^4=b^4=c^4=d^2=a^2*b*a^2*b=a^2*c*b^2*c^{-1}=\\
&a^2*(c^{-1})*(b^2)*c=a^2*d*(a^2)*d=a*b^2*(a^{-1})*(b^2)=\\	
&b^2*d*b^2*d=(b*c^{-1})^3=c^2*d*c^2*d=\\
&a*b*a*b*(a^{-1})*( b^{-1})*( a^{-1})*b=\\
&a*b*c*b*(c^{-1})*( a^{-1})*b*c=a*c^2*a*(c^2)*(b^2 )=\\
&a*d*a*d*( a^{-1})*d*(a^{-1})*d=\\
&b*c*(b^{-1})*(d^{-1})*(b^{-1})*(c^{-1})*b*(d^{-1})=\\
&b*d*b*d*(b^{-1} )*d* (b^{-1})*d=c*d*c*d*( c^{-1} )*d* ( c^{-1})*d=\\
&a^2*b*c*d*( c^{-1} )*d*b*d=\\
&a*b*d*a*(c^{-1} )*a*c*d*( a^{-1}) *b=\\
&a*b*( d^{-1})*b*(a^{-1}) *( d^{-1})*( a ^{-1}) * ( b^{-1}) * ( d^{-1}) *\\
& \qquad 	(b^{-1})*a*(d^{-1})=1.
\end{align*}
This group has 44 classes, a trivial center, and the following order structure:\\

\begin{tabular}{|c|c|c|}\hline
order of element& number of elements& number of classes\\ \hline  
2&	399	&20\\
3&	512	&1\\
4&	2544	&16\\
6&	1536	&3\\
8&	1152	&3\\ \hline
\end{tabular}
\newpage

b. The order 49,152 group has the following presentation:
\begin{align*}
&a^4=b^4=c^4=e^4=f^2=g^2=\\
&a*b*(a^{-1})*(b^{-1})=a*d*(a^{-1})*(d^{-1})=\\
&a*g*(a^{-1} )*g=b*g*(b^{-1})*g=\\
&c^2*(d^2 )=c*d*( c^{-1} )*(d^{-1} )=c*g*( c^{-1}) *g=\\
&d*g*(d^{-1})*g=e*g*(e^{-1})*g=a^2*c*(a^2 ) *c=\\
&a^2*(e^{-1} )*( c^2)*e=a^2*f*(a^2)*f=\\
&a^2*( g^{-1})*( f^{-1} )*( g^{-1} ) *( f^{-1} )=\\
&a*e*(a^{-1} )*b^2*e=b^2*c*(b^2)*( c^{-1} )=\\
&b^2*d*(b^2)*d=(b*(d^{-1})*f)^2=\\
&b*f*(b^{-1})*e^2*f=c^2*f*(c^2)*f=(c*(e^{-1}))^3=\\
&a*c*e*c*e^{-1}*a^{-1}*c*e=\\
&a*( c^{-1} ) *b* ( e^2) * (a^{-1}) *b*c=\\
&a*( f^{-1})*a*(f^{-1})*(d^{-1})*(b^{-1})*(d^{-1} )*(b^{-1} )=\\
&b*c*(b^{-1})*(e^{-1})*(b^{-1})*(c^{-1})*b*(e^{-1})=\\
& b*(c^{-1})*e*c*b*(c^{-1})*(e^{-1})*(c^{-1})=\\
&b*e*c*d*b*e*(d^{-1})*( c^{-1} )= \\
&b*e*(f^{-1})*e*(b^{-1})*c*(f^{-1})*(c^{-1})=\\ 
&c*e*(c^{-1})*(f^{-1})*(c^{-1})*(e^{-1})*c*(f^{-1})=\\
&c*f*c*f*( c^{-1} )*f*(c^{-1})*f=\\
&a^2*c*e*f*( e^{-1} ) *f*c*f=1.
\end{align*}
\newpage

c. The next member in this tower has order $2^{18}*3 = 786,432$. The derived group of this group has order 49,152, but is not isomorphic to the previous order 49,152 group. A presentation for this derived group is:
\begin{align*}
&a^4=b^4=a^2*(e^{-2})=a*b*(a^{-1})*(b^{-1})=\\
&(a,e)=(b*e^{-1})^2=\\
&a^2*( c^{-1} )*b^2*c=a^2*(d^{-1}) *b^2*d=\\
&a^2*c*b^2*( c^{-1})*b^2=a^2*d*b^2*(d^{-1})*b^2=\\
&(a*c*(a^{-1})*(c^{-1}))^2= \\
&a*c*d^2*a*(d^{-2})*(c^{-1})=\\ 
&a*c*(d^{-1})*b*c*(a^{-1})*(d^{-1})*(b^{-1})=\\
&(a*c*(e^{-1})*(c^{-1}) )^2=\\ 
&a*(c^{-1})*b*c*a*(c^{-1})*(b^{-1})*c=\\
& a*(c^{-1})*d*(b^{-1})*(d^{-1})*(a^{-1})*c*(b^{-1})=\\ 
&(a*d*(a^{-1})*(d^{-1}))^2= \\
&a*(d^{-1} )*c*(e^{-1})*(d^{-1} )*e*(a^{-1})*c=\\
& b*c*b*( c^{-1} )*(b^{-1} )*c*b*( c^{-1} )=\\ 
&(b*(c^{-1})*(d^{-1})*(c^{-1}))^2=\\
& b*d*(c^{-1})*(e^{-1})*(b^{-1})*c*e*(d^{-1})=\\ 
&b*d^2*c*d^2*(b^{-1})*c=(b*(d^{-3}) )^2=\\
& b*e*c^2*(e^{-1} )*(b^{-1}) *( c^{-2} )=( c^2*(d^{-2}) )^2=\\ 
&(c*d*(c^{-1})*(d^{-1}) )^2= \\
&a*(c^{-2} )*(d^{-1})*(a^{-1} )*b^2*d*c^2=1.
\end{align*}
 
This group has 244 classes with the following order structure:\\

\begin{tabular}{|c|c|c|}\hline
order of element& number of elements& number of classes\\ \hline  
2&	1343& 51\\
3&	512& 2\\
4&	11968& 104\\
6&	24064& 62\\
8&	3072& 8\\
12&	8192 &16\\ \hline
\end{tabular}
\newpage

Several attempts to find the automorphism group of this group have failed. Running times were over 2 days without completion.\\

\subsection{\underline{Group \#181}} The automorphism tower for this group starts out with the following group of order 1536:
\begin{align*}
&b^4=c^4=d^2=b*d*(b^{-1})*d=(c*d)^2=\\
&a^3*d*(a^{-1})*d=a^2*b*a^2*(b^{-1})=\\ 
&(a*(c^{-1}))^3=b*c^2*b*c^2=\\
& a*b*c*b*( c^{-1}) * (a^{-1} ) * ( b^{-1})=\\ 
&a*c*a*c*b^2*a*(c^{-1})=1.
\end{align*}
This group has 33 classes and the order structure:\\

\begin{tabular}{|c|c|c|}\hline
order of element& number of elements& number of classes\\ \hline   
2& 159&10\\	   
3& 128&1\\   
4& 480&13\\	   
6& 384&3\\   
8& 384&5\\ \hline	 
\end{tabular}
\linebreak\\
An alternate presentation for this order 1536 group is:
\begin{align*}
& b^4=c^4=d^4=b*d*b*(d^{-1})=(c*(d^{-1}))^2=\\
& a^4*(d^{-2})=a^2*b*a^2*(b^{-1})=a*b*d*a*d*(b^{-1})=\\
& a*b*(d^{-1})*a*(d^{-1})*(b^{-1})=\\
&a*(b^{-1})*d*a*(b^{-1})*(d^{-1})=\\
& (a*(c^{-1}))^3=b^2*d*(c^{-1})*(d^{-1})*(c^{-1})=\\
& b*c^2*b*(c^{-2})=a*b*c*b*(c^{-1})*(a^{-1})*(b^{-1})=\\
& a*c*a*c*(b^{-2})*a*(c^{-1})=1.\\
\end{align*}
The second member of this tower is a group of order 6144 = $2^{11} * 3$. 
\begin{align*}
&a^2=b^6=c^4=d^2=a*b*a*(b^{-1} )=(a*c)^2=\\
& a*c^2*d*a*c^{-2}*d=(a*d)^4=(a*d*(c^{-1})*d)^2= \\
&b^3*c*(b^{-3})*(c^{-1})=(b*c*b*(c^{-1}) )^2=\\ 
&b*c*d*(c^{-1})*b*(c^{-1})*d*(c^{-1})=\\
& c*d*c*d*(c^{-1})*d*(c^{-1})*d=\\ 
&a*b*(c^{-1})*(b^{-1})*d*(b^{-1})*(c^{-1})*a*b*d=\\
&b*( c^{-1} ) *d* (b^{-1} ) *d* ( b^{-1} ) *c*d*b*d=1.
\end{align*}
This group has 72 classes and the order structure:\\

\begin{tabular}{|c|c|c|}\hline
order of element& number of elements& number of classes\\ \hline  
2&	431	&17\\
3&	128	&1\\
4&	1744	&28\\
6&	1408	&7\\
8&	1920	&16\\
12&	512	&2\\ \hline
\end{tabular}
\linebreak\\

The final, complete group, in this tower is a group of order $12,288 = 2^{12} * 3$ and has the presentation:\\
\begin{align*}
&b^4=c^4=a^2*d^2=a*d*(a^{-1})*(d^{-1} )=\\
&a^2*(d^{-4})=b^2*c*b^2*c=b*c^2*(b^{-1} )*c^2=\\
&a^3*c*a^{-3}*c^{-1}=(a*b*a*b^{-1})^2=\\
&(a*b*d*c^{-1})^2=(a*b*d^{-1}*b^{-1})^2=\\
&a*(b^{-1})*d*c*a*(b^{-1})*d*(c^{-1})=(a*c*a*(c^{-1}))^2=\\
&a*c*d^{-1}*c^{-1}*a*c^{-1}*d^{-1}*c=\\
&a*d*b*c^2* (d^{-1} )*(a^{-1} ) *b=(c*d)^4=\\
&(c*(d^{-1}) )^4=a^2*c*a*b*(a^{-1})*(c^{-1})*a*(b^{-1})=\\
&a*b*c*b*( c^{-1} )*(a^{-1} )*b*(a^{-1} ) *b=1.
\end{align*}
This group has 78 classes and an order structure:\\

 \begin{tabular}{|c|c|c|}\hline
order of element& number of elements& number of classes\\ \hline  
2&	559 &18\\	   
3&	128 &1\\	   
4&	4560 &34\\	   
6&	2432 &10\\	   
8&	3072 &11\\	   
12&	1536 &3\\ \hline
\end{tabular} 
\linebreak\\

\subsection{\underline{Group \#183}} The automorphism tower for this group has three steps with the following presentations and order structure:\\

a. Group of order 9216 = $2^{10} * 3^2$:\\
\begin{align*}
&a^2=b^{12}=c^2=(a*b*a*(b^{-1}))^2=a*(b^{-2})*c*a*c*b^2=\\
&(a*c)^4=a*b^2*a*c*a*c*(b^{-2})=(a*b*a*(b^{-1})*c)^2=\\
& a*b*c*a*b*c*a*c*(b^{-1})*c*(b^{-1})=\\
& a*b*a*b*c*b^4 * c*b^2=(b^2*c)^4=(b*c)^6=\\
& a*b*c*( b^{-1} )*a*c*b*c*(b^{-1} )*c*b*c* (b^{-1} )*c=1.
\end{align*}  
This group has 50 classes and the order structure:\\

\begin{tabular}{|c|c|c|}\hline
order of element& number of elements& number of classes\\ \hline   
2&	495 &18\\	   
3&	800 &3\\	   
4&	3600 &18\\	   
6&	3168 &7\\   
12&	1152 &3\\ \hline
\end{tabular}	 
\linebreak\\

b. The next group has order 18,432 and presentation:
\begin{align*}
&a^2=c^2=d^2=e^2=(a*c)^2=(a*e)^2=b*c*(b^{-1} )*c=\\
&(c*e)^2=a*b*a*e*b*e=b^6*c=b*d*c*e*d*( b^{-1} )*e=\\ 
&(a*b*a*(b^{-1}))^2=a*(b^{-2})*d*a*d*b^2=(a*d)^4=\\
&b^2*e*d*e*(b^{-2})*d=(c*d)^4=a*b^2*a*d*a*d*(b^{-2})=\\
&(b*d^{-1})^6=a*b*d*(b^{-1} )*a*d*b*d* (b^{-1} )*d*b*d*(b^{-1} )*d=1.
\end{align*}
This group has 88 classes and the order structure:\\

\begin{tabular}{|c|c|c|}\hline
order of element& number of elements& number of classes\\ \hline  
2&	943	&26\\
3&	800	&3\\
4&	7248	&41\\
6&	6752	&12\\
12&	2688	&5\\ \hline
\end{tabular}
\linebreak\\

c. The final group is complete of order 36,864 and has a presentation:
\begin{align*}
&a^2=c^2=d^2=e^2=f^4=(a*c)^2=(a*e)^2=\\
&a*f*a*( f^{-1} )=b*c*(b^{-1} )*c=( c*e)^2=\\
&a*b*a*e*b*e=c*d*f*c* ( f^{-1})*d=\\ 
&c*d*(f^{-1})*c*f*d=b*f*(b^{-1})*e*d*b*(f^{-1})=\\
&(a*b*a*(b^{-1}) )^2=a*(b^{-2} )*d*a*d*b^2=\\ 
&a*c*f*(b^{-1})*f^2*b*f=a*c*f*e*f^2*e*f=\\ 
&(a*d)^4=b^2*e*d*e*(b^{-2})*d=(b*e*(b^{-1})*d)^2=\\
&(e*f*e*( f^{-1}) )^2=a*b*a*d*c*( f^{-1} )*b*d*f=\\ 
&a*b^2*a*d*a*d*(b^{-2})=a*c*d*a*e*f*d*e*f=1.
\end{align*}
\newpage
This complete group has 98 classes and the order structure:\\

\begin{tabular}{|c|c|c|}\hline
order of element& number of elements& number of classes\\ \hline  
2&1167&23\\
3&800&3\\
4&12144&41\\
6&13152&19\\
8&3072&5\\
12&4992&5\\
24&1536&1\\ \hline
\end{tabular}

\begin{thebibliography}{99}
\bibitem{12} W. Becker, \textquotedblleft Automorphism groups for the groups of the form $C_p\times C_p$ @ $X$, where $X$ is a group of order 32 and has $D_4$ or $Q_2$ action on the $p$-group".\\
\end{thebibliography}
\end{document}